\documentclass[11pt]{amsart}
\usepackage{dsfont}
\usepackage{times}
\usepackage{enumerate}
\usepackage[euler-digits]{eulervm}
\theoremstyle{plain}
\usepackage{graphicx,color} 
\usepackage{mathtools}

\usepackage{hyperref}
\usepackage{subcaption}
\usepackage{amsmath, amssymb, graphics}
\usepackage{float}

\newtheorem{theorem}{Theorem}[section]
\newtheorem*{theorem*}{Theorem}
\newtheorem{cor}{Corollary}[section]
\newtheorem{prop}{Proposition}[section]
\newtheorem{lemma}{Lemma}[section]

\newtheorem{remark}{Remark}[section]

\newtheorem*{conjecture*}{Conjecture}

\begin{document}

\title{The Euler-Helfrich Functional}
\author{Bennett Palmer and \'Alvaro P\'ampano}
\date{\today}
\maketitle

\begin{abstract}
We investigate equilibrium configurations for surface energies which contain the squared $L^2$ norm of the difference of the mean curvature $H$ and the spontaneous curvature $c_o$ coupled with the elastic energy of the boundary curve, which we studied previously in \cite{PP2}.
 
It is shown that if a critical surface for this type of functional is  axially symmetric, then it satisfies a simpler second order variational problem. Many examples of solutions of this are given.
\\

\noindent K{\tiny EY} W{\tiny ORDS}.\, Helfrich Energy, Bending Energy, Energy Minimization.

\noindent MSC C{\tiny LASSIFICATION} (2020).\, 49Q10, 53A04, 53A05.
\end{abstract}

\begin{figure}[H]
\includegraphics[width=0.25\linewidth,angle=0]{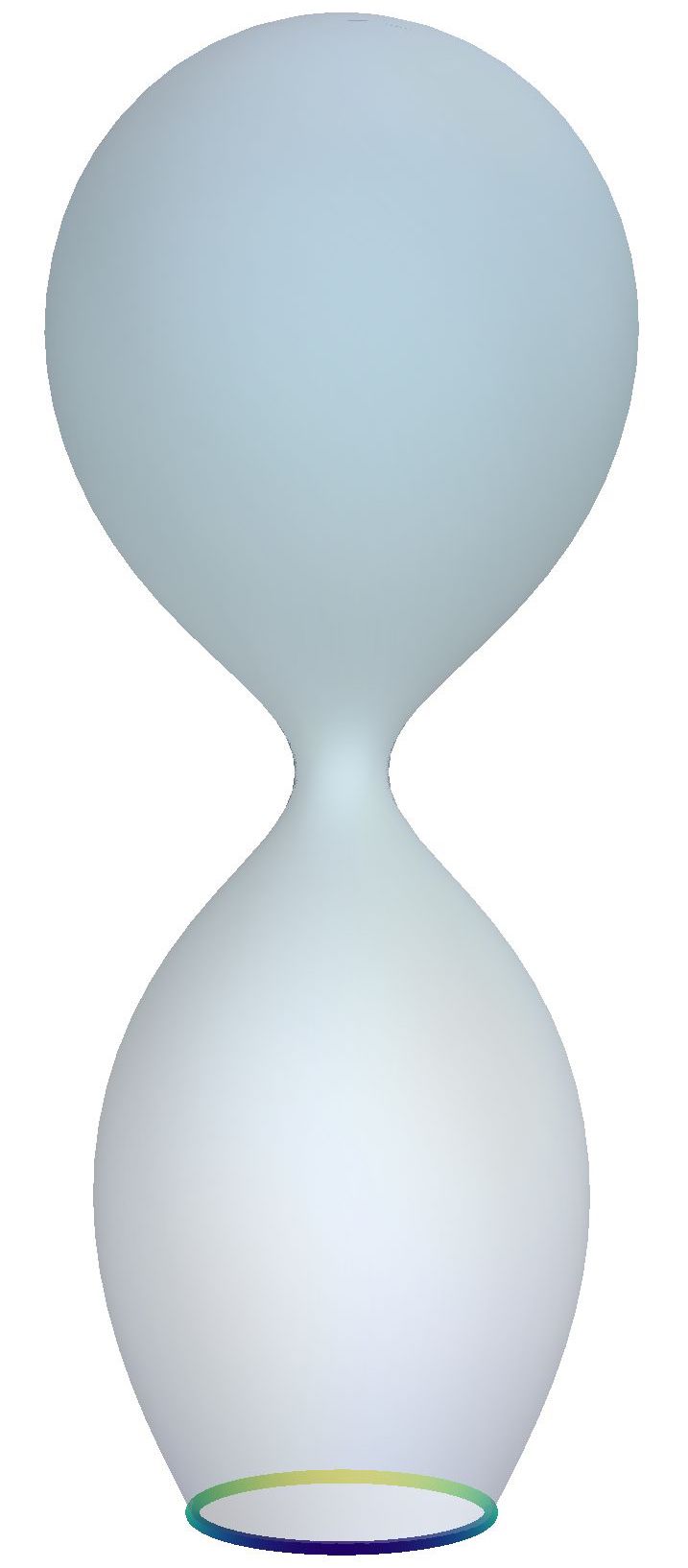}
\end{figure}

\section{Introduction}

We study a variational model for a bilipid membrane having elastic boundary components. The surface energy carries the widely studied elastic energy introduced by Helfrich (\cite{Helfrich}), while the boundary curves carry the classical bending energy studied by Euler, \cite{Euler} (see \cite{Levien} for a historical background), hence we refer to the total energy by the name \emph{Euler-Helfrich}. 

The main constituents of cellular membrane in most living organisms are lipid bilayers formed from a double layer of phospholipids, which have, in general, distinct compositions and tensions.  Lipid bilayers are frequently modeled as mathematical surfaces since they are very thin compared to the size of the observed cells or vesicles.  The difference in the layers' tensions forces the membrane to bend toward the interface having the lower tension \cite{RL}, making it  energetically favorable for the surface to attain a specific value of its mean curvature. This tendency of the surface to favor a specific value of the mean curvature is encoded in the Lagrangian by using the spontaneous curvature  $c_o$   (our definition for $c_o$ differs from the classical spontaneous curvature by the sign and a coefficient two).

When suspended in an aqueous solution, phospholipid membranes spontaneously aggregate in order to shield their hydrophobic tails from the solvent, resulting in a closed bilayer. However, the formation of pores in biological membranes is quite common due to different types of stimuli, resulting in stabilized   open bilipid membranes with edges,  \cite{BMF}. These have been modeled with the inclusion of a variety of types of boundary energies \cite{MF}, \cite{TOY}, \cite{TOY2} and \cite{W}. Other studies of elastic surfaces having elastic boundaries can be found in  \cite{AB}, \cite{BMF}, \cite{G17} and \cite{PP2}. In addition, investigations of elastic surfaces with inelastic boundaries can be found in \cite{BR}, \cite{C}, \cite{CZT}, \cite{Tu2}, \cite{TOY03} and \cite{Z}.

In this paper, we will discuss critica for the Euler-Helfrich energy among compact surfaces with boundary. We are primarily interested in those axially symmetric critical surfaces having the topological type of a disc. In our previous paper \cite{PP2}, we studied the case of topological annuli and established that, in many cases, the absolute minimizer has constant mean curvature. However, in the same paper we showed that, for many choices of parameters in the energy functional, constant mean curvature surfaces of disc type cannot occur. This  complicates the determination of extrema.

The Helfrich energy is quadratic in the second derivative and its critical surfaces satisfy a fourth order elliptic partial differential equation (see \eqref{EL1} below). However, the topological restriction, together with the symmetry, imply that our surfaces are critical for a Lagrangian quadratic in the first derivatives of the surface. Its critica are either constant mean curvature surfaces or they are closely related to singular minimal surfaces, a class of surfaces that model soap films under a constant gravitational force, \cite{Dierkes}. 

In order to clarify this, we study a generalization of the conformal Gauss map which is used to study Willmore surfaces. Although the modified map is no longer conformal, it is still critical for a variational principle that defines a harmonic map with potential.

Axially symmetric surfaces which are critical for the Helfrich energy have already been studied by many authors, for example \cite{DDG}, \cite{MS} and \cite{CZT}, but with boundary conditions distinct from those considered here. The reference \cite{CZT} shows that axially symmetric critical surfaces for the Helfrich energy satisfy a third order differential equation which is equivalent to equation \eqref{fir} below. In the case where the surface is an axially symmetric topological disc, we take this one step further and show that either the surface has constant mean curvature or it must satisfy a second order differential equation which expresses the mean curvature $H$ as a function of the coordinates of the the surface $X$ and its normal $\nu$:
$$H+c_o=-\frac{\nu\cdot V}{X\cdot V}\,,$$
for a constant vector $V\in{\bf R}^3$. After a rigid motion, this equation is precisely \eqref{nonCMC} below.  

Exploiting this allows us to compute many examples some of which can be viewed, after a conformal transformation, as sessile or pendant drops in hyperbolic space (c.f. \cite{Rafa2}). We calculate their energies numerically, which, in some cases, gives convincing evidence that they are minimizers.

\section{The Variational Problem}

Let $\Sigma$ be a compact, connected, oriented surface with smooth  boundary, and consider an immersion of $\Sigma $ in an Euclidean 3-space ${\bf R}^3$, 
$$X:\Sigma\rightarrow{\bf R}^3\,.$$ 
Throughout this paper, we assume $X$ to be of class $\mathcal{C}^4$ up to the boundary. Let $\nu$ denote the unit normal vector field along $\Sigma$ with the orientation such that $\nu$  points out of any convex domain. We always consider the boundary with its  positive orientation relative to the surface. Motivated by the possible applications of the theory we will always restrict our examples to be embedded.

For a sufficiently smooth curve $C:I\rightarrow{\bf R}^3$, we denote by $s\in I=[0,\mathcal{L}]$ the arc length parameter of $C$, while $\mathcal{L}$ stands for its length. Then, if $\left(\,\right)'$ represents the derivative with respect to the arc length, the vector field $T(s):=C'(s)$ is the unit tangent to $C$. Moreover, the (Frenet) curvature of $C$ is defined by $\kappa(s):=\lVert T'(s)\rVert\geq 0$. The Darboux frame along $C$ is the orthonormal frame $\{n,T,\nu\}$, consisting of the unit tangent of the boundary $T$, the normal $\nu$ to the surface and the (outward pointing) conormal $n:=T\times \nu$ of the boundary. The derivative of this frame with respect to the arc length parameter $s$ is given by
\begin{equation*}
\begin{pmatrix}
n\\T\\\nu \end{pmatrix}'=
 \begin{pmatrix}
0 &-\kappa_g& \tau_g\\
\kappa_g&0&\kappa_n\\
-\tau_g&-\kappa_n&0
\end{pmatrix}\begin{pmatrix}
n\\T\\\nu \end{pmatrix},
\end{equation*}
where the functions involved, $\kappa_g$, $\kappa_n$ and $\tau_g$ are, respectively, the geodesic curvature, the normal curvature and the geodesic torsion of the boundary relative to the immersion $X$. From the definitions above, we get $\kappa^2(s)=\lVert T'(s)\rVert^2=\kappa_g^2(s)+\kappa_n^2(s)$. \emph{Our sign convention for $\kappa_g$ is opposite from the usual one}.

For an immersion $X:\Sigma\rightarrow{\bf R}^3$, the \emph{Euler-Helfrich energy} functional is the total energy
\begin{equation*}
E_{a,c_o,b,\alpha,\beta}=E:=\int_\Sigma\left(a\left[H+c_o\right]^2+b K\right)d\Sigma+\oint_{\partial\Sigma}\left(\alpha\kappa^2+\beta\right)ds\,,
\end{equation*}
where $a>0$, $\alpha>0$, $\beta>0$ and $c_o$, $b$ are any real constants. Here, $H$ denotes the mean curvature of the surface and $K$ is its Gaussian curvature. The coefficient $a$ is the bending rigidity of the surface, $c_o$ represents the spontaneous curvature and $b$ is the saddle-splay modulus. \emph{Note that $c_o$ differs by the sign and a coefficient two from the classical spontaneous curvature}. The constant $\alpha$ is called the flexural rigidity (\cite{L}), while $\beta \mathcal{L}$ is the line tension. Under rescaling $X\mapsto tX$, $t>0$, these coefficients change according to $a\mapsto a$, $c_o\mapsto c_o/t$, $b\mapsto b$, $\alpha \mapsto t\alpha$ and $\beta\mapsto \beta/t$.

For simplicity, we assume that all connected components of the boundary are made of the same material, so that the flexural rigidity $\alpha$ and the constant $\beta$ are the same constants for all boundary components. 

We will next review the first variation formula for the energy $E$. For details, the reader is referred to Section 2 of \cite{PP2} (see also Appendix A of the same reference). Unlike most variational problems for surfaces with boundary, there is no restriction on an admissible variation field $\delta X$ other than it's regularity and that regards we assume $\delta X$ to be of class ${\mathcal C}^4$. We associate to it a vector field tangent to the surface by
\begin{equation}\label{J}
\mathcal{J}_{\delta X}:=\left(H+c_o\right)\nabla\left[\nu\cdot\delta X\right]-\left(\nu\cdot\delta X\right)\nabla\left[H+c_o\right]+\left(H+c_o\right)^2\delta X^T\,,
\end{equation}
and we define a vector field $J$, defined along the boundary, by
$$J:=2\alpha T''+\left(3\alpha\kappa^2-\beta\right)T\:.$$
The field $J$ is the Noether current associated to translational invariance of the boundary energy. Then the first variation of the energy $E$ can be expressed
\begin{eqnarray}\label{var1}
\delta_{\delta X}E_{a,c_o, b,\alpha, \beta}&=& a\int_\Sigma \left(\Delta H+2\left[H+c_o\right]\left(H\left[H-c_o\right]-K\right)\right)\nu\cdot\delta X\:d\Sigma\nonumber\\
&&+\oint_{\partial \Sigma} \left(\left[a\mathcal{J}_{\delta X}+bK\delta X\right]\cdot n+J'\cdot \delta X\right)ds\nonumber\\
&&+\oint_{\partial\Sigma}b\left(\left[d\nu+2H\,{\rm Id}\right]\nabla\left[\nu\cdot\delta X\right]\right)\cdot n\,ds\:.   
\end{eqnarray}

This results in the following necessary and sufficient conditions for criticality of the surface:
\begin{eqnarray}
\Delta H+2\left(H+c_o\right)\left(H\left[H-c_o\right]-K\right)&=&0\,,\quad\quad\quad\text{on $\Sigma$}\,,\label{EL1}\\
a\left(H+c_o\right)+b\kappa_n&=&0\,,\quad\quad\quad\text{on $\partial\Sigma$}\,,\label{EL2}\\
J'\cdot\nu-a\partial_n H+b\tau_g'&=&0\,,\quad\quad\quad\text{on $\partial\Sigma$}\,,\label{EL3}\\
J'\cdot n+a\left(H+c_o\right)^2+bK&=&0\,,\quad\quad\quad\text{on $\partial\Sigma$}\,.\label{EL4}
\end{eqnarray}
The partial differential equation \eqref{EL1} arises from considering compactly supported variations, equation \eqref{EL2} results from normal variations which fix the boundary, equation \eqref{EL3} from normal variations that fix the tangent plane along the boundary and \eqref{EL4} results from variations tangent to the surface.

Using variation by rescaling, $X\mapsto tX$ and differentiating the energy with respect to $t$ at $t=1$, we obtain a necessary condition for equilibrium
\begin{equation}\label{rescalings}
2a c_o\int_\Sigma\left(H+c_o\right)d\Sigma=\oint_{\partial\Sigma}\left(\alpha\kappa^2-\beta\right)ds\, .
\end{equation}

\section{The Governing PDE for the Helfrich Energy}

For an immersion $X:\Sigma\rightarrow{\bf R}^3$ critical for $E$, the governing equation in the interior of the surface $\Sigma$ is the fourth order elliptic partial differential equation \eqref{EL1}, which we write as $L[H+c_o]=0$, where $L$ is the second order elliptic operator
\begin{equation}\label{L}
L[f]:=\Delta f+2\left(H\left[H-c_o\right]-K\right)f\,.
\end{equation}
By expressing the equation $L[H+c_o]=0$ in a nonparametric form we get a fourth order elliptic partial differential equation with coefficients that depend analytically on the height function and its derivatives. The real analyticity of any $\mathcal{C}^4$ critical surface then follows from Theorem 6.6.1 of \cite{Morrey}.

We will call a smooth immersion $X:\Sigma\rightarrow{\bf R}^3$ a critical surface for the Helfrich energy 
$$\mathcal{H}[\Sigma]:=\int_\Sigma\left(H+c_o\right)^2d\Sigma\,,$$
if it is critical with respect to compactly supported variations, i.e. if \eqref{EL1} holds in the surface's interior. Then for an arbitrary variation of the immersion $X:\Sigma\rightarrow{\bf R}^3$, i.e. $X+\epsilon\delta X+\mathcal{O}(\epsilon^2)$, the first variation formula for $\mathcal{H}$ reads (c.f. \eqref{var1})
$$\delta\mathcal{H}[\Sigma]=\oint_{\partial\Sigma}\mathcal{J}_{\delta X}\cdot n\: ds\,.$$

We take a variation field $\delta X\equiv E_i:=\nabla X_i+\nu_i\nu$, $i=1,2,3$, where $E_i$ are the constant unit vector fields in the direction of the coordinates axes. Since $\delta X\equiv E_i$ is a constant vector field which generates a symmetry (in this case, a translation) of the Lagrangian, it must be the case that $\delta\mathcal{H}[\Omega]=0$ for any subdomain $\Omega\subset\Sigma$ and, therefore
\begin{equation}\label{divprev}
0=\oint_{\partial\Omega}\mathcal{J}_{E_i} \cdot n\,ds\,,
\end{equation}
which implies, by Stokes' Theorem, that
$$\nabla\cdot \mathcal{J}_{E_i}=\nabla\cdot\left(\left[H+c_o\right]\nabla\nu_i-\nu_i\nabla\left[H+c_o\right]+\left[H+c_o\right]^2\nabla X_i\right)=0\,,$$
holds on $\Omega$, for $i=1,2,3$. Wherever $H+c_o\neq 0$, this equation can be expressed as
$$\nabla\cdot\left(\left[H+c_o\right]^2\nabla\left[\frac{\nu_i}{H+c_o}+X_i\right]\right)=0\,.$$

Next, from the definition \eqref{L} and the product rule for the Laplacian, it follows that 
$$L[fg]=gL[f]+\frac{\nabla\cdot\left(f^2\nabla g\right)}{f}\,,$$
for smooth functions $f$ and $g$, wherever $f\neq 0$. We apply this with $f=H+c_o$ and $g=\nu_i/(H+c_o)+X_i$ to obtain
\begin{equation}\label{L1}
L[(H+c_o)X_i+\nu_i]=0\,,
\end{equation}
for $i=1,2,3$. 

We take now a variation field $\delta X\equiv X=q\nu+\nabla X^2/2$, where $q:=X\cdot \nu$ is the support function. This variation field generates a one parameter family of rescalings. Using again the first variation formula for the Helfrich energy $\mathcal{H}$, we get for the subdomain $\Omega\subset\Sigma$,
\begin{equation*}
2c_o\int_\Omega \left(H+c_o\right)d\Sigma=\oint_{\partial\Omega} \mathcal{J}_X \cdot n\,ds\,.
\end{equation*}
As before, since $\Omega$ is arbitrary, we get after applying Stokes' Theorem,
$$\nabla\cdot  \mathcal{J}_X= \nabla\cdot\left(\left[H+c_o\right]\nabla q-q\nabla\left[H+c_o\right]+\left[H+c_o\right]^2\nabla X^2/2\right)=2c_o\left(H+c_o\right),$$
which can be rewritten, wherever $H+c_o\neq 0$, as
$$\nabla\cdot\left(\left[H+c_o\right]^2\nabla\left[\frac{q}{H+c_o}+\frac{X^2}{2}\right]\right)=2c_o\left(H+c_o\right).$$
Then, a similar argument as before shows that
\begin{equation}\label{L2}
L\left[\left(H+c_o\right)\frac{X^2}{2}+q\right]=2c_o\,,
\end{equation}
holds on $\Omega$.

We organize the formulas given above in a manner similar to that used in the conformal geometry of surfaces.  We denote by $\left({\bf E}^5_1,\langle\,,\,\rangle\right)$ the five-dimensional Minkowski space with signature $(+,+,+,+,-)$ and consider the hyperquadric ${\bf S}_1^4:=\{Y\in{\bf E}^5_1\,\lvert\,\langle Y,Y\rangle=1\}$, which endowed with the induced metric, is a model for the four-dimensional de Sitter space. The conformal Gauss map $Y:\Sigma\rightarrow {\bf S}_1^4\subset{\bf E}_1^5$ is defined by
$$Y:=H{\underline X}+\left(\nu,q,q\right),$$
where
$${\underline X}:=\left(X,\frac{X^2-1}{2},\frac{X^2+1}{2}\right).$$
The Euclidean space ${\bf R}^3$ can be isometrically embedded in the light cone of ${\bf E}^5_1$ by means of the map $X\mapsto{\underline X}$. It is well known that \emph{an immersion $X$ is Willmore if and only if its conformal Gauss map $Y$ is a zero mean curvature map}. The map $Y$ is spacelike except at the umbilic points of $X$. For more details about the conformal Gauss map of Willmore immersions, see \cite{Palmer}. Here, we extend this notion to study immersions critical for the Helfrich energy $\mathcal{H}$. 

For a real constant $c_o$, we modify the conformal Gauss map $Y$ and define
$$Y^{c_o}:=\left(H+c_o\right){\underline X}+\left(\nu,q,q\right).$$
Note that $Y^{c_o}=Y+c_o{\underline X}$ may be regarded as a parallel surface to the conformal Gauss map, since ${\underline X}$ is a section of the normal bundle of $Y$.

Then, we have the following characterization of critical immersions for the Helfrich energy $\mathcal{H}$.

\begin{theorem}\label{confGmthe} The immersion $X:\Sigma\rightarrow{\bf R}^3$ is critical for the Helfrich energy $\mathcal{H}$ with respect to compactly supported variations if and only if the modified conformal Gauss map $Y^{c_o}$ satisfies
\begin{equation*}
\Delta Y^{c_o}+\lVert dY^{c_o}\rVert^2 Y^{c_o}=2c_o l^T,
\end{equation*}
where $l^T$ denotes the projection of the constant light like vector  $l:=(0,0,0,1,1)$ to ${\bf S}_1^4$. (The metric on $\Sigma$ is the one induced by the immersion $X$.)

Consequently, the immersion $X:\Sigma\rightarrow{\bf R}^3$ is critical for the Helfrich energy $\mathcal{H}$ if and only if the map $Y^{c_o}$ is critical for the energy
$$\mathcal{F}[Z]:=\int_\Sigma\left(\lVert dZ\rVert^2+4c_oU(Z)\right)d\Sigma\:,$$
where $U(Z):=Z_4-Z_5$, with respect to compactly supported variations.
\end{theorem}
{\it Proof.\:} Assume first that the immersion $X:\Sigma\rightarrow{\bf R}^3$ is critical for the Helfrich energy $\mathcal{H}$ for compactly supported variations. Then, as proved above, equations \eqref{L1} and \eqref{L2} hold. Combining them with the definition of the map $Y^{c_o}$ we conclude that $Y^{c_o}$ satisfies the equation
\begin{equation}\label{L3}
L[Y^{c_o}]=2c_o l\,.
\end{equation}

In what follows, we compute $\lVert dY^{c_o}\lVert^2=4e^{-\rho}\langle Y^{c_o}_z,Y^{c_o}_{\bar{z}}\rangle$, where $z$ is a local complex coordinate for which the metric induced by $X$ is $ds_X^2:=e^\rho\lvert dz\rvert^2$. Denoting by $\Phi dz^2$ the Hopf differential, we get that
$$Y^{c_o}_z=\left(H+c_o\right)_z{\underline X}+c_o{\underline X}_z-\Phi e^{-\rho} {\underline X}_{\bar{z}}.$$
From this and the Gauss equation $\lvert \Phi\rvert^2=(H^2-K)e^{2\rho}$, it follows that
$$\lVert dY^{c_o}\rVert^2=2\left(H^2-K+c_o^2\right).$$
The metric induced by $Y^{c_o}$ is
$$\langle dY^{c_o}, \,dY^{c_o}\rangle=2\Re\left(-c_o\Phi dz^2+\left[c^2_o+\lvert\Phi\rvert^2 e^{-2\rho}\right]\frac{e^\rho}{2}dz d{\bar z}\right).$$
This metric has the area form
$$d{\mathcal A}=\frac{i}{2}\, \lvert  c_o^2-\lvert\Phi\rvert^2 e^{-2\rho} \rvert \,e^\rho dz\wedge d{\bar z}\:,$$
so $Y^{c_o}$ fails to be a definite immersion at the points where $H^2-K=c_o^2$ holds. These are precisely the points where one of the principal curvatures of $\Sigma$ satisfies $k_i=H+c_o$. (With regard to this, there is a particularly interesting family of surfaces, the \emph{circular biconcave discoids}, \cite{NOOY}, for which $H^2-K\equiv c_o^2$ holds throughout the surface.)

Finally, since $\langle Y^{c_o},l\rangle=-\left(H+c_o\right)$, we get from \eqref{L3},
\begin{eqnarray}\label{Y}
\Delta Y^{c_o}+\lVert dY^{c_o}\rVert^2Y^{c_o}&=&L[Y^{c_o}]+2c_o\left(H+c_o\right)Y^{c_o}\\&=&2c_o\left(l-\langle Y^{c_o},l\rangle Y^{c_o}\right)=2c_o l^T \nonumber\,.
\end{eqnarray}

This proves one direction. To prove the other direction, just note that if \eqref{Y} holds, then, in particular,
\begin{equation*}\label{tension}
\Delta \langle Y^{c_o},l\rangle+2\left(H^2-K+c_o^2\right)\langle Y^{c_o},l\rangle=2c_o\langle l^T,l\rangle=-2c_o\left(H+c_o\right)^2,
\end{equation*}
which reduces to $L[H+c_o]=0$, proving the statement. {\bf q.e.d.}
\\

\begin{remark}\normalfont Since the left hand side of \eqref{Y} is the tension field of the map $Y^{c_o}$, Theorem \ref{confGmthe} says that the immersion $X:\Sigma\rightarrow{\bf R}^3$ is critical for the Helfrich energy $\mathcal{H}$ with respect to compactly supported variations if and only if the tension field $Y^{c_o}$ is a specific timelike conformal field on ${\bf S}^4_1$.
\end{remark}

Using this characterization of critical immersions $X:\Sigma\rightarrow{\bf R}^3$ for the Helfrich energy $\mathcal{H}$, we will consider some special solutions. These solutions will depend on the causal character of the hyperplane where the associated map $Y^{c_o}$ lies:
\begin{enumerate}[(i)]
\item Take the null vector $l:=(0,0,0,1,1)$ and consider the case in which $Y^{c_o}$ lies in the hyperplane $\langle Y^{c_o},l\rangle=0$. In this case, we have
$$\hspace{1cm}0=\langle Y^{c_o},l\rangle=\left(H+c_o\right)\langle {\underline X},l\rangle=-\left(H+c_o\right),$$
so, we obtain $H+c_o\equiv 0$. This case was studied in detail in \cite{PP2}.
\item Suppose that $Y^{c_o}$ lies in a hyperplane $\langle Y^{c_o},w\rangle=0$, where $w$ is a spacelike vector. After a rigid motion, we may assume that $w=(0,0,1,0,0)$. Then, we have
$$\hspace{1cm}0=\langle Y^{c_o},w\rangle=\left(H+c_o\right)z+\nu_3\,.$$
After a vertical translation, this condition is satisfied by axially symmetric discs. (See Section \ref{4}.)
\item Finally,  consider any timelike vector $w=(0,0,0,0,A)$ for some real constant $A\neq0$ and assume that $Y^{c_o}$ lies in the hyperplane $\langle Y^{c_o},w\rangle=0$. Then, we have
\begin{eqnarray*}
\hspace{1cm}0&=&\langle Y^{c_o},w\rangle=\left(H+c_o\right)\langle{\underline X},w\rangle-A q\\&=&-A\left(H+c_o\right)\frac{X^2+1}{2}-A q.
\end{eqnarray*}
We then apply the linear operator $L$, \eqref{L}, to obtain
\begin{eqnarray*}
\hspace{1cm}0&=&L[0]=AL[-\left(H+c_o\right)\frac{X^2+1}{2}-q]\\&=&-AL\left[\left(H+c_o\right)\frac{X^2}{2}+q\right]-\frac{A}{2}L[H+c_o]=-2Ac_o\,,
\end{eqnarray*}
where, in the last line, we have used that $L[H+c_o]=0$ and \eqref{L2} holds by the criticality of the immersion $X$. Therefore, this case is only possible if the initial immersion $X:\Sigma\rightarrow{\bf R}^3$ is critical for the Willmore energy, i.e. $c_o=0$ in $\mathcal{H}$, which corresponds with minimal surfaces in ${\bf S}^3$.
\end{enumerate}

\subsection{The Flux Homomorphism} 

Let $\Sigma$ be a surface with smooth boundary that is critical for the functional $\mathcal{H}$, i.e. \eqref{EL1} holds on $\Sigma$. In this case there is a well defined {\it flux homomorphism} $\Psi$ from the vector space ${\mathcal K}$ of Killing fields on ${\bf R}^3$ into the first cohomology group $H^1(\Sigma, \bf{R})$ defined for $\mathcal{Z} \in \mathcal{K}$, $[C]\in H_1(\Sigma, \bf R)$ by
\begin{equation}\label{flux}
\mathcal{Z}\mapsto \Psi_\mathcal{Z}  ,\qquad \langle \Psi_\mathcal{Z}  ,[C] \rangle=\oint_C \mathcal{J}_\mathcal{Z} \cdot n\:ds\:,
\end{equation}
where $\mathcal{J}_\mathcal{Z} $ is defined as in \eqref{J}. In this case $\mathcal{J}_\mathcal{Z} $ 
is exactly the Noether current corresponding to $\mathcal{Z} $. That $\Psi_\mathcal{Z} $  is well defined on $H_1$ follows from the first variation formula and the fact that the functional $\mathcal{H}$ is invariant under all isometries.  

\begin{theorem}\label{zeroflux} Let $X:\Sigma \rightarrow {\bf R}^3$ be the immersion of a compact surface with boundary. If the immersion is critical for a functional $E=E_{a,c_o,b,\alpha,\beta}$, then for all $\mathcal{Z} \in {\mathcal K}$, $\Psi_\mathcal{Z} = 0$ holds on $\{[C]\in H_1(\Sigma,{\bf R})\:|\:C\subset \partial \Sigma\}$.
\end{theorem}
{\it Proof.\:} Let $C$ be any cycle contained in $\partial \Sigma$. Using \eqref{EL2}, \eqref{EL3} and \eqref{EL4}, we can replace the factors $H+c_o$, $\partial_nH$ and $(H+c_o)^2$ in \eqref{flux}, to obtain
\begin{eqnarray*}
&&\langle \Psi_\mathcal{Z} ,[C] \rangle=\\
&&\frac{-1}{a}\oint_C\left(b\kappa_n\partial_n[\nu\cdot \mathcal{Z} ]+\left[J'\cdot \nu+b\tau_g'\right]\nu\cdot\mathcal{Z} +\left[J'\cdot n+bK\right]n\cdot\mathcal{Z}\right)ds\,.
\end{eqnarray*}
Expanding out the first term, we get
\begin{eqnarray*}
\langle \Psi_\mathcal{Z} ,[C] \rangle&=&\frac{-1}{a}\oint_C\left( b\kappa_n\left[-(2H-\kappa_n)n\cdot\mathcal{Z} -\tau_g T\cdot\mathcal{Z} +\nu\cdot\partial_n\mathcal{Z} \right]\right.\\
&&\hspace{1cm}\left.+\left[J'\cdot \nu +b\tau_g'\right]\nu\cdot\mathcal{Z}+\left[J'\cdot n +bK\right]n\cdot\mathcal{Z} \right)ds\:.
\end{eqnarray*}
If we integrate the term containing $\tau_g'$ by parts and rearrange terms, we arrive at
\begin{eqnarray*}
\langle \Psi_\mathcal{Z}  ,[C] \rangle&=&\frac{-1}{a}\oint_C\left( b\kappa_n\left[-(2H-\kappa_n)n\cdot\mathcal{Z} -\tau_g T\cdot\mathcal{Z} +\nu\cdot\partial_n\mathcal{Z} \right]\right).\\
&&\hspace{1cm}\left.\left[(J'\cdot n)n+(J'\cdot\nu)\nu\right]\cdot\mathcal{Z} +bKn\cdot\mathcal{Z} -b\tau_g\left[\nu\cdot\mathcal{Z} \right]'\right)ds\\
&=&\frac{-1}{a}\oint_C\left( b\kappa_n\left[-(2H-\kappa_n)n\cdot\mathcal{Z}-\tau_g T\cdot\mathcal{Z}+\nu\cdot\partial_n\mathcal{Z}\right]\right.\\
&&\hspace{1cm}\left[(J'\cdot n)n+(J'\cdot\nu)\nu\right]\cdot\mathcal{Z}+bKn\cdot\mathcal{Z}\\
&&\hspace{1cm}\left. +b\tau_g\left[\kappa_n T\cdot\mathcal{Z}+\tau_g n\cdot\mathcal{Z}-\nu\cdot\partial_s\mathcal{Z}\right]\right)ds\:.
\end{eqnarray*}
The terms containing the factor $n\cdot\mathcal{Z}$ cancel because $K=\kappa_n(2H-\kappa_n)-\tau_g^2$ along $C$, so we are left with
\begin{equation*}\label{xyz}
\langle \Psi_\mathcal{Z} ,[C] \rangle=\frac{-1}{a}\oint_C\left(b\left[\kappa_n \nu\cdot\partial_n \mathcal{Z}-\tau_g\nu\cdot\partial_s \mathcal{Z}\right]+\left[(J'\cdot n)n+(J'\cdot \nu)\nu\right]\cdot \mathcal{Z}\right)ds\:.
\end{equation*}
For the first term above, recall that $\mathcal{K}$ is generated by constant vector fields $E_i\in {\bf R}^3$, which are infinitesimal translations, and fields of the form $E_i\times X$, which are infinitesimal rotations. In the first case, the first term above clearly vanishes. In the second case, the first term in brackets can be written
\begin{eqnarray*}
\kappa_n \nu\cdot E_i\times n -\tau_g\nu\cdot E_i\times T&=&\left(-\kappa_n\nu \times n+\tau_g\nu\times T\right)\cdot E_i\\
&=&\left(-\kappa_nT-\tau_g n\right)\cdot E_i=\nu'\cdot E_i\:,
\end{eqnarray*}
which integrates to zero. So in either case, the first term in above expression of $\langle\Psi_\mathcal{Z},[C]\rangle$ vanishes. For the second term in brackets, note that 
\begin{eqnarray*}
\left[(J'\cdot n)n+(J'\cdot \nu)\nu\right] \cdot \mathcal{Z} =J'\cdot \mathcal{Z}-\left(J'\cdot T\right)T\cdot \mathcal{Z}\:.
\end{eqnarray*}
A straightforward calculation using the definition of $J$ shows that $J'\cdot T \equiv 0$. Then we use that the boundary energy is invariant under isometries and so
$$0=\delta_\mathcal{Z} \left(\oint_C \left[\alpha \kappa^2+\beta\right] ds\right)=\oint_C J'\cdot \mathcal{Z}\:ds\,$$
implying that $\langle \Psi_\mathcal{Z} ,[C] \rangle =0$ and the result follows. {\bf q.e.d.}
\\

\begin{remark}\normalfont Theorem \ref{zeroflux} shows that the boundary energy is related to  the surface energy $\mathcal{H}$ in a natural way. If a surface satisfies \eqref{EL1}, the non vanishing of $\langle \Psi_\mathcal{Z} ,[C] \rangle$ for a boundary cycle represents an obstruction to finding constants $\alpha$ and $\beta$ so that the surface is critical for $E_{a,c_o, b, \alpha, \beta}$.  In Appendix A, we will show the existence of annuli satisfying \eqref{EL1} for which the flux homomorphism does not vanish.
\end{remark}

\section{Axially Symmetric Equilibria}\label{4}

Let  $X:\Sigma\rightarrow{\bf R}^3$ be an axially symmetric immersion. For these immersions, the topology of the surface is restricted and $\Sigma$ must be either a topological disc or a topological annulus. For axially symmetric immersions of topological discs critical for $\mathcal{H}$, we will see that the fourth order Euler-Lagrange equation \eqref{EL1} always reduces to a simpler second order differential equation and the same is true for topological annuli which are critical for a functional $E=E_{a,c_o,b,\alpha,\beta}$. This reduction will lead to the generation of many examples of critica for the boundary value problem, i.e. for the energy $E$.

\subsection{Critica of ${\mathcal H}$}

We begin by obtaining a first integral of \eqref{EL1}. Consider an axially symmetric immersion critical for $\mathcal{H}$, i.e. \eqref{EL1} holds on the interior of $\Sigma$. After a rigid motion, we may assume that the $z$-axis is the axis of rotation. We denote by $C_j$, $j=1,2$ the two horizontal parallels bounding any domain $\Omega$ in $\Sigma$ and take $\delta X\equiv E_3$. Next, arguing as in previous section, we obtain from \eqref{divprev} that for $j=1,2$,
\begin{equation}\label{fir}
r\left(\left[H+c_o\right]\partial_n \nu_3-\nu_3\partial_n H\, +\left[H+c_o\right]^2\partial_n z \right)=(-1)^j\bar{A}\,,
\end{equation}
for some real constant $\bar{A}$ independent of $j$. This first integral is another form of the integral obtained in \cite{CZT}. In particular, if $\bar{A}=0$, this implies that, along every parallel 
either $H+c_o\equiv 0$ or
\begin{eqnarray}
0&\equiv& \left(H+c_o\right)\partial_\varsigma \nu_3- \nu_3 \partial_\varsigma H+\left(H+c_o\right)^2\partial_\varsigma z\label{=0}\\ &=&\left(H+c_o\right)^2\partial_\varsigma \left(\frac{\nu_3}{H+c_o}+z\right),\nonumber
\end{eqnarray}
holds, where $\varsigma$ denotes the arc length parameter along the generating curve. Consequently, either $H+c_o\equiv 0$, or 
\begin{equation}\label{nonCMC}
H+c_o=\frac{\nu_3}{A-z}\,,
\end{equation}
for a real constant $A$. Note that after a suitable vertical translation, if needed, the constant of integration $A$ can be taken to be zero and hence, from now on, we assume $A=0$ holds.

It turns out that, when the surface $\Sigma$ is a topological disc, we can let $C_1$ shrink to a point so that, from \eqref{fir}, $\bar{A}=0$, and we have the following result.

\begin{theorem}\label{discsnew} Let $X:\Sigma\rightarrow{\bf R}^3$ be an axially symmetric immersion of a topological disc critical for $\mathcal{H}$. Then, the critical domain is either a planar disc $(H\equiv -c_o=0)$, a spherical cap $(H\equiv -c_o\neq 0)$ or, after a suitable vertical translation, the mean curvature of the immersion satisfies,
$$H+c_o=-\frac{\nu_3}{z}\,,$$
on $\Sigma$.
\end{theorem}
{\it Proof.\:} From the geometric argument above, shrinking one of the boundary components, the result follows. For completeness, we include here an analytical proof. 

If the surface has constant mean curvature $H=-c_o$, the result is direct. In the rest of the cases, by equations \eqref{EL1} and \eqref{L1}, both the functions $H+c_o\neq 0$ and $(H+c_o)z+\nu_3$ satisfy the equation $L[f]=0$, \eqref{L}, and since the surface is axially symmetric, this reduces to a second order linear ordinary differential equation. If $\varsigma$ denotes the arc length parameter along the generating curve measured from the cut with the vertical axis to the boundary, then in order for both functions to be regular, the $\varsigma$ derivatives must vanish at $\varsigma=0$. This means that there is a linear relation $A(H+c_o)=(H+c_o)z+\nu_3$, which gives us \eqref{nonCMC}. {\bf q.e.d.}
\\

In the following proposition we show that any immersion satisfying \eqref{nonCMC} also satisfies \eqref{EL1}, whether the surface is axially symmetric or not and regardless of the topological type of the surface. Non-axially symmetric surfaces satisfying \eqref{nonCMC} exist, the simplest examples being cylinders over an \emph{elastic curve circular at rest} (defined, for instance, in \cite{AGM}) located in the plane $\{y=0\}$.

\begin{prop} Let $X:\Sigma\rightarrow{\bf R}^3$ be an immersion whose mean curvature satisfies
$$H+c_o=-\frac{\nu_3}{z}\,,$$
on $\Sigma$. Then 
$$\Delta H+2\left(H+c_o\right)\left(H\left[H-c_o\right]-K\right)=0\,,$$ 
also holds on $\Sigma$. In other words, the immersion is critical for the Helfrich energy $\mathcal{H}$ for compactly supported variations.
\end{prop}
{\it Proof.\:} We begin by computing the Laplacian of the mean curvature. Since the surface satisfies \eqref{nonCMC}, applying the product rule for the Laplacian, we have
\begin{equation}\label{Laplacian}
\Delta H=-\Delta\left(\frac{\nu_3}{z}\right)=-\left(\frac{1}{z}\Delta\nu_3+2\nabla\nu_3\cdot\nabla\frac{1}{z}+\nu_3\Delta\frac{1}{z}\right).
\end{equation} 
Next, if we consider a variation of the immersion $X:\Sigma\rightarrow{\bf R}^3$ through translations, its mean curvature remains invariant. Therefore, for a variation vector field $\delta X\equiv E_3$, we get (for details, see Appendix A of \cite{PP2})
$$0=\delta H=\frac{1}{2}\Delta \nu_3+\frac{1}{2}\lVert d\nu\rVert^2\nu_3+\nabla H\cdot \nabla z\,.$$
Using this equation for $\Delta\nu_3$ together with
$$\Delta\frac{1}{z}=\nabla\cdot\left(\nabla\frac{1}{z}\right)=\nabla\cdot\left(-\frac{\nabla z}{z^2}\right)=-\frac{\Delta z}{z^2}+2\frac{\lVert \nabla z\rVert^2}{z^3}\,,$$
in \eqref{Laplacian} we conclude that 
\begin{eqnarray*}
\Delta H&=&\lVert d\nu\rVert^2\frac{\nu_3}{z}+2\frac{\nabla H\cdot \nabla z}{z}+2\frac{\nabla\nu_3\cdot\nabla z}{z^2}+\frac{\nu_3}{z^2}\Delta z-2\frac{\nu_3}{z^3}\lVert \nabla z\rVert^2\\
&=&\lVert d\nu\rVert^2\frac{\nu_3}{z}+2\frac{\nabla H\cdot \nabla z}{z}+2\frac{\nabla\nu_3\cdot\nabla z}{z^2}+2H\frac{\nu_3^2}{z^2}-2\frac{\nu_3}{z^3}\lVert \nabla z\rVert^2\,,
\end{eqnarray*}
where in the last equality we have used that $\Delta z=\Delta X\cdot E_3=2H\nu\cdot E_3=2H\nu_3$.

We then compute, from \eqref{nonCMC}, the gradient of the mean curvature, obtaining
$$\nabla H=-\nabla\left(\frac{\nu_3}{z}\right)=\frac{\nu_3\nabla z}{z^2}-\frac{\nabla\nu_3}{z}\,.$$
Thus,
$$\frac{\nabla H\cdot \nabla z}{z}+\frac{\nabla\nu_3\cdot \nabla z}{z^2}-\frac{\nu_3}{z^3}\lVert\nabla z\rVert^2=0\,,$$
and we conclude, using once again \eqref{nonCMC}, that
\begin{eqnarray*}
\Delta H&=&\lVert d\nu\rVert^2\frac{\nu_3}{z}+2H\frac{\nu_3^2}{z^2}=-\left(H+c_o\right)\lVert d\nu\rVert^2+2H\left(H+c_o\right)^2\\
&=&-2\left(H+c_o\right)\left(H\left[H-c_o\right]-K\right),
\end{eqnarray*}
where in the second line we have used that $\lVert d\nu\rVert^2=4H^2-2K$. This finishes the proof. {\bf q.e.d.}
\\

For $c_o=0$, equation \eqref{nonCMC} characterizes \emph{$(-2)$-singular minimal surfaces} of ${\bf R}^3$, \cite{Dierkes}. In particular, axially symmetric ones were studied in \cite{Rafa}. These surfaces are known to be minimal in the hyperbolic space, i.e. critical for the area functional. The inclusion of an arbitrary constant $c_o$ yields a second term in this energy, in addition to the area.

\begin{theorem} Let $X:\Sigma\rightarrow{\bf R}^3$ be an axially symmetric immersion satisfying $H\neq -c_o$ on $\Sigma$ and assume that the surface is contained in a domain on which $z\neq 0$ holds. Then, the immersion satisfies
$$H+c_o=-\frac{\nu_3}{z}\,,$$
and, hence, it is critical for $\mathcal{H}$ if and only if it is critical for the functional
$$\mathcal{G}[\Sigma]:=\int_\Sigma\frac{1}{z^2}\,d\Sigma+2c_o\int_\Sigma\frac{\nu_3}{z}\,d\Sigma\,,$$
for variations vanishing on $\partial\Sigma$.
\end{theorem}
{\it Proof.\:} We will prove the result by comparing the Euler-Lagrange equation for $\mathcal{G}$ with \eqref{nonCMC}. For convenience, let $G(z):=-1/z$, so the energy $\mathcal{G}$ becomes
$$\mathcal{G}[\Sigma]=\int_\Sigma G_z\,d\Sigma-2c_o\int_\Sigma G\,\nu_3\,d\Sigma=\int_\Sigma\left(G_z-2c_oG\nu_3\right)d\Sigma\,.$$
For a normal variation $\delta X=\psi\nu$, we have $\delta z=\psi \nu_3$ and $\delta \nu_3=-\nabla \psi\cdot \nabla z$. We then get
\begin{eqnarray*}
\delta \mathcal{G}[\Sigma]&=&\int_\Sigma \left(G_{zz}\psi \nu_3 -2H\psi G_z\right)d\Sigma\\
&&-2c_o\int_\Sigma \left(G_z\psi\nu^2_3-\left[\nabla \psi \cdot \nabla z\right] G-2H\psi G\nu_3\right)d\Sigma\\&=&\int_\Sigma \left(G_{zz}\psi \nu_3 -2H\psi G_z\right)d\Sigma-\oint_{\partial \Sigma}\psi G \nabla z\cdot n\:ds\\ 
&&-2c_o\int_\Sigma \left(G_z\psi\nu^2_3+2H\psi \nu_3G+\psi G_z \lVert\nabla z\rVert^2-2H\psi G\nu_3\right)d\Sigma\:,
\end{eqnarray*}
where we have integrated by parts in the second equality. Next, using that $\lVert\nabla z\rVert^2=1-\nu_3^2$, we get
$$\delta \mathcal{G}[\Sigma]=\int_\Sigma \left(G_{zz} \nu_3-2G_z\left[H+c_o\right]\right)\psi\,d\Sigma+2c_o\oint_{\partial \Sigma}\psi G \nabla z\cdot n\:ds\:.$$
If $\psi\equiv 0$ on $\partial \Sigma$, the boundary integral vanishes and we obtain the Euler-Lagrange condition 
$G_{zz} \nu_3-2G_z(H+c_o)\equiv 0$ which implies 
$$H+c_o=-\frac{\nu_3}{z}$$ 
on $\Sigma$, i.e. precisely \eqref{nonCMC}. {\bf q.e.d.}
\\

\begin{remark}\normalfont The energy $\mathcal{G}$ can be expressed as
$$\mathcal{G}[\Sigma]=\widetilde{\mathcal{A}}[\Sigma ]-2c_o\int_{\widetilde{V}} z\:d{\widetilde V}\:,$$
where ${\widetilde{\mathcal{A}}}$ and ${\widetilde V}$ are, respectively, the area of $\Sigma$ and the volume enclosed by $\Sigma$ in the hyperbolic space ${\bf H}^3$. In other words, the critical points in ${\bf H}^3$ have mean curvature which is a linear function of their height. A more general functional for surfaces in ${\bf H}^3$ was studied in \cite{Rafa2}, while in ${\bf R}^3$ it was analyzed in \cite{KP} and \cite{KP2}.
\end{remark}

\subsection{Critica of $E$}

For axially symmetric immersions critical for $E$, we also need to consider the Euler-Lagrange equations \eqref{EL2}-\eqref{EL4} along the boundary circles. Let $X:\Sigma\rightarrow{\bf R}^3$ be an axially symmetric immersion critical for $E$. Using the boundary conditions \eqref{EL2}-\eqref{EL4} we will show that $\bar{A}=0$ and, hence, \eqref{nonCMC} holds on $\Sigma$, not only for topological discs (as in Theorem \ref{discsnew}) but also for annuli.

First, if $b=0$, it was shown in Proposition 2.3 of \cite{PP2} that either the surface has constant mean curvature $H=-c_o$ in which case the radii of the boundary circles are $\sqrt{\alpha/\beta}$, or the boundary is composed by closed geodesics. In the latter, from \eqref{EL2} we have that $H=-c_o$ along the geodesic boundary components, which combined with \eqref{fir}, shows that $\bar{A}=0$ and, as a consequence, \eqref{nonCMC} holds on $\Sigma$. We summarize this in the following result.

\begin{prop}\label{b=0prop} Let $X:\Sigma\rightarrow{\bf R}^3$ be an axially symmetric immersion critical for $E$ with $b=0$. Then, either the critical domain is a part of a Delaunay surface with $H\equiv -c_o$ in which case the radii of the boundary circles are $\sqrt{\alpha/\beta}$, or, after a suitable vertical translation, the mean curvature of the immersion satisfies,
\begin{equation*}\label{H+c}
H+c_o=-\frac{\nu_3}{z}\,,
\end{equation*}
on $\Sigma$ and the boundary $\partial\Sigma$ is composed by geodesic circles.
\end{prop}

The case $b\neq 0$ is richer. Indeed, among axially symmetric domains critical for $E$ satisfying \eqref{nonCMC}, we can prove the existence of three essentially different types of boundary circles.

\begin{theorem}\label{4.1} Let $X:\Sigma\rightarrow{\bf R}^3$ be an axially symmetric immersion critical for $E$ with $b\neq 0$. Then, either the critical domain is a part of a Delaunay surface bounded by circles of radii $\sqrt{\alpha/\beta}$, or, after a suitable vertical translation, the mean curvature of the immersion satisfies, 
$$H+c_o=-\frac{\nu_3}{z}\,,$$
on $\Sigma$ and each boundary component is one of the followings:
\begin{enumerate}[(i)]
\item A geodesic circle located in the plane $\{z=0\}$.
\item A non geodesic circle of radius $\sqrt{\alpha/\beta}$.
\item A non geodesic circle of radius $r\neq\sqrt{\alpha/\beta}$ along which $a\kappa_g r\pm b\kappa_n z\equiv 0$ holds.
\end{enumerate}
Moreover, along boundary components of case (i) and (ii), $a\left(H+c_o\right)^2+bK \equiv 0$ holds.
\end{theorem}
{\it Proof.\:} It was shown in \cite{PP2} that if a critical surface has constant mean curvature, then either $a=-b>0$, $c_o=0$ and the surface is a compact domain in a sphere (Theorem 3.1 of \cite{PP2}) or the mean curvature must be $-c_o$. In both cases the boundary circles have radii $\sqrt{\alpha/\beta}$. This shows the first statement.

Since the surface is axially symmetric \eqref{fir} holds on each parallel and since every parallel is homologous to a boundary cycle, we conclude from Theorem \ref{zeroflux} that $\bar{A}=0$ holds. As shown above, for non constant mean curvature cases, this implies that \eqref{nonCMC} holds also.

Assume now that the mean curvature of the immersion is not constant. Since the surface is axially symmetric, the geodesic torsion of each latitude vanishes, i.e. $\tau_g\equiv 0$, while $\kappa_g$ and $\kappa_n$ are constants. Therefore, the Euler-Lagrange equations \eqref{EL2}-\eqref{EL4} along the boundary $\partial\Sigma$ become:
\begin{eqnarray}
a\left(H+c_o\right)+b\kappa_n&=&0\,, \label{nsys1}\\
\left(\alpha r^{-2}-\beta\right)\kappa_n-a\partial_nH&=&0\,, \label{nsys2} \\
\left(\alpha r^{-2}-\beta\right)\kappa_g+a\left(H+c_o\right)^2+bK&=&0\,.\label{nsys3}
\end{eqnarray}

Next, we denote by $\varsigma$ the arc length parameter along the generating curve so that the derivative with respect to $\varsigma$ along $C_j$, $j=1,2$, coincides with the derivative in the conormal direction $n$. Then, for any value of $\varsigma$ for which $z(\varsigma)\ne 0$, we use equation \eqref{nonCMC} to compute
$$H_\varsigma=-\left(\frac{\nu_3}{z}\right)_\varsigma= \left(\kappa_1-\left[H+c_o\right]\right)\frac{z_\varsigma}{z}\,,$$
where $\kappa_1$ is the principal curvature such that $\kappa_1=2H-\kappa_n$ along $\partial\Sigma$. By continuity of $H_\varsigma$, this formula is also valid when $z=0$. Using this together with \eqref{nsys1}, equations \eqref{nsys2} and \eqref{nsys3} can be expressed as
\begin{eqnarray*}
\left(\alpha r^{-2}-\beta\right)\kappa_n+a\left(\kappa_n-H+c_o\right)\frac{\partial_nz}{z}&=&0\,,\\
\left(\alpha r^{-2}-\beta\right)\kappa_g+a\left(H+c_o\right)\left(\kappa_n-H+c_o\right)&=&0\,.
\end{eqnarray*}
From this, the rest of the statements in the theorem are clear. {\bf q.e.d.}
\\

We end this section proving that, under certain conditions, the symmetries of the boundary are inherited by the interior of the critical surface. Since the Euler-Lagrange equation \eqref{EL1} is fourth order, it is expected that these conditions should depend on higher order boundary values of the surface.  Indeed, if one boundary component of a critical immersion for $E$ is a circle and the surface meets it in a constant angle,  then the surface is axially symmetric. A first result in this direction was proven in Proposition 5.1 of \cite{PP2} (\emph{a constant mean curvature immersion whose boundary contains at least one circle on which $\tau_g\equiv 0$ holds, is axially symmetric}). Here, we extend this result to more general solutions of \eqref{EL1}.

\begin{lemma}\label{axsym} Let $X:\Sigma\rightarrow{\bf R}^3$ be an immersion of a compact surface with boundary satisfying
$$\Delta H+2\left(H+c_o\right)\left(H\left[H-c_o\right]-K\right)=0\,,$$
on $\Sigma$. If any connected component of $\partial\Sigma$ is a circle on which $\tau_g\equiv 0$ holds  and along which $H$ and $\partial_nH$ are constant, then the surface is axially symmetric.
\end{lemma}
{\it Proof.\:} Let $C$ be a circle representing a connected component of $\partial\Sigma$. We can assume that the circle $C$ lies in a horizontal plane. Let $\mathcal{R}_t$ denote the one parameter family of rotations about a vertical axis passing through the circle's center. We define the function
$$\psi:=\partial_t\left(\mathcal{R}_tX\right)_{t=0}\cdot\nu=E_3\times X\cdot \nu\,,$$
which is the normal component of the variation of $X$ obtained by rotating the surface about the vertical axis.

Next, using the notation introduced in \eqref{L}, the PDE in the statement reads $L[H+c_o]=0$. We regard this equation as being a condition on the immersion and denote it by
$F[X]=0$.  The linearization of this equation, at a solution of $F[X]=0$, is defined by 
$$DF_X[u]:=\partial_\epsilon F[X+\epsilon u\nu]_{\epsilon=0}\,,$$
which is a fourth order elliptic operator with analytic coefficients, since the surface is real analytic by elliptic regularity (see Theorem 6.6.1 of \cite{Morrey}). Then, since $\psi$ is the normal component of the Killing field $E_3\times X$, it is clear that the function $\psi$ is a solution of the linearized equation $DF_X[\psi]=0$.

Note that $\psi\equiv 0$ on $C$, since the boundary circle is preserved by the rotation $\mathcal{R}_t$. This implies that $\psi'=\nabla\psi\cdot T\equiv 0$ along $C$. Also, along $C$, we have
\begin{eqnarray*}
\partial_n\psi&=&E_3\times n\cdot \nu+E_3\times X\cdot d\nu(n)\\
&=&-T\cdot E_3-\tau_g E_3\times X\cdot T-\left(2H-\kappa_n\right)E_3\times X\cdot n\\
&\equiv&0\,,
\end{eqnarray*}
since $E_3$ is normal to $C$, $\tau_g\equiv 0$ along $C$ and $E_3\times X$ is tangent to $C$. Combining with $\psi'\equiv0$, we conclude that $\nabla\psi\equiv 0$ along $C$.

We next claim that $\partial_n^2\psi\equiv 0$ also holds along $C$. To see this, we note that the pointwise variation of the mean curvature with respect to the Killing field $E_3\times X$ is zero. This means (Appendix A of \cite{PP2}):
\begin{equation}\label{varH}
0=\frac{1}{2}\left(\Delta \psi+\lVert d\nu\rVert^2\psi\right)+\nabla H\cdot E_3\times X\:.
\end{equation}
Along $C$, $\psi\equiv 0$ holds. Moreover, since $H$ is constant along $C$ and $E_3\times X$ is tangent to $C$, $\nabla H\cdot E_3\times X$ also vanishes. Therefore, we conclude that
$$0=\Delta \psi=\nabla\cdot\left(\nabla\psi\right)=\nabla _n\nabla \psi \cdot n+\nabla _T\nabla \psi\cdot T=\nabla _n\nabla \psi \cdot n\,,$$
where the last equality holds since $\nabla\psi\equiv 0$ on $C$ and, hence, $\nabla_T \nabla\psi\equiv 0$ along $C$. The claim follows from, $\nabla _n\nabla \psi \cdot n=\partial_n^2\psi-\nabla \psi\cdot \nabla_nn =\partial_n^2\psi$.

We show now that $\partial_n^3\psi\equiv 0$ along $C$. From \eqref{varH}, we obtain along $C$
\begin{eqnarray*}
0&=&\partial_n\left(\Delta \psi\right)+\partial_n\left(\lVert d\nu\rVert^2\right)\psi+\lVert d\nu\rVert^2\partial_n\psi+2E_3\times n\cdot \nabla H\\&&+2E_3\times X\cdot \nabla_n\nabla H=\partial_n\left(\Delta \psi\right)+2E_3\times X\cdot \nabla_n\nabla H\,,
\end{eqnarray*}
since, along $C$, $\psi=\partial_n\psi\equiv 0$ and $H$ is constant. Therefore,
$$\partial_n\left(\Delta \psi\right)=-2E_3\times X\cdot \nabla_n\nabla H\,,$$
which is a multiple of $\nabla_T\nabla H\cdot n$. However, $\nabla H=\partial_nH n$ along $C$ and, since $\partial_nH$ is constant along $C$, $\nabla_T\nabla H\cdot n=0$ holds. Consequently, 
$$0=\partial_n\left(\Delta \psi\right)=\partial_n\left(\nabla_n\nabla\psi\cdot n\right)+\partial_n\left(\nabla_T\nabla\psi\cdot T\right)=\partial_n^3\psi+\nabla_n\nabla_T\nabla\psi\cdot T\,,$$
since $\nabla_T\nabla\psi\equiv 0$ along $C$. We focus on the last term, obtaining
$$\nabla_n\nabla_T\nabla\psi\cdot T=\nabla_T\nabla_n\nabla\psi+\nabla_{[n,T]}\nabla\psi+R(n,t)\nabla\psi\,,$$
where $R$ is the curvature tensor of $\Sigma$. Observe that the three terms in above equation vanish since $\nabla\psi\equiv 0$ and $\nabla\nabla\psi=0$ along $C$. Thus, we conclude that $\partial_n^3\psi\equiv 0$ along $C$.

Finally, we apply the Cauchy-Kovalevskaya Theorem for the Cauchy problem $DF_X[u]=0$ with initial conditions $0=u=\partial_nu=\partial_{n}^2u=\partial_{n}^3u$ on the real analytic set $C$, obtaining that this problem, locally, has a unique real analytic solution, $u\equiv 0$, so we obtain $\psi\equiv 0$ locally. Once again, using analyticity of the surface, we obtain that $\psi\equiv 0$ \emph{globally} on $\Sigma$. From this we can conclude that the surface is axially symmetric. {\bf q.e.d.}
\\

Applying Lemma \ref{axsym}, we now prove sufficient conditions for the boundary of a critical surface, in order to obtain axial symmetry.

\begin{theorem}\label{thmconv} Let $X:\Sigma\rightarrow{\bf R}^3$ be an immersion critical for $E$. If one boundary component is a circle and the surface meets the plane of the circle in a constant angle, then the surface is axially symmetric.
\end{theorem}
{\it Proof.\:} Let $X:\Sigma\rightarrow{\bf R}^3$ be a critical immersion for $E$. Then, the Euler-Lagrange equations \eqref{EL1}-\eqref{EL4} hold. From the classical Joachimsthal Theorem, along the circular boundary component $C$, $\kappa_g$ and $\kappa_n$ are constants, while $\tau_g\equiv 0$ holds. Therefore, along $C$, \eqref{EL2}-\eqref{EL4} reduce, respectively, to \eqref{nsys1}-\eqref{nsys3}. In particular, we deduce that $H$ and $\partial_nH$ are constants along $C$. Thus, we can apply Lemma \ref{axsym}, proving the result. {\bf q.e.d.}
\\

We point out  that the result above cannot be improved. Indeed, to get axially symmetric domains, it is clear that boundary components must be circles. Moreover, requesting constant contact angle along a circular boundary component is also essential. Without assuming this, there are non-axially symmetric domains bounded by circles critical for $E$, for instance, suitable domains in Riemann's minimal examples (see Figure 4 in \cite{PP2}).

\section{Examples of Axially Symmetric Critical Discs}\label{Examples}

In this section, specific examples of axially symmetric discs critical for $E$ are obtained. Let $X:\Sigma\rightarrow{\bf R}^3$ be the immersion of a disc type surface critical for $E$. We will focus on the case where \eqref{nonCMC} holds, after a suitable vertical translation, since the ground state $H+c_o\equiv 0$ was studied in detail in \cite{PP2}. 

Denote by $\gamma(\sigma)=\left(r(\sigma),z(\sigma)\right)$ the arc length parameterized profile curve of the axially symmetric surface satisfying \eqref{nonCMC}. The mean curvature can be expressed, up to a sign, as
\begin{equation}\label{Haxsym}
H=\frac{z'(\sigma)}{2r(\sigma)}-\frac{r''(\sigma)}{2z'(\sigma)}\,,
\end{equation}
where $\left(\,\right)'$ denotes the derivative with respect to the arc length parameter $\sigma$. Then, a straightforward computation shows that \eqref{nonCMC} can be rewritten as
$$\frac{z'(\sigma)}{2r(\sigma)}-\frac{r''(\sigma)}{2z'(\sigma)}\pm c_o=-\frac{r'(\sigma)}{z(\sigma)}\,.$$
Since $\gamma(\sigma)$ is parameterized by arc length, we can introduce a function $\varphi(\sigma)$, so that $r'(\sigma)=\cos\varphi(\sigma)$ and $z'(\sigma)=\sin\varphi(\sigma)$. The function $\varphi(\sigma)$ represents the angle between the positive part of the $r$-axis and the tangent vector to $\gamma(\sigma)$. With this notation, we conclude from \eqref{Haxsym} that, up to a sign,
\begin{equation}\label{HAS}
H=\frac{1}{2}\left(\varphi'+\frac{\sin\varphi}{r}\right).
\end{equation}
Therefore, using \eqref{HAS}, the profile curve $\gamma(\sigma)$ satisfies the following system of first order ordinary differential equations:
\begin{eqnarray}
r'(\sigma)&=&\cos\varphi(\sigma)\,,\label{dif1}\\
z'(\sigma)&=&\sin\varphi(\sigma)\,,\label{dif2}\\
\varphi'(\sigma)&=&-2\frac{\cos\varphi(\sigma)}{z(\sigma)}-\frac{\sin\varphi(\sigma)}{r(\sigma)}\mp 2c_o\,.\label{dif3}
\end{eqnarray}

\begin{remark}\label{5.1}\normalfont Up to the transformation $z\mapsto -z$,  the sign in front of $c_o$ in \eqref{dif3} can be fixed to be negative. Considering the opposite sign in \eqref{dif3} also carries a different sign in front of the terms with $c_o$ in the boundary conditions \eqref{bcdif1}-\eqref{bcdif2} and \eqref{bcnew2}.
\end{remark}

We now impose the initial conditions ($\sigma=0$). Since we want to obtain a disc type surface, the profile curve $\gamma(\sigma)$ must cut the $z$-axis, i.e. $r(0)=0$. Moreover, the regularity of the surface implies that this cut with the $z$-axis must be perpendicular. This completely describes the initial tangent vector, which is equivalent to fixing the initial value $\varphi(0)=0$. Finally, the initial height can be assumed to be a parameter, that is, $z(0)=z_o\neq 0$.

\begin{remark}\label{co}\normalfont If $c_o=0$, for any initial height $z_o\neq 0$, the uniqueness part of the Cauchy-Kovalevskaya Theorem for the Cauchy problem \eqref{dif1}-\eqref{dif3} with initial conditions $r(0)=0$, $z(0)=z_o$ and $\varphi(0)=0$ shows that the sphere of suitable radius is the only solution. Indeed, in Theorem 1.1 of \cite{P}, a more general result was obtained for this case using Bryant's quartic differential.
\end{remark} 

Next, in order to obtain critical discs for $E$, boundary conditions at $\sigma=\mathcal{L}$, where $\mathcal{L}$ denotes the length of $\gamma$, must also be imposed. These boundary conditions differ depending on whether $b=0$ or not (c.f. Proposition \ref{b=0prop} and Theorem \ref{4.1}).

First consider the case $b=0$. In this case, for non constant mean curvature surfaces, we know from Proposition \ref{b=0prop} that the boundary is a geodesic circle whose radius is determined by \eqref{EL3}. Moreover, since \eqref{nonCMC} holds on $\Sigma$, the remaining Euler-Lagrange equations are automatically satisfied. 

\begin{figure}[h!]
\hspace{0.75cm}
\begin{subfigure}[b]{0.32\linewidth}
\includegraphics[width=\linewidth]{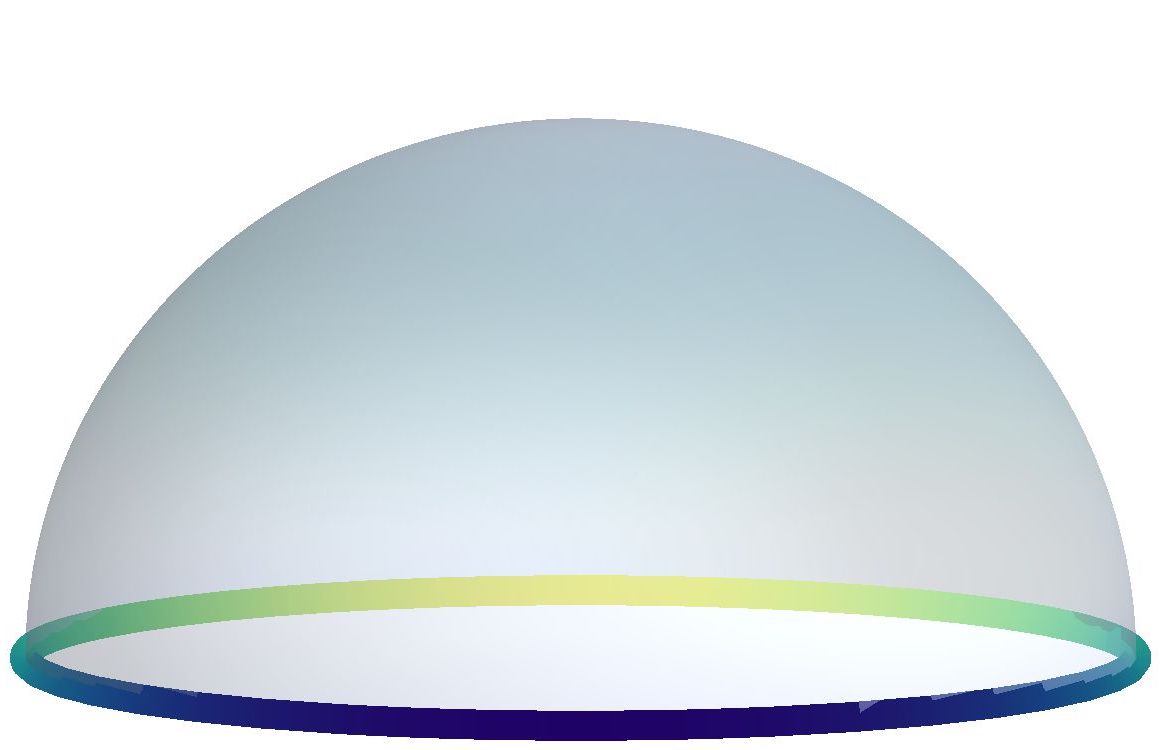}
\caption{$c_o=1.1$}
\end{subfigure}
\hspace{0.3cm}
\begin{subfigure}[b]{0.32\linewidth}
\includegraphics[width=\linewidth]{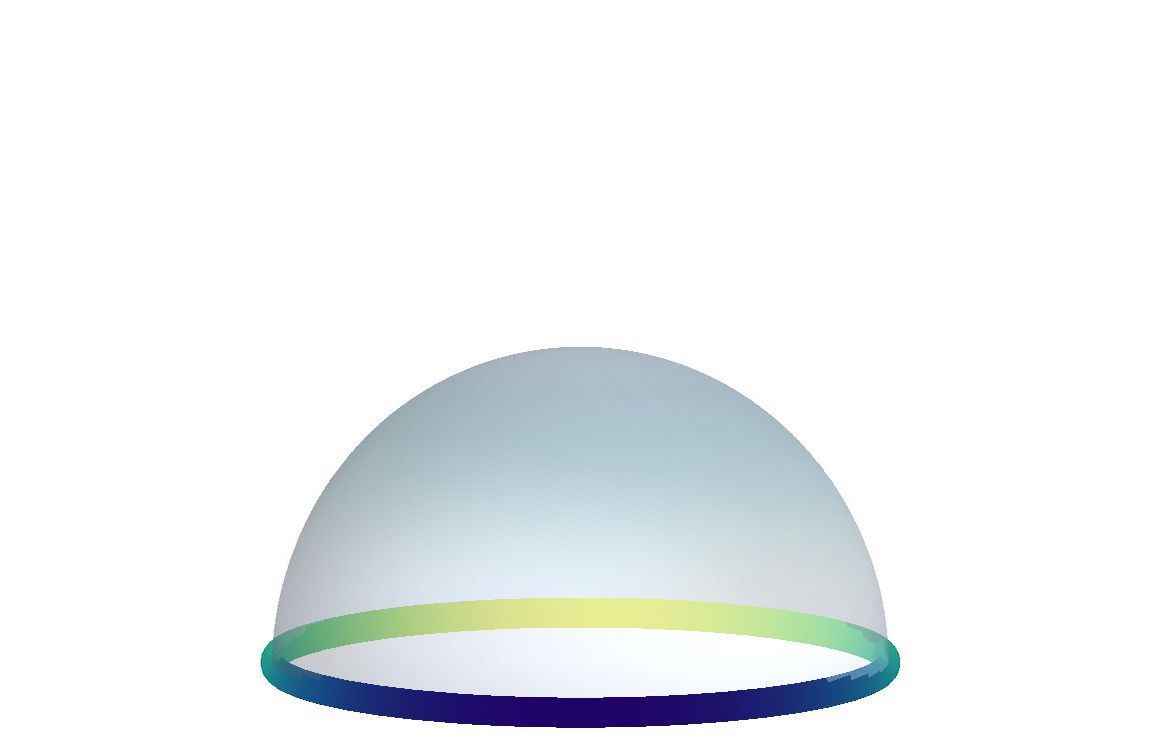}
\caption{$c_o=2$}
\end{subfigure}
\hspace{-1cm}
\begin{subfigure}[b]{0.32\linewidth}
\includegraphics[width=\linewidth]{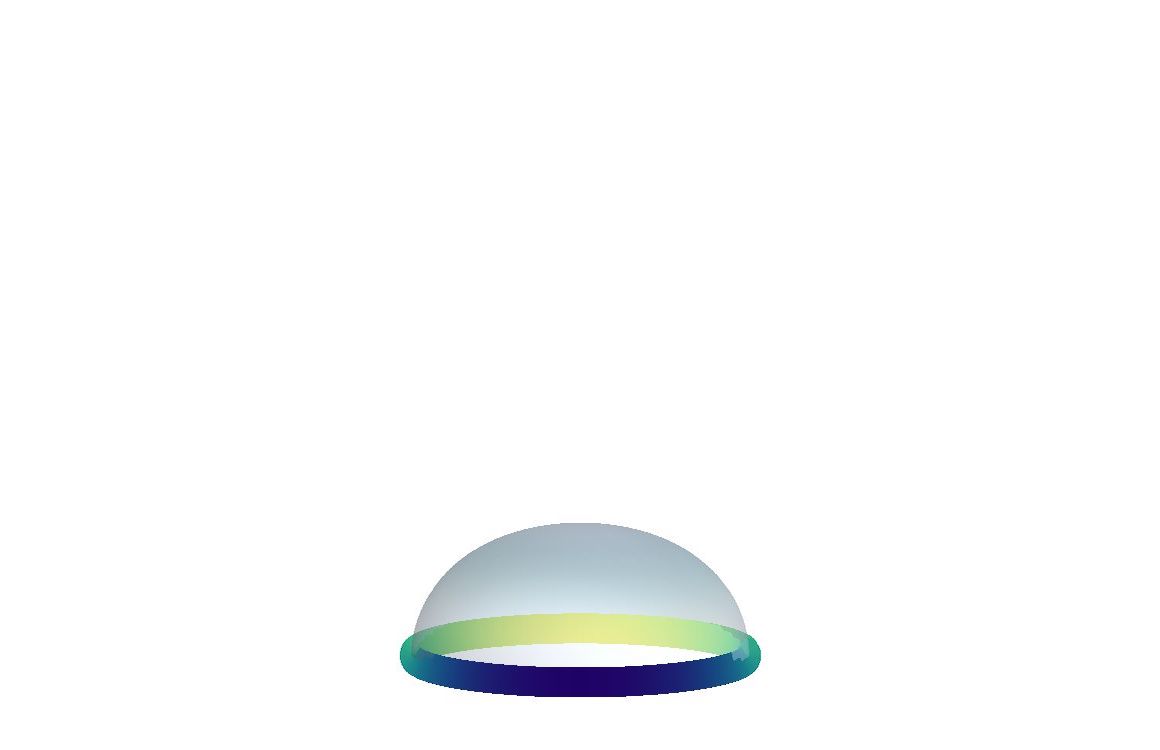}
\caption{$c_o=5$}
\end{subfigure}
\caption{Axially symmetric critical discs for the energy $E_{a=1,c_o,b=0,\alpha=1,\beta=1}$ and different values of $c_o>0$ with $c_o^2>\beta/\alpha$. In all the cases the boundary is a geodesic circle of radius $r<\sqrt{\alpha/\beta}$.}
\label{b=0}
\end{figure}

For fixed energy parameters $a>0$, $b=0$, $\alpha>0$ and $\beta>0$, we show in Figure \ref{b=0} three axially symmetric discs with non constant mean curvature critical for the energy $E$. Each of these domains corresponds with different values of $c_o>0$ satisfying $c_o^2>\beta/\alpha$. As $c_o$ increases, critical domains get more ``planar" and the radius of their boundary circle, which is always smaller than $\sqrt{\alpha/\beta}$, decreases.

Consider now the case $b\neq 0$ and assume that $a(H+c_o)^2+bK=0$ holds on $\partial\Sigma$. Then, from Theorem \ref{4.1}, the boundary parallel is either a geodesic whose radius is given by \eqref{EL3} or $r(\mathcal{L})=\sqrt{\alpha/\beta}$, necessarily. In addition to fixing the radius of the boundary circle, the following two equations must also hold along $\partial\Sigma$,
\begin{eqnarray}
\kappa_n-H+c_o&=&0\,,\label{Hknco}\\
\left(a+b\right)\kappa_n+2ac_o&=&0\,.\label{sec}
\end{eqnarray}
Equation \eqref{Hknco} is a consequence of $a(H+c_o)^2+bK=0$ along the boundary, while \eqref{sec} arises after combining \eqref{Hknco} with \eqref{EL2}. 

We rewrite \eqref{Hknco}-\eqref{sec} in terms of the variables of the system of differential equations, respectively, as
\begin{eqnarray}
\varphi'(\mathcal{L})=\sin\varphi(\mathcal{L})/r(\mathcal{L})+2c_o\,,&&\label{bcdif1}\\
\left(a+b\right)\sin\varphi(\mathcal{L})+2a c_or(\mathcal{L})=0\,.&&\label{bcdif2}
\end{eqnarray}
These equations give some restrictions for the existence of critical domains. For instance, if the boundary is a geodesic circle, then from \eqref{bcdif2},
\begin{equation}\label{radiusgeodesic}
r^2(\mathcal{L})=\frac{\left(a+b\right)^2}{4a^2c_o^2}
\end{equation}
holds. On the other hand, if the boundary circle is not a geodesic, $r(\mathcal{L})=\sqrt{\alpha/\beta}$ holds, which combined with \eqref{bcdif2}, shows that the following restriction must be satisfied by the energy parameters:
\begin{equation}\label{radiusnongeodesic}
c_o^2<\frac{\left(a+b\right)^2\beta}{4a^2\alpha}\,.
\end{equation}

\begin{remark}\normalfont If $c_o\neq 0$ and $a=-b>0$ holds, there are no axially symmetric discs critical for $E$ satisfying $a(H+c_o)^2+bK\equiv 0$ along the boundary $\partial\Sigma$.
\end{remark}

To get specific examples, we solve the system of differential equations \eqref{dif1}-\eqref{dif3} with the initial conditions $r(0)=0$, $z(0)=z_o\neq 0$, $\varphi(0)=0$, together with the boundary conditions \eqref{bcdif1}-\eqref{bcdif2} and such that $r(\mathcal{L})$ satisfies \eqref{EL3}. 

\begin{figure}[h!]
\makebox[\textwidth][c]{
\begin{subfigure}[b]{0.22\linewidth}
\includegraphics[width=\linewidth]{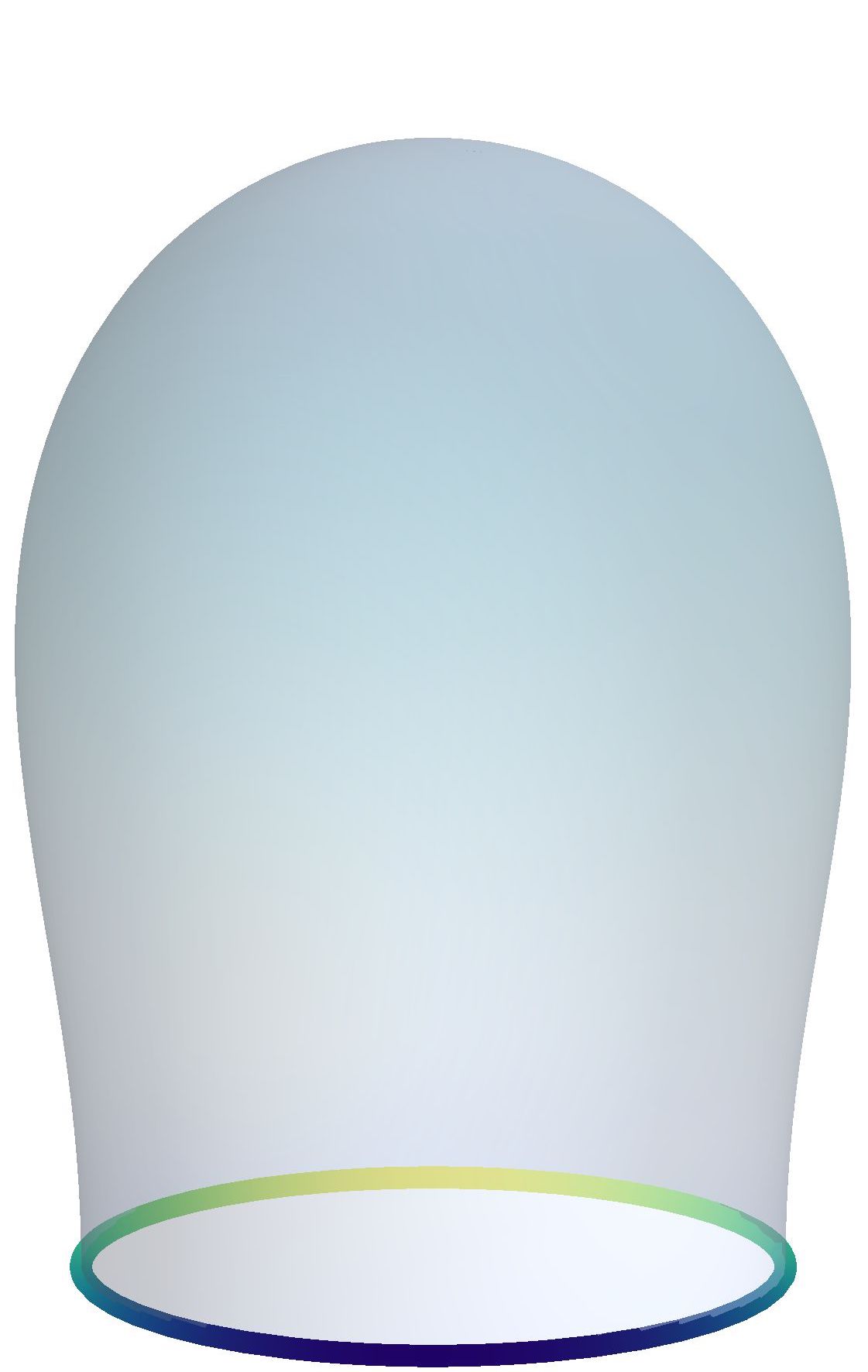}
\caption{$E_{1,2,0.25,1,10}$}
\end{subfigure}
\quad\quad\quad\,
\begin{subfigure}[b]{0.22\linewidth}
\includegraphics[width=\linewidth]{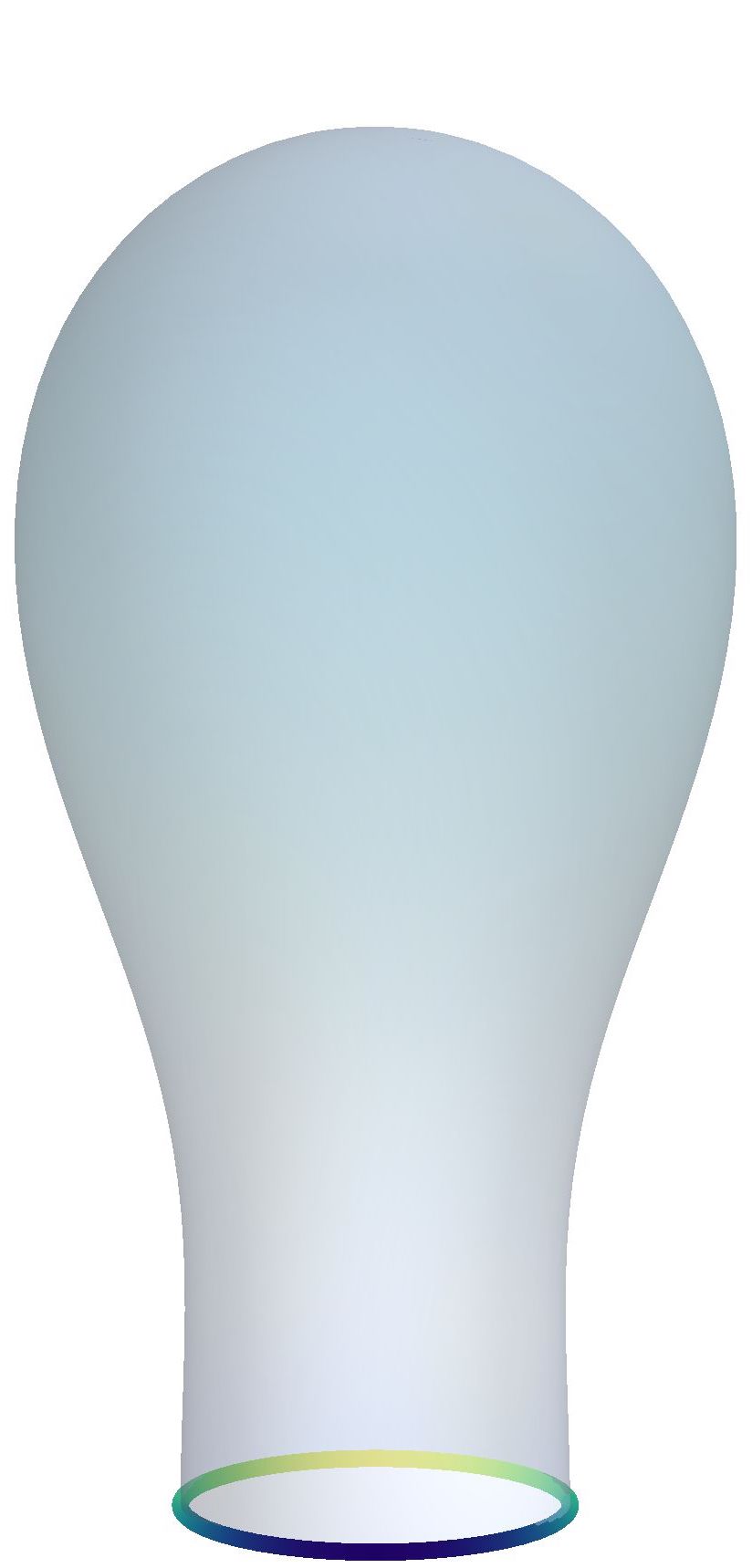}
\caption{$E_{1,2,-0.11,1,20}$}
\end{subfigure}
\quad\quad\quad
\begin{subfigure}[b]{0.22\linewidth}
\includegraphics[width=\linewidth]{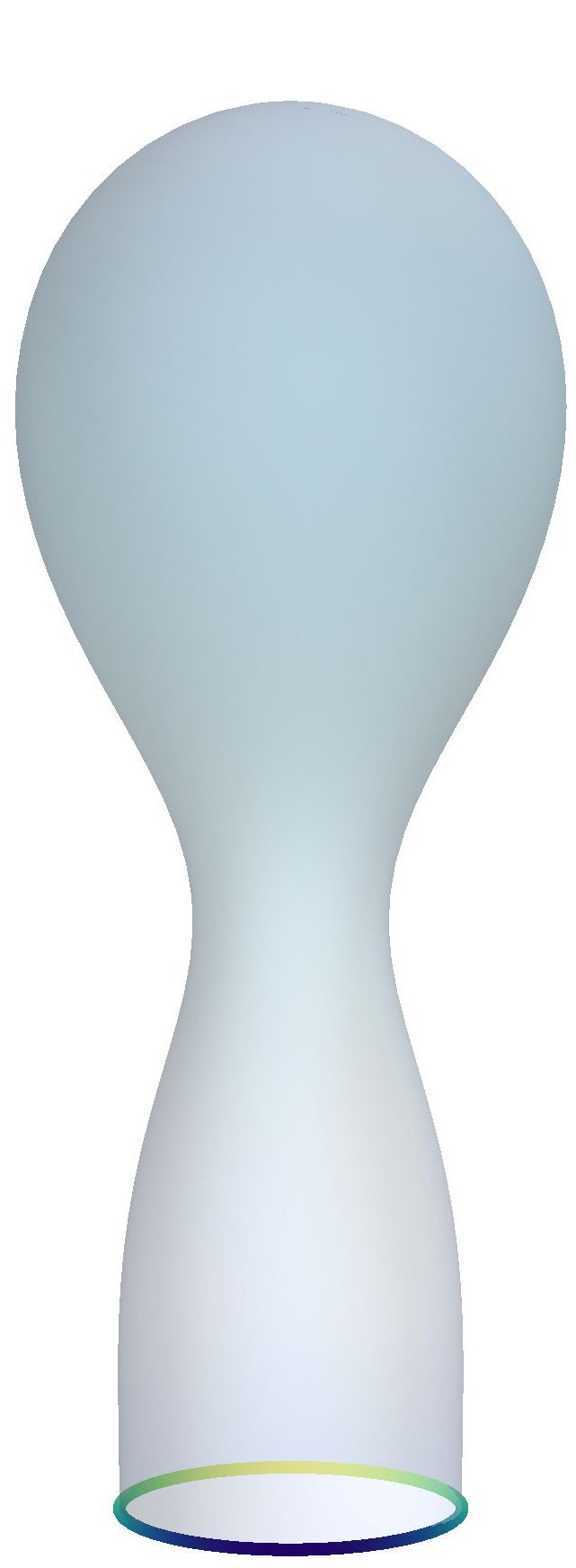}
\caption{$E_{1,2,0.08,1,14}$}
\end{subfigure}
\quad\quad\quad
\begin{subfigure}[b]{0.22\linewidth}
\includegraphics[width=\linewidth]{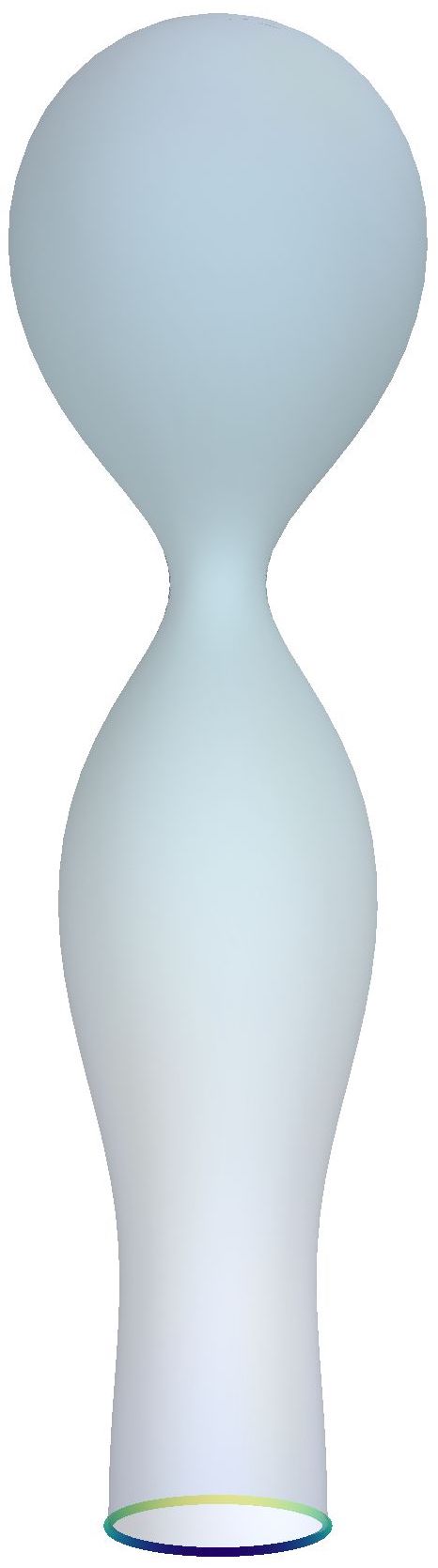}
\caption{$E_{1,2,-0.05,1,18}$}
\end{subfigure}}
\caption{Axially symmetric critical discs for the energy $E_{a,c_o,b,\alpha,\beta}$ bounded by a geodesic (parallel) circle located in the plane $\{z=0\}$ whose radius satisfies \eqref{radiusgeodesic}.}
\label{geo}
\end{figure}

In Figure \ref{geo},  different axially symmetric discs critical for the energy $E$ with $c_o>0$ are shown, whose boundary is a geodesic parallel located in the plane $\{z=0\}$. 

\begin{figure}[h!]
\makebox[\textwidth][c]{
\hspace{-0.5cm}
\begin{subfigure}[b]{0.23\linewidth}
\includegraphics[width=\linewidth]{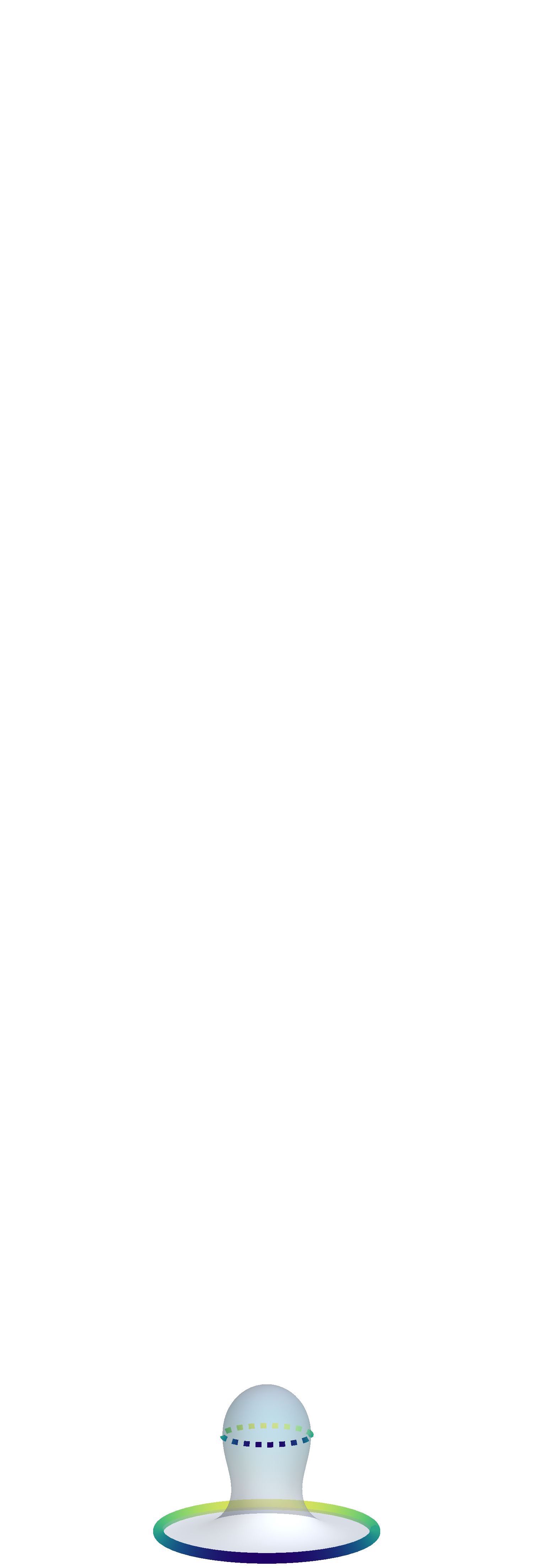}
\caption{$E_{1,2,-5,1,20}$}
\end{subfigure}
\hspace{-0.6cm}
\begin{subfigure}[b]{0.23\linewidth}
\includegraphics[width=\linewidth]{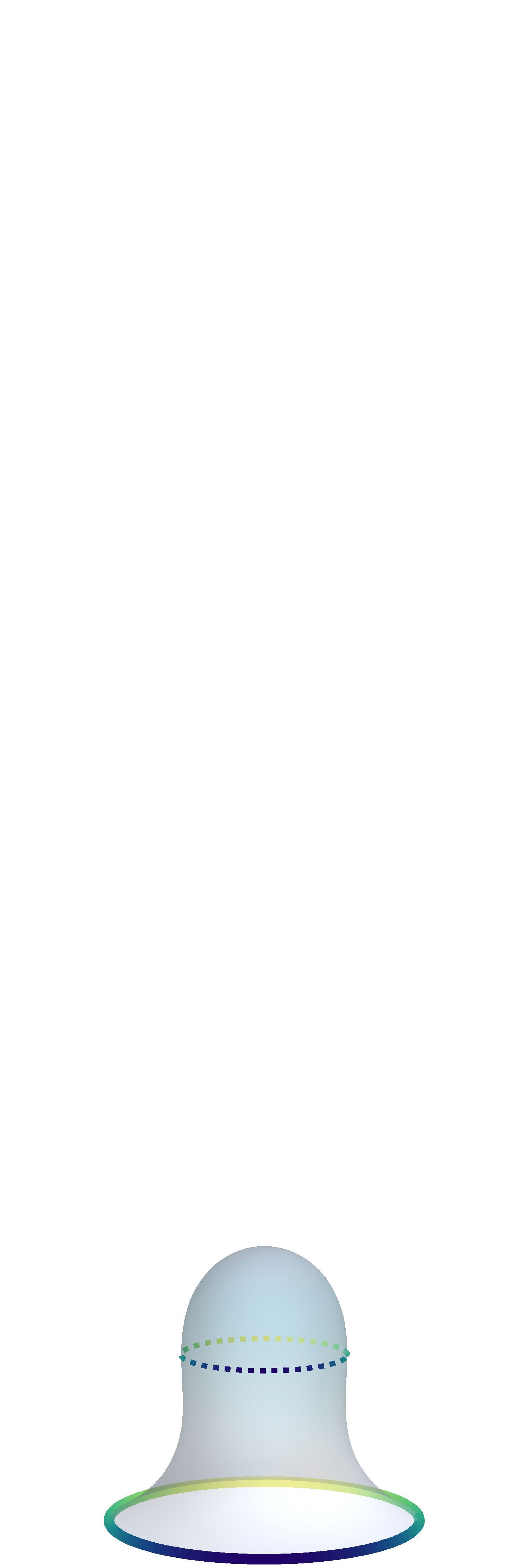}
\caption{$E_{1,2,2,1,10}$}
\end{subfigure}
\hspace{-0.3cm}
\begin{subfigure}[b]{0.23\linewidth}
\includegraphics[width=\linewidth]{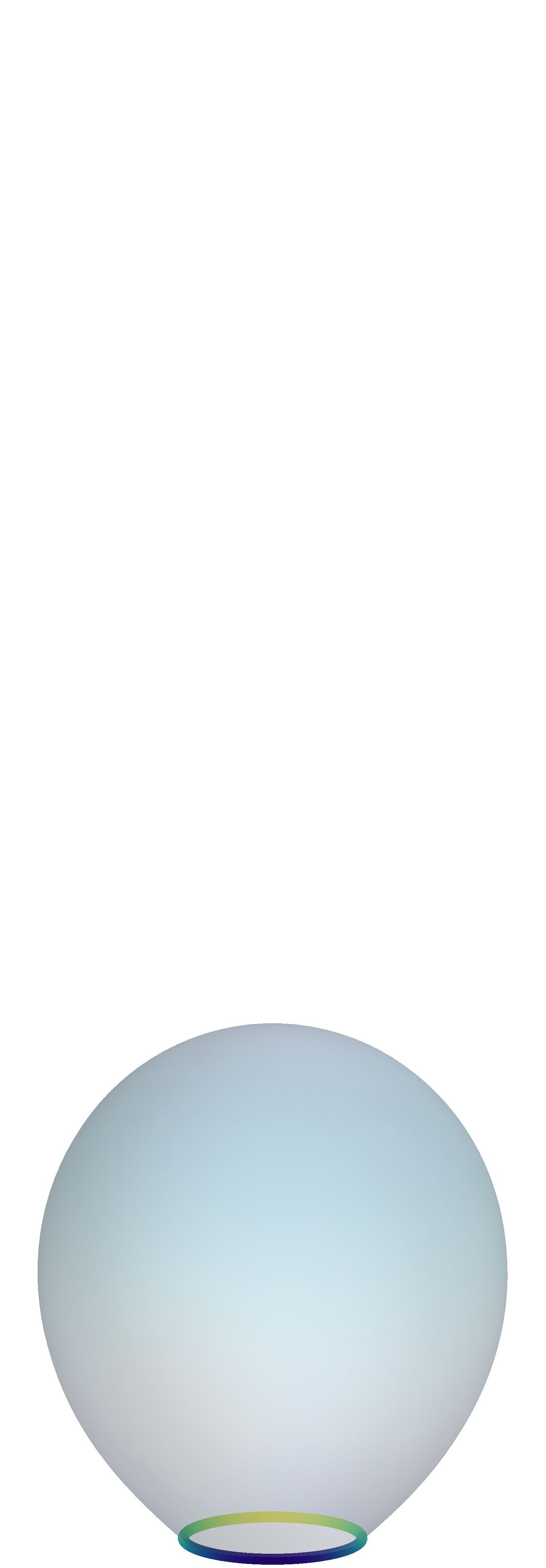}
\caption{$E_{1,2,0.05,1,31}$}
\end{subfigure}
\begin{subfigure}[b]{0.23\linewidth}
\includegraphics[width=\linewidth]{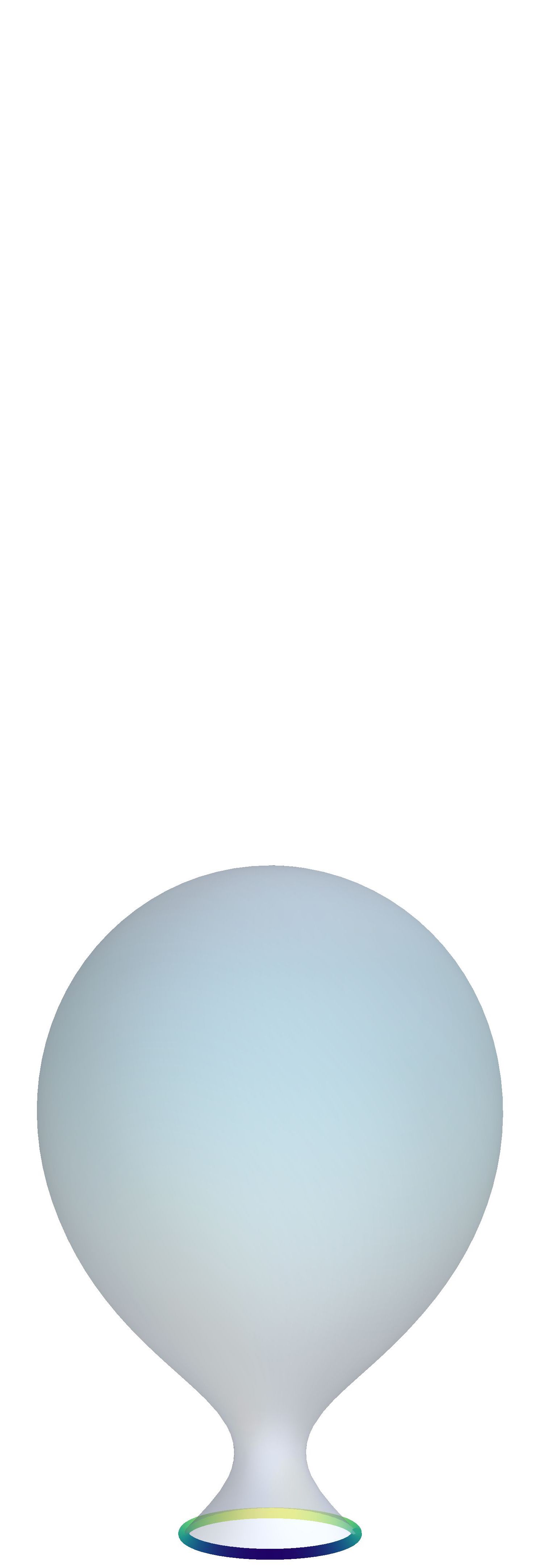}
\caption{$E_{1,2,-0.05,1,33}$}
\end{subfigure}
\begin{subfigure}[b]{0.23\linewidth}
\includegraphics[width=\linewidth]{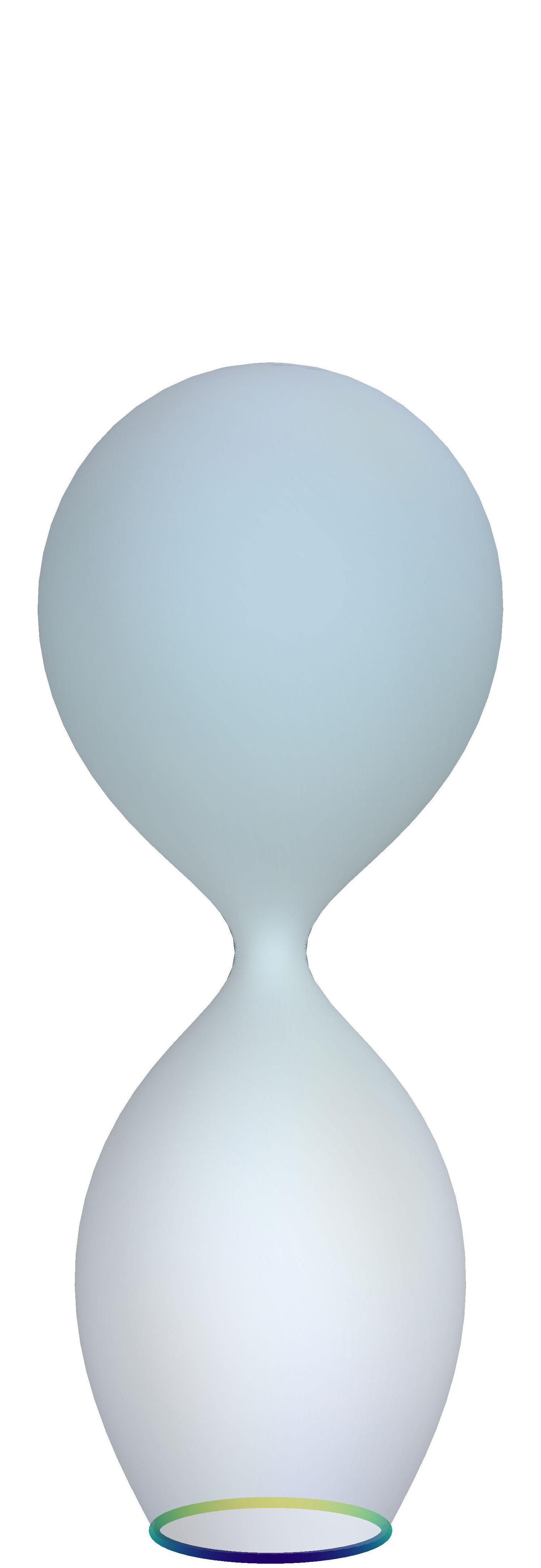}
\caption{$E_{1,2,0.05,1,18}$}
\end{subfigure}
\begin{subfigure}[b]{0.23\linewidth}
\includegraphics[width=\linewidth]{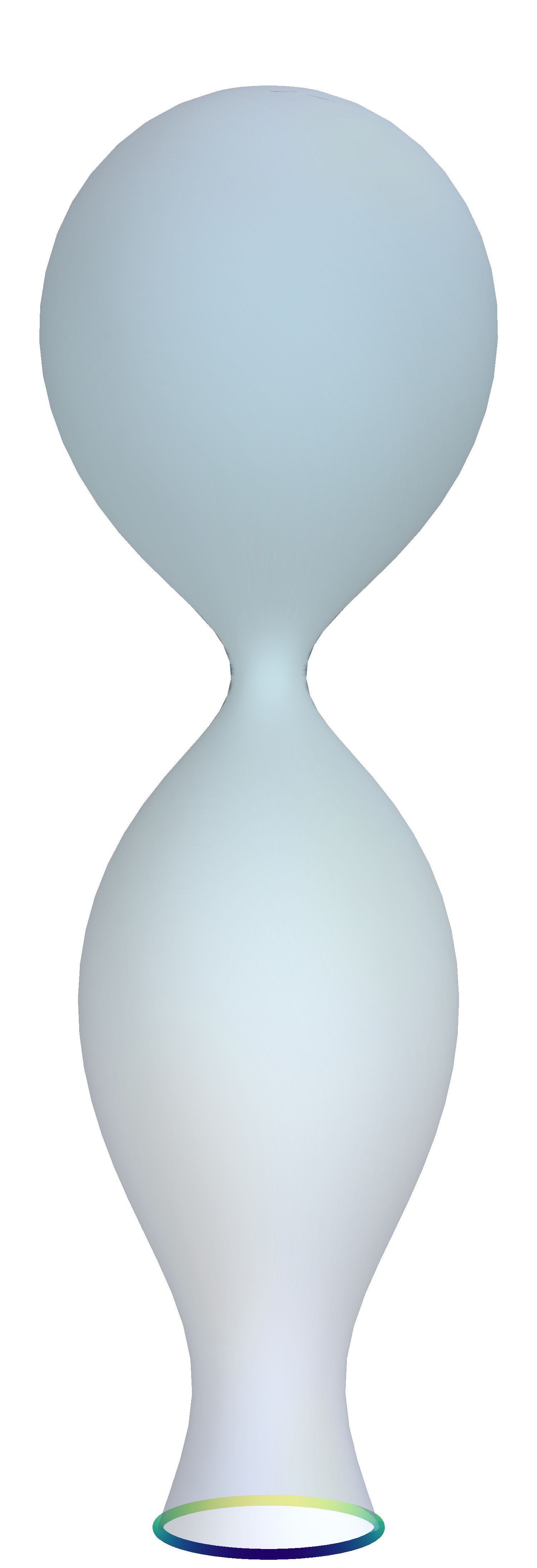}
\caption{$E_{1,2,-0.05,1,20}$}
\end{subfigure}}
\caption{Axially symmetric critical discs for the energy $E_{a,c_o,b,\alpha,\beta}$ bounded by a (non geodesic) circle of radius $\sqrt{\alpha/\beta}$.}
\label{bal}
\end{figure}

Similarly, in Figure \ref{bal} we show several axially symmetric discs critical for $E$ with $c_o>0$ bounded by a (non geodesic) circle of radius $\sqrt{\alpha/\beta}$. In this case, $\alpha\kappa^2=\beta$ holds along $\partial\Sigma$ and, the right hand side of equation \eqref{rescalings} vanishes. Consequently, since $c_o\neq 0$, we have
$$0=\int_\Sigma\left(H+c_o\right)d\Sigma=-\int_\Sigma\frac{\nu_3}{z}\,d\Sigma\,.$$
In particular, if the sign of $z$ does not change, necessarily the sign of $\nu_3$ must change, proving that in these cases critical domains cannot be graphs. This is the case of Figure \ref{bal}, (C)-(F). Moreover, if there is a change of sign in $z$, this must happen along a geodesic parallel and, hence, the normal along that parallel is horizontal so the surface cannot be a graph either. This phenomena can be appreciated in Figure \ref{bal}, (A)-(B) (see the dashed parallel).

Finally, we consider the remaining case of Theorem \ref{4.1}, i.e. the boundary is a non geodesic circle of radius $r(\mathcal{L})\neq\sqrt{\alpha/\beta}$ along which $a\kappa_gr+b\kappa_n z=0$ holds (the opposite sign corresponds with the transformation explained in Remark \ref{5.1}). With the notation introduced in \eqref{dif1}-\eqref{dif2}, this boundary condition rewrites
\begin{equation}\label{bcnew1}
a r(\mathcal{L})\cos\varphi(\mathcal{L})-b z(\mathcal{L})\sin\varphi(\mathcal{L})=0\,,
\end{equation}
while \eqref{EL4} gives us the boundary condition for the first derivative of the angle, namely,
\begin{equation}\label{bcnew2}
\varphi'(\mathcal{L})=2\left(\alpha r^{-2}(\mathcal{L})-\beta\right)\frac{z(\mathcal{L})}{a r(\mathcal{L})}+\frac{\sin\varphi(\mathcal{L})}{r(\mathcal{L})}+ 2c_o\,.
\end{equation}
Differentiating \eqref{nonCMC}, we can see that the Euler-Lagrange equation \eqref{EL3} is an identity.

\begin{figure}[h!]
\makebox[\textwidth][c]{
\hspace{-0.7cm}
\begin{subfigure}[b]{0.35\linewidth}
\includegraphics[width=\linewidth]{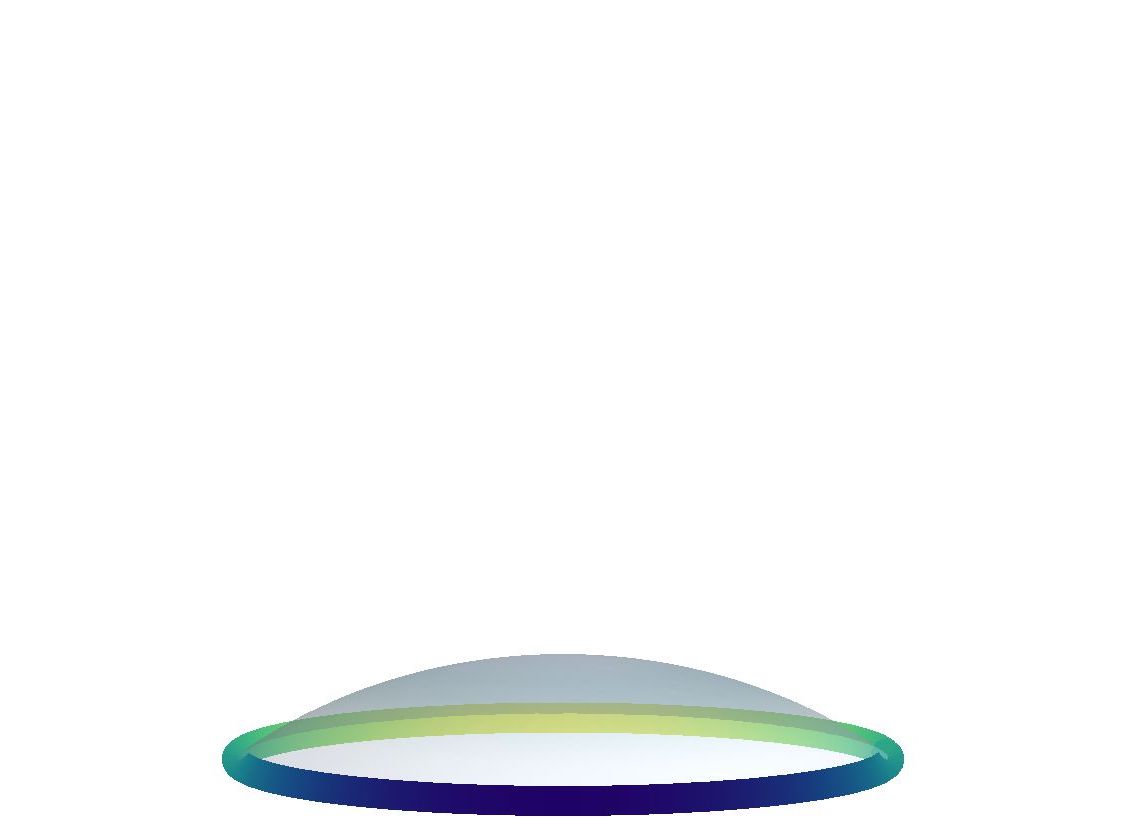}
\caption{$E_{1,2,0.5,1,4}$}
\end{subfigure}
\hspace{-0.4cm}
\begin{subfigure}[b]{0.35\linewidth}
\includegraphics[width=\linewidth]{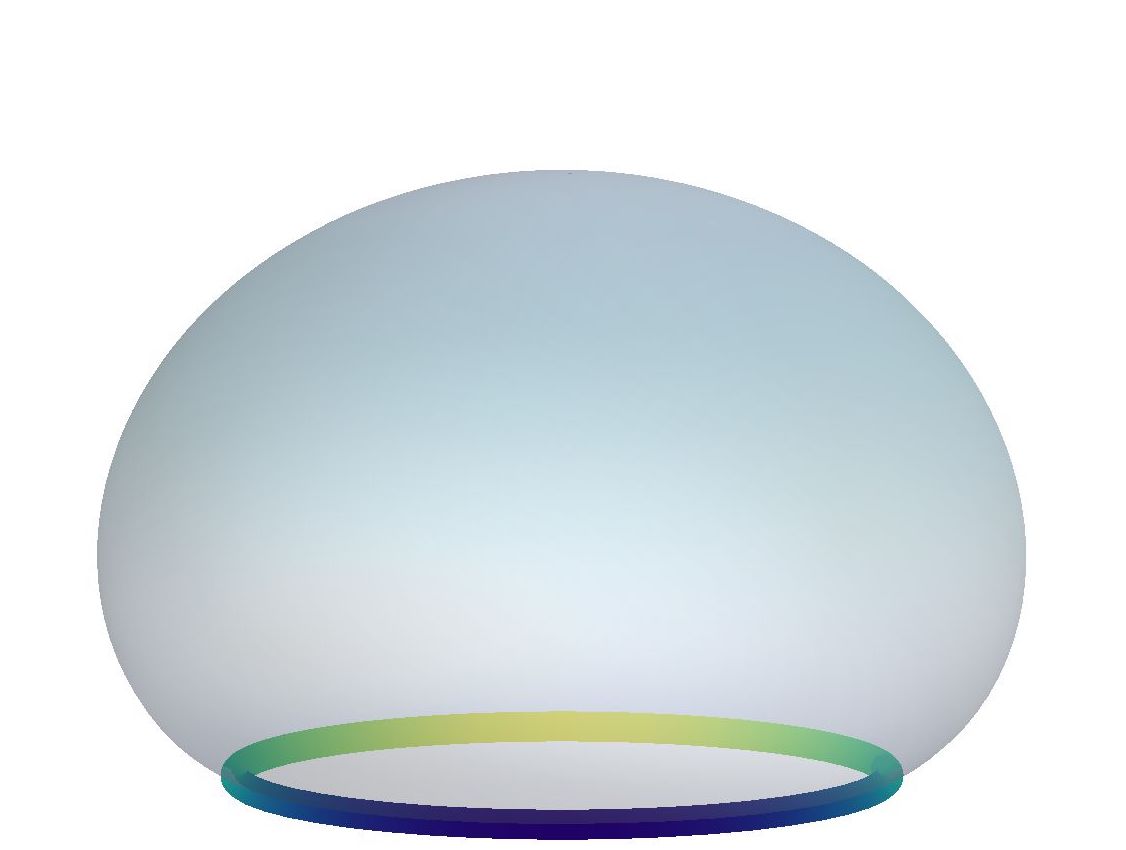}
\caption{$E_{1,2,-0.5,1,4}$}
\end{subfigure}
\,\,
\begin{subfigure}[b]{0.35\linewidth}
\includegraphics[width=\linewidth]{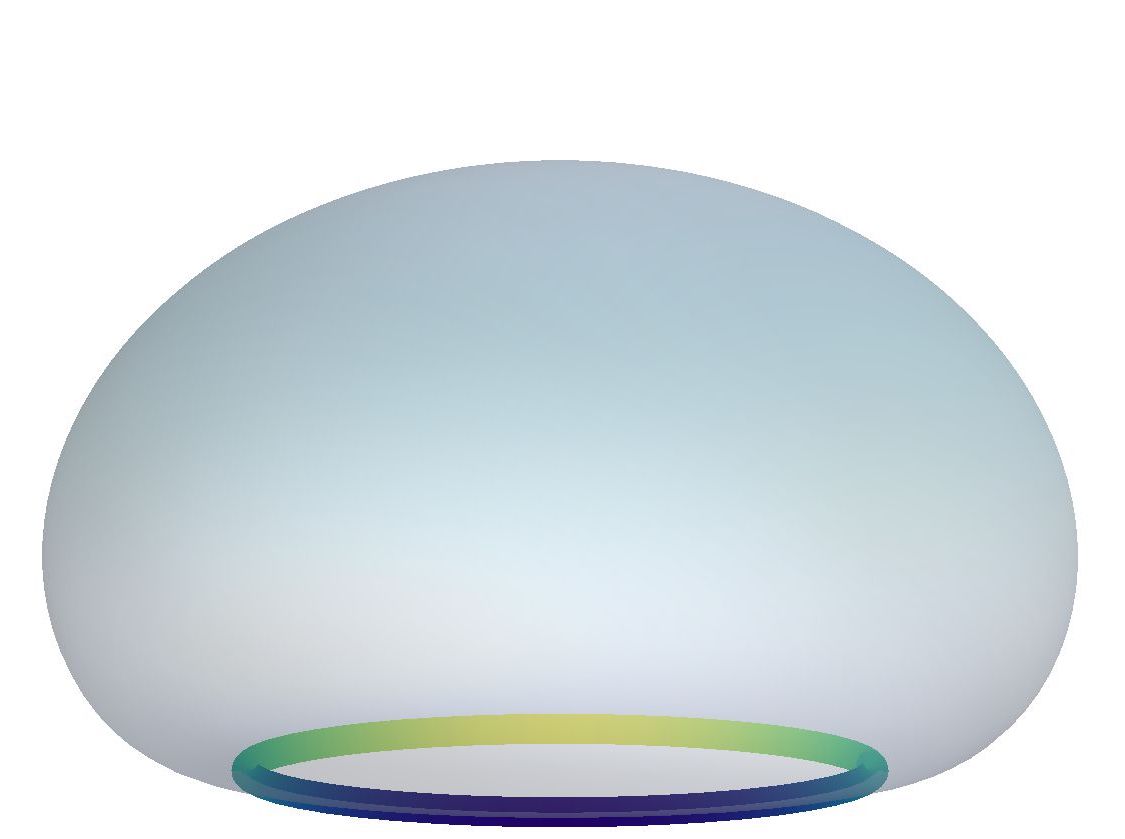}
\caption{$E_{1,2,2,1,4}$}
\end{subfigure}
\quad
\begin{subfigure}[b]{0.35\linewidth}
\includegraphics[width=\linewidth]{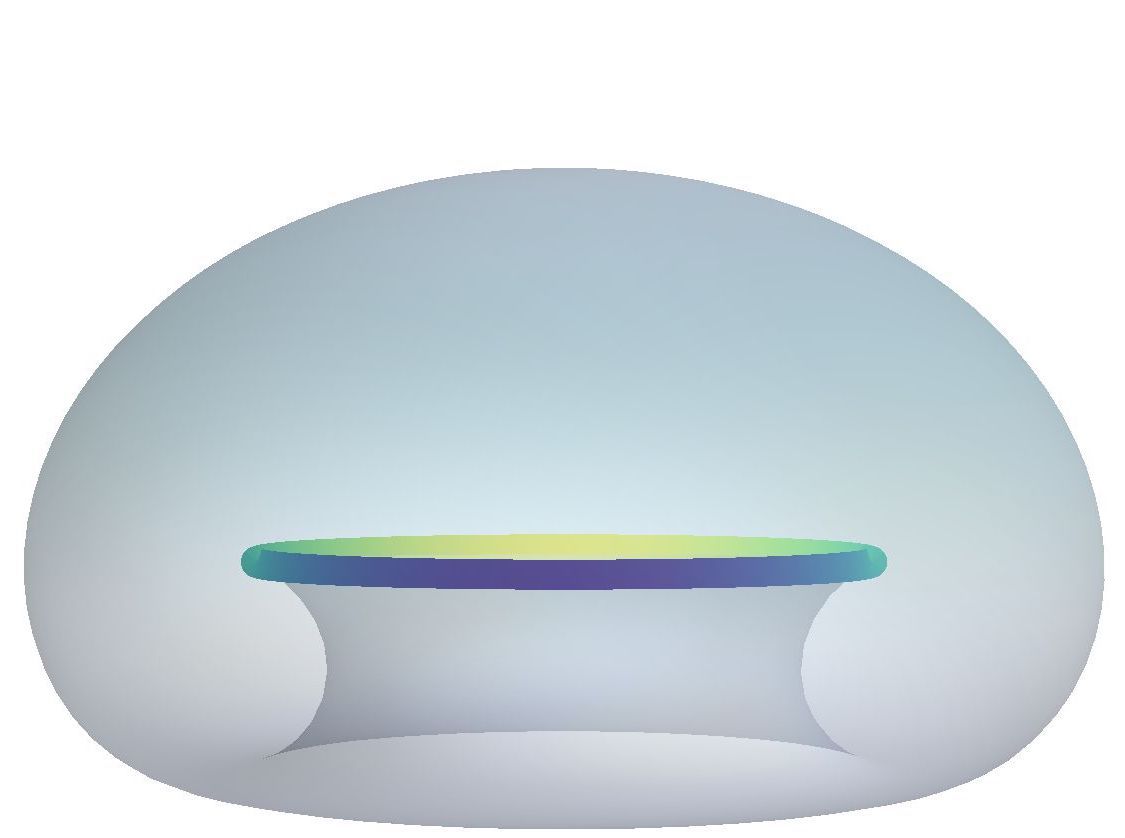}
\caption{$E_{1,2,-0.5,1,4}$}
\end{subfigure}}
\caption{Axially symmetric critical discs for the energy $E_{a,c_o,b,\alpha,\beta}$ bounded by a non geodesic (parallel) circle of radius $r\neq\sqrt{\alpha/\beta}$.}
\label{newcase}
\end{figure}

Consequently, as above, we solve the system \eqref{dif1}-\eqref{dif3} with the initial conditions $r(0)=0$, $z(0)=z_o\neq 0$, $\varphi(0)=0$ and the boundary conditions \eqref{bcnew1} and \eqref{bcnew2}, obtaining several examples of axially symmetric critical domains. We show some of these axially symmetric critical domains for the energy $E$ with $c_o>0$ in Figure \ref{newcase}.

\begin{figure}[h!]
\makebox[\textwidth][c]{
\hspace{0.5cm}
\begin{subfigure}[b]{0.3\linewidth}
\includegraphics[width=\linewidth]{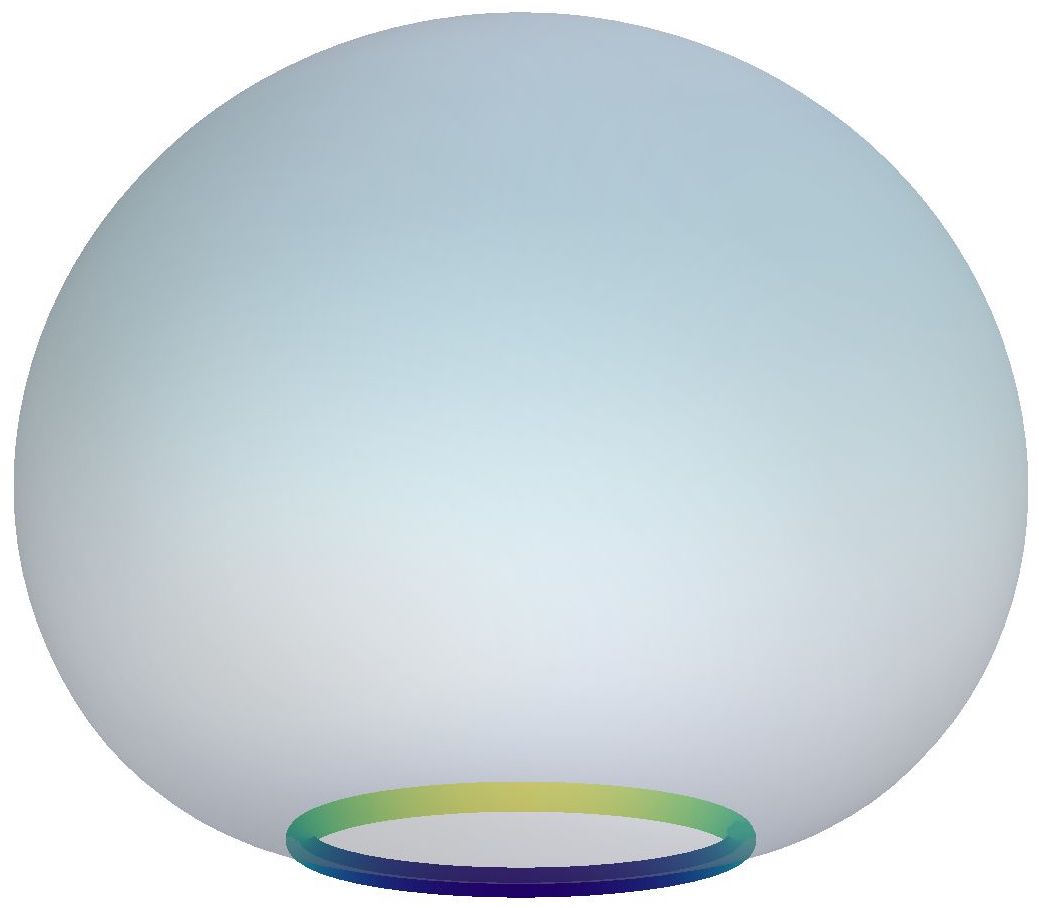}
\caption{$c_o=0.5$}
\end{subfigure}
\quad\,\,
\begin{subfigure}[b]{0.3\linewidth}
\includegraphics[width=\linewidth]{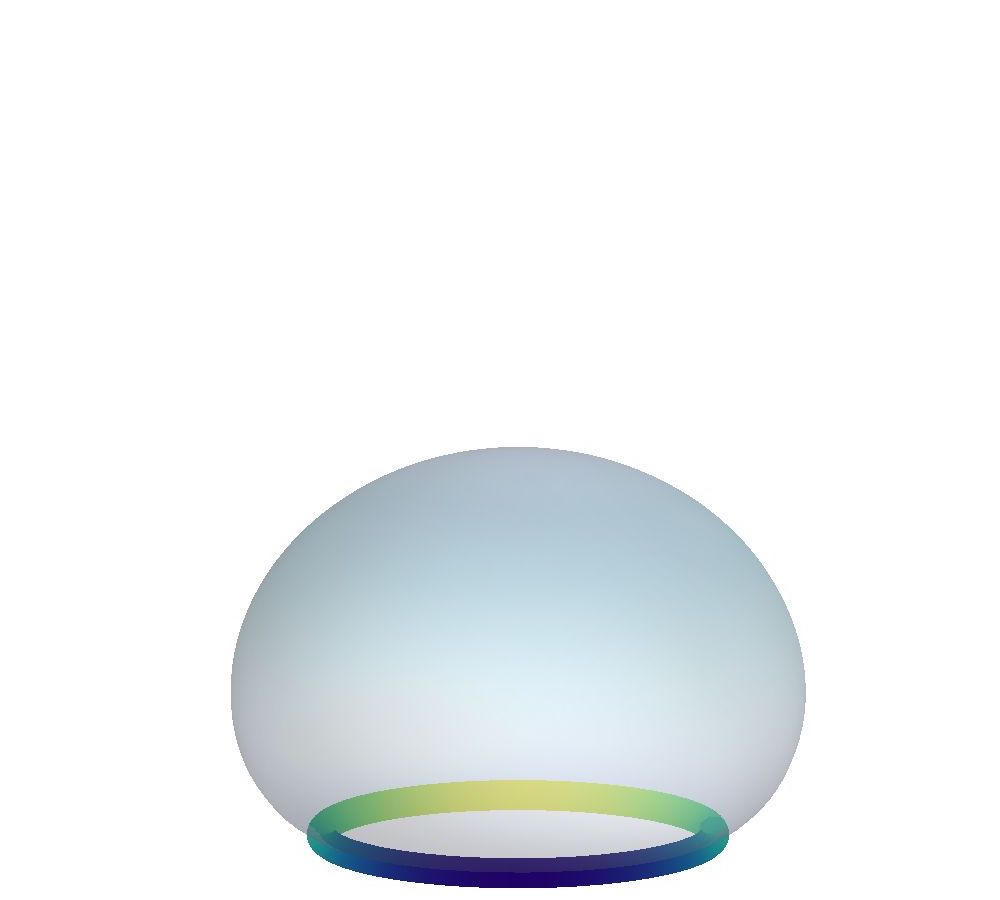}
\caption{$c_o=1$}
\end{subfigure}
\hspace{-0.6cm}
\begin{subfigure}[b]{0.3\linewidth}
\includegraphics[width=\linewidth]{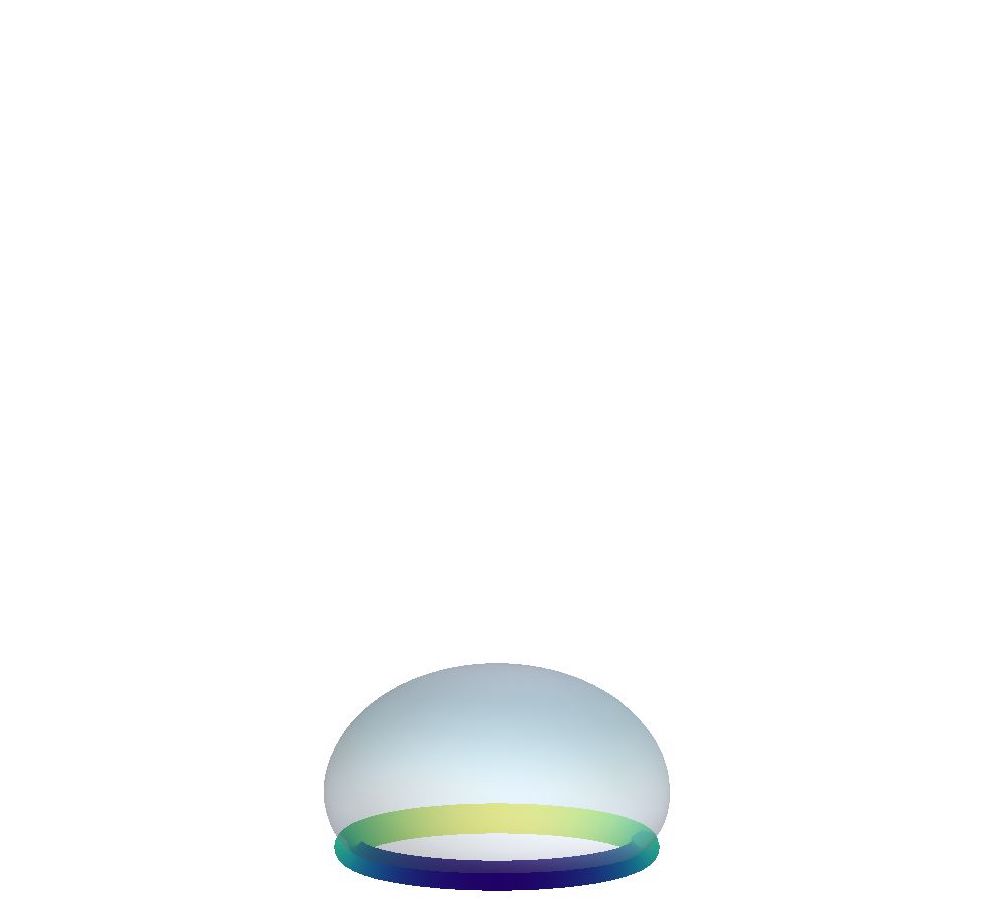}
\caption{$c_o=2$}
\end{subfigure}
\hspace{-1.3cm}
\begin{subfigure}[b]{0.3\linewidth}
\includegraphics[width=\linewidth]{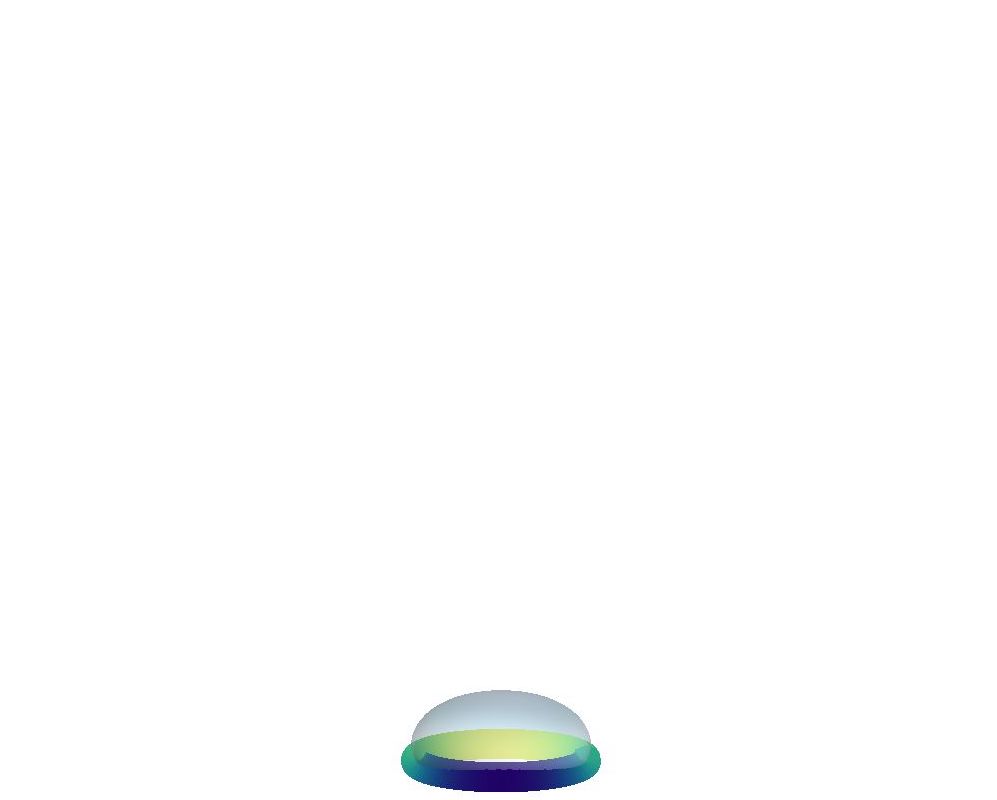}
\caption{$c_o=5$}
\end{subfigure}}
\caption{Axially symmetric critical discs for the energy $E_{a=1,c_o,b=-0.5,\alpha=1,\beta=1}$ for different values of $c_o>0$. In all the cases the boundary is a non geodesic circle of radius $r<\sqrt{\alpha/\beta}$.}
\label{b<0}
\end{figure}

\begin{figure}[h!]
\makebox[\textwidth][c]{
\hspace{0.7cm}
\begin{subfigure}[b]{0.3\linewidth}
\includegraphics[width=\linewidth]{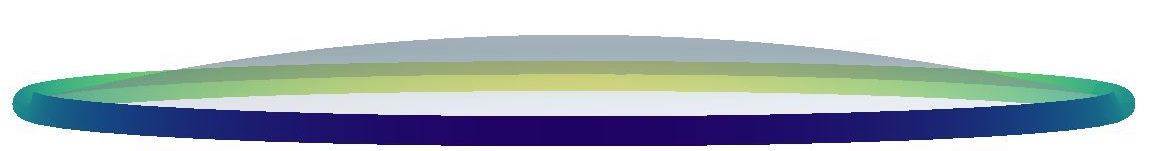}
\caption{$c_o=0.5$}
\end{subfigure}
\quad\quad
\begin{subfigure}[b]{0.3\linewidth}
\includegraphics[width=\linewidth]{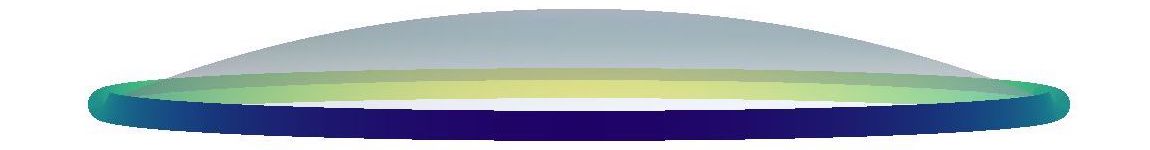}
\caption{$c_o=1$}
\end{subfigure}
\,
\begin{subfigure}[b]{0.3\linewidth}
\includegraphics[width=\linewidth]{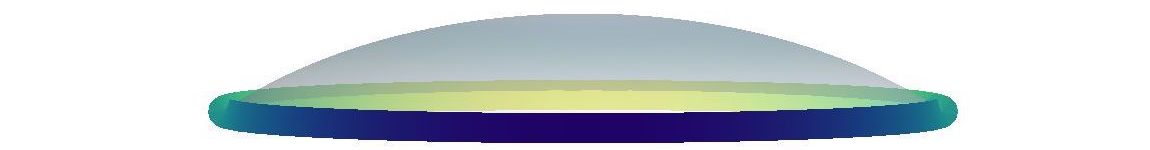}
\caption{$c_o=2$}
\end{subfigure}
\hspace{-0.8cm}
\begin{subfigure}[b]{0.3\linewidth}
\includegraphics[width=\linewidth]{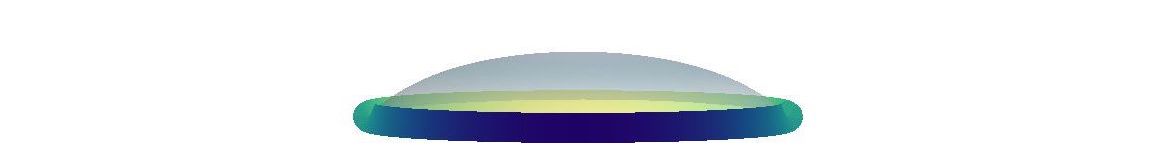}
\caption{$c_o=5$}
\end{subfigure}}
\caption{Axially symmetric critical discs for the energy $E_{a=1,c_o,b=1,\alpha=1,\beta=1}$ for different values of $c_o>0$. In all the cases the boundary is a non geodesic circle of radius $r<\sqrt{\alpha/\beta}$.}
\label{b>0}
\end{figure}

As in the case $b=0$ (Figure \ref{b=0}), for fixed $a>0$, $\alpha>0$ and $\beta>0$ we show in Figure \ref{b<0} ($b<0$ fixed) and Figure \ref{b>0} ($b>0$ fixed) some axially symmetric discs with non constant mean curvature bounded by non geodesic circles of radii $r\neq\sqrt{\alpha/\beta}$ which are critical for the energy $E$ with different values of $c_o$. As the spontaneous curvature $c_o$ increases, these critical domains get smaller in an attempt to decrease the energy. In this process the radius of the boundary circle also decreases.

\section{Stability and Minimization Results}

In this section we investigate when axially symmetric critical discs can minimize the energy $E$. For this purpose, we first recall that on Theorem 4.3 of \cite{PP2} all the cases where the spontaneous curvature is zero were discussed. In particular, it was proved that \emph{when the energy $E$ for $c_o=0$ has a finite infimum, this infimum is either attained or approached by a sequence of totally umbilical surfaces bounded by a circle of radius $\sqrt{\alpha/\beta}$}.

Consequently, we assume from now on that $c_o\neq 0$ holds. In this case, as proved in Lemma 4.2 of \cite{PP2} (see also Remark 4.1 of \cite{PP2}), \emph{the infimum of the energy $E$ is finite if and only if the quantity ${\underline E}:=2\sqrt{\alpha\,\beta}-\lvert b\rvert$ is nonnegative}. In what follows, we also assume that ${\underline E}\geq 0$ holds.

Under these restrictions, we begin by proving a result regarding instability of axially symmetric critical domains for the energy $E$. We do not derive a general second variation formula but instead rely on variations of the surface obtained from simply varying a boundary curve within a fixed surface satisfying \eqref{EL1}  and seeing how the energy changes to second order. As one would expect, this can only yield partial results. 

\begin{prop}\label{instability} Let $X:\Sigma\rightarrow{\bf R}^3$ be an axially symmetric immersion critical for $E_{a,c_o,b,\alpha,\beta}$ and assume one of the following cases is satisfied: 
\begin{enumerate}[(i)]
\item Case $b=0$ and any boundary component is a geodesic circle along which
$$\hspace{1cm}\left(\frac{\alpha}{r^2}-\beta\right)K<0$$
holds.
\item Case $b<0$ and any boundary component is a geodesic circle of radius $r>\sqrt{\alpha/\beta}$.
\item Case $b>0$ and any boundary component is a geodesic circle of radius $r<\sqrt{\alpha/\beta}$.
\item Case
$$\hspace{1cm}\alpha>\frac{\beta^3(a+b)^2}{4c_o^2\left(b^2\left[a+b\right]^2c_o^2+a^2\beta^2\right)}$$
and any boundary component is a circle of radius $\sqrt{\alpha/\beta}$ along which $b\kappa_g>0$ holds.
\end{enumerate}
Then the surface is unstable. %Moreover, in cases (ii) and (iii), if $r=\sqrt{\alpha/\beta}$ the surface is not minimizing.
\end{prop}
{\it Proof.\:} We prove the result by considering variations of the axially symmetric critical immersion $X:\Sigma\rightarrow{\bf R}^3$ which fix the surface while vary one boundary component through parallels.

For $\lvert \epsilon\rvert\geq 0$ sufficiently small, we denote $\Sigma_\epsilon=\left[0,\mathcal{L}-\epsilon\right]\times\left[0,2\pi\right)$ the abstract surfaces, where $\mathcal{L}$ denotes the length of the profile curve $\gamma(\sigma)$ for the critical immersion $X:\Sigma\rightarrow{\bf R}^3$. The energy $E$ of these domains $\Sigma_\epsilon$ is then given by
$$f(\epsilon):=E[\Sigma_\epsilon]=\int_{\Sigma_\epsilon} \left(a\left[H+c_o\right]^2+b K\right)d\Sigma+\oint_{\partial\Sigma_\epsilon}\left(\alpha\kappa^2+\beta\right)ds\,.$$
Now, using that the surface is axially symmetric and coordinates $(\sigma,s)\in\Sigma_\epsilon$, we rewrite the energy as
$$f(\epsilon)=2\pi\int_0^{\mathcal{L}-\epsilon}\left(a\left[H(\sigma)+c_o\right]^2+bK(\sigma)\right)r(\sigma)\,d\sigma+2\pi\left(\frac{\alpha}{r^2}+\beta\right) r+P\,,$$
where $P$ is a constant representing the integral $\int(\alpha\kappa^2+\beta)ds$ over the fixed boundary component. In the case of a topological disc, $P=0$. All functions are evaluated at $\sigma=\mathcal{L}-\epsilon$. For simplicity, we are avoiding to explicitly write it. 

We next proceed differentiating $f$ with respect to the variation parameter $\epsilon$ (we denote by upper ``dots" derivatives with respect to $\epsilon$ and by $\left(\,\right)'$ the derivatives with respect to the arc length parameter $\sigma$). For the first derivative we have,
$$\dot{f}(\epsilon):=\frac{df}{d\epsilon}(\epsilon)=-2\pi \left(a\left[H+c_o\right]^2+bK\right)r+2\pi\left(\frac{\alpha}{r^2}-\beta\right)r'\,.$$

From the boundary condition \eqref{EL4} it follows that $\dot{f}(0)=0$, as expected from the criticality of the immersion.

The second derivative of $f$ with respect to $\epsilon$ evaluated at $\epsilon=0$, after using once again \eqref{EL4} and the definitions $K:=-r''/r$ and $\kappa_g:=-r'/r$, is given by
\begin{equation*}
\ddot{f}(0)=2\pi r\left(2a \left[H+c_o\right]H'+bK'+\left[\alpha r^{-2}-\beta\right] K+\left[3\alpha r^{-2}-\beta\right]\kappa_g^2\right).
\end{equation*}

We now compute $K'$ along the boundary. To this end, denote $S^2=H^2-K$ the skew curvature of the surface. Differentiating it with respect to $\varsigma$ and combining with the Codazzi equation $(r^2S)'=r^2H'$, we obtain at $\varsigma=\mathcal{L}$,
$$K'=2\left(H-S\right)H'-4\kappa_g\left(H^2-K\right).$$

Assume first that the boundary is a circle of radius $\sqrt{\alpha/\beta}$. Thus, the third term in $\ddot{f}(0)$ vanishes. Moreover, from \eqref{EL3}, we also get that $H'=0$ holds along the boundary. In this case, $S^2=H^2-K=c_o^2$ holds along the boundary, so that
$$K'=-4c_o^2\kappa_g\,.$$
Consequently, it is then clear that if 
$$\ddot{f}(0)=4\pi r\left(\beta\kappa_g^2-2bc_o^2\kappa_g\right)<0\,,$$
then the surface is unstable. Using the boundary condition \eqref{bcdif2} we obtain
$$\lvert r'(\mathcal{L})\rvert=\lvert\nu_3(\mathcal{L})\rvert=\lvert \cos\varphi(\mathcal{L})\rvert=\sqrt{1-\frac{4a^2c_o^2\alpha}{(a+b)^2\beta}}\,,$$
from which we compute $\kappa_g=-r'/r$. This proves case (iv).

Now, if the boundary circle is a geodesic, we use the boundary conditions \eqref{EL2}-\eqref{EL4} to conclude that
$$\ddot{f}(0)=2\pi \left(\alpha r^{-2}-\beta\right) K r\,.$$
Therefore, if $\left(\alpha r^{-2}-\beta\right) K<0$ holds along the boundary geodesic, the surface is unstable. This proves case (i). If $b\neq 0$, we use \eqref{EL4} to rewrite $K$ in terms of $a$, $b$ and $\left(H+c_o\right)^2$, concluding with cases (ii) and (iii). {\bf q.e.d.}
\\ 

From parts (ii) and (iii) of Proposition \ref{instability}, we conclude that the surfaces shown in Figure \ref{geo}, (A) and (D), are unstable.

The result of Proposition \ref{instability} for geodesic boundary circles can be rewritten in terms of geometric objects of the profile curve $\gamma$.

\begin{cor} Let $X:\Sigma\rightarrow{\bf R}^3$ be an axially symmetric disc type surface critical for $E$ and assume that it is stable. Denote by $\gamma(\sigma)$, $\sigma\in\left[0,\mathcal{L}\right]$, the profile curve of the surface. If $\gamma'(\mathcal{L})$ is vertical, then the squared curvature $\mu^2$ of $\gamma$ is non-decreasing at $\sigma=\mathcal{L}$.
\end{cor}
{\it Proof.\:} The proof follows immediately after noticing that $r''/z'$ is the curvature of the profile curve $\gamma$, up to perhaps the sign. In any case, using equation \eqref{EL3}, we have that
$$-\left(\alpha\kappa^2-\beta\right)r''=\frac{a}{2} rr''r'''=\frac{a}{2} r\mu\mu'=\frac{a}{4} r(\mu^2)'$$
holds. In the calculation above, we have used that $z''=\pm \mu r'=0$ since the tangent vector is vertical. Therefore, from Proposition \ref{instability}, (i)-(iii), we conclude that for a stable domain, the quantity appearing above must be nonnegative. {\bf q.e.d.}
\\

Since our second variation formula only involves deforming the surface through subdomains, we make the following definition. We define a {\it minimizing domain} to be a smoothly bounded subdomain of a surface $\Sigma$ which minimizes the energy $E_{a, c_o, b, \alpha, \beta}$ among all smoothly bounded subdomains of $\Sigma$. In this definition, it will always be assumed that $\Sigma$ satisfies \eqref{EL1}. Of course any minimizer for the functional must be a minimizing domain.

We conclude this section with results restricting the geometry of minimizing domains. 

For vanishing saddle-splay modulus $b$, we will see that whenever a constant mean curvature surface does not attain the minimum, i.e. when $c_o^2>\beta/\alpha$ the lower bound $E_D:=4\pi\sqrt{\alpha\,\beta}$ (see Lemma 4.2 of \cite{PP2}), then a minimizing domain of an axially symmetric surface must be bounded by a geodesic circle of radius $r<\sqrt{\alpha/\beta}$.

\begin{theorem}\label{nonminb=0} Let $X:\Sigma\rightarrow{\bf R}^3$ be an axially symmetric disc type surface and assume that the surface is a minimizing domain for an energy $E$ with $b=0$. Then:
\begin{enumerate}[(i)]
\item If $c_o=0$ holds, the surface is a planar disc bounded by a circle of radius $\sqrt{\alpha/\beta}$.
\item If $0<c_o^2\leq\beta/\alpha$ holds, the surface is a spherical cap with $H\equiv -c_o$ bounded by a circle of radius $\sqrt{\alpha/\beta}$.
\item If  $c_o^2>\beta/\alpha$ holds, the mean curvature of the surface satisfies
$$\hspace{1cm} H+c_o=-\frac{\nu_3}{z}\,$$
and the boundary is a geodesic circle of radius $r<\sqrt{\alpha/\beta}$.
\end{enumerate}
\end{theorem}
{\it Proof.\:} For the choice of the saddle-splay modulus $b=0$, the condition ${\underline E}\geq 0$ is clearly satisfied and, hence, the infimum of the energy $E$ is finite. From Theorem 4.3 of \cite{PP2}, \emph{the minimum of the energy $E$ among any sufficiently smooth topological disc (whether it is axially symmetric or not) is attained by a planar disc bounded by a circle of radius $\sqrt{\alpha/\beta}$}. In particular, we obtain (i). Moreover, from Proposition 4.1 of \cite{PP2}, we have that \emph{spherical caps with constant mean curvature $H=-c_o$ bounded by a circle of radius $\sqrt{\alpha/\beta}$ are the absolute minimizers of $E$ with $b=0$ among any topological discs, provided that $c_o^2\leq\beta/\alpha$}. This proves the case (ii). If $c_o^2>\beta/\alpha$ then no such spherical cap exists and, from Proposition \ref{b=0prop}, the surface is bounded by a closed geodesic of radius $r\neq\sqrt{\alpha/\beta}$.

In what follows, by contradiction, we assume that the surface is bounded by a geodesic circle of radius $r>\sqrt{\alpha/\beta}$. Denote $\Sigma=[0,r]\times [0,2\pi)$ the abstract surface. Since $r>\sqrt{\alpha/\beta}$, there exists a subdomain $\widetilde{\Sigma}=[0,\sqrt{\alpha/\beta}]\times [0,2\pi)\subset\Sigma$. Then, since the boundary energy attains its minimum at a circle of radius $\sqrt{\alpha/\beta}$ (see Lemma 4.1 of \cite{PP2}), we have
\begin{eqnarray*}
E[\Sigma]&=&a\int_\Sigma\left(H+c_o\right)^2 d\Sigma+\oint_{\partial\Sigma}\left(\alpha\kappa^2+\beta\right)ds\\
&>&a\int_{\Sigma\setminus\widetilde{\Sigma}}\left(H+c_o\right)^2 d\Sigma+a\int_{\widetilde{\Sigma}}\left(H+c_o\right)^2 d\Sigma+4\pi\sqrt{\alpha\,\beta}\\
&\geq&a\int_{\widetilde{\Sigma}}\left(H+c_o\right)^2 d\Sigma+4\pi\sqrt{\alpha\,\beta}=E[\widetilde{\Sigma}]\,.
\end{eqnarray*}
This proves that an axially symmetric surface bounded by a circle of radius $r>\sqrt{\alpha/\beta}$ cannot be a minimizing domain of $E$ with $b=0$ among topological discs, concluding with (iii). {\bf q.e.d.}
\\

For the case $b\neq 0$, it seems reasonable to think that surfaces minimizing the energy among axially symmetric ones are bounded by a circle of radius $\sqrt{\alpha/\beta}$ since these circles attain the minimum of the boundary energy (see Lemma 4.1 of \cite{PP2}). However, in the following result we obtain a negative answer to the question regarding the existence of minimizing surfaces bounded by circles of radii $\sqrt{\alpha/\beta}$ for many choices of the energy parameters.

\begin{prop}\label{nonminr} Let $X:\Sigma\rightarrow{\bf R}^3$ be an axially symmetric disc type surface which is critical for $E$ with $-3a\leq b \leq a$. If the boundary is a circle of radius $\sqrt{\alpha/\beta}$, then the surface does not minimize the energy $E$ among axially symmetric surfaces.
\end{prop}
{\it Proof.\:} Assume that $-3a\leq b\leq a$ (see Figure \ref{bal}, (C)-(F)). Then, $(a+b)^2\leq 4a^2$ holds. Moreover, from \eqref{radiusnongeodesic}, we obtain that $c_o^2\leq \beta/\alpha$. In these cases, there exists a sphere with constant mean curvature $H=-c_o$ such that it contains circles of radii $\sqrt{\alpha/\beta}$. Any of those circles divides the sphere in two caps, a ``big" spherical cap (i.e. more than half sphere) and a ``small" spherical cap (i.e. less than half sphere). We denote them by $\Omega_B$ and $\Omega_S$, respectively. If equality holds both domains represent a half sphere.

Next, for the axially symmetric disc $X:\Sigma\rightarrow{\bf R}^3$ we compute its energy. Since the boundary is a circle of radius $\sqrt{\alpha/\beta}$, the boundary term of the energy $E$ gives $4\pi\sqrt{\alpha\,\beta}$ and, hence,
\begin{eqnarray*}
E[\Sigma]&=&a\int_\Sigma \left(H+c_o\right)^2 d\Sigma+b\int_\Sigma K\,d\Sigma+4\pi\sqrt{\alpha\,\beta}\\
&>&2\pi b-2\pi b r'(\mathcal{L})+4\pi\sqrt{\alpha\,\beta}\\
&\geq&2\pi b-2\pi\lvert b\rvert \lvert r'(\mathcal{L})\rvert +4\pi\sqrt{\alpha\,\beta}\,,
\end{eqnarray*}
where, in the first line, we have used that $(H+c_o)^2>0$ holds for these axially symmetric surfaces.

As in the proof of Proposition \ref{instability}, using the boundary condition \eqref{bcdif2} we compute $\lvert r'(\mathcal{L})\rvert=\lvert\nu_3(\mathcal{L})\rvert$. Moreover, from $(a+b)^2\leq 4a^2$ we conclude that
$$\lvert r'(\mathcal{L})\rvert\leq \sqrt{1-c_o^2\frac{\alpha}{\beta}}\,,$$
and, hence,
$$E[\Sigma]>2\pi b-2\pi \lvert b\rvert \sqrt{1-c_o^2\frac{\alpha}{\beta}}+4\pi\sqrt{\alpha\,\beta}\,.$$
The right hand side of above equation is the energy of the ``big" spherical cap $\Omega_B$ if $b<0$; or if $b>0$ it represents the energy of the ``small" spherical cap $\Omega_S$. {\bf q.e.d.}

\section*{Appendix A. Axially Symmetric Annuli}

In this appendix we will show the existence of axially symmetric domains for which \eqref{nonCMC} does not hold, i.e. we will produce annular surfaces for which the relation \eqref{fir} holds for $\bar{A}\neq 0$. We utilize the notation of Section \ref{Examples}.

With this notation and introducing a new variable $\zeta:=\varphi'$, equation \eqref{fir} can be expressed as
\begin{equation}\label{zprime}
\zeta'=2(H+c_o)(H+c_o-\zeta)\tan\varphi+2(H-\zeta)\frac{\cos\varphi}{r}-\frac{2\bar{A}}{r\cos\varphi}\,,
\end{equation}
where $H$ is given in \eqref{HAS}.

We impose initial conditions for $r$, $z$, $\varphi$ and $\zeta$ so that the right hand side of \eqref{zprime} is defined, then standard existence and uniqueness theorems for first order systems of differential equations ensure the existence of a short time solution for the system consisting of \eqref{dif1}, \eqref{dif2}, $\varphi'=\zeta$ and \eqref{zprime}, from which the desired surface can be produced. The resulting surfaces cannot be continued to give smooth discs since, otherwise, Theorem \ref{discsnew} would give a contradiction. 

Examples of axially symmetric annular surfaces obtained this way are domains in circular biconcave discoids (\cite{NOOY}).

\section*{Appendix B. Numerical Calculations of Energies}

\begin{center}
\begin{table}[h!]
\caption{Numerical values of the energy $E$ for the critical domains of Figures \ref{b=0}-\ref{b>0}.}\label{tab1}
\makebox[\textwidth][c]{
\begin{tabular}{|c|c|c|c|c|} 
\hline
Surface & Parameters ($E=E_{a,c_o,b,\alpha,\beta}$) & Radius & Energy & $E_D$ \\
\hline
Fig. \ref{b=0}, (A) & $a=1, c_o=1.1, b=0, \alpha=1, \beta=1$ & 0.96 & 12.6 & 12.57 \\ 
Fig. \ref{b=0}, (B) & $a=1, c_o=2, b=0, \alpha=1, \beta=1$ & 0.7 & 14.1 & 12.57 \\
Fig. \ref{b=0}, (C) & $a=1, c_o=5, b=0, \alpha=1, \beta=1$ & 0.38 & 22.2 & 12.57 \\
Fig. \ref{geo}, (A) & $a=1, c_o=2, b=0.25, \alpha=1, \beta=10$ & 0.31 & 42.2 & 39.74 \\
Fig. \ref{geo}, (B) & $a=1, c_o=2, b=-0.11, \alpha=1, \beta=20$ & 0.22 & 55.9 & 54.82 \\
Fig. \ref{geo}, (C) & $a=1, c_o=2, b=0.08, \alpha=1, \beta=14$ & 0.27 & 47.9 & 47.02 \\
Fig. \ref{geo}, (D) & $a=1, c_o=2, b=-0.05, \alpha=1, \beta=18$ & 0.24 & 53.2 & 52.69 \\
Fig. \ref{bal}, (A) & $a=1, c_o=2, b=-5, \alpha=1, \beta=20$ & 0.22 & 64.8 & -6.63 \\
Fig. \ref{bal}, (B) & $a=1, c_o=2, b=2, \alpha=1, \beta=10$ & 0.32 & 45 & 39.74 \\
Fig. \ref{bal}, (C) & $a=1, c_o=2, b=0.05, \alpha=1, \beta=31$ & 0.18 & 71.3 & 69.97 \\
Fig. \ref{bal}, (D) & $a=1, c_o=2, b=-0.05, \alpha=1, \beta=33$ & 0.17 & 72.5 & 71.56 \\
Fig. \ref{bal}, (E) & $a=1, c_o=2, b=0.05, \alpha=1, \beta=18$ & 0.24 & 54.4 & 53.31 \\
Fig. \ref{bal}, (F) & $a=1, c_o=2, b=-0.05, \alpha=1, \beta=20$ & 0.22 & 55.8 & 55.57 \\
Fig. \ref{newcase}, (A) & $a=1, c_o=2, b=0.5, \alpha=1, \beta=4$ & 0.46 & 26.2 & 25.13 \\
Fig. \ref{newcase}, (B) & $a=1, c_o=2, b=-0.5, \alpha=1, \beta=4$ & 0.46 & 19.6 & 18.85 \\
Fig. \ref{newcase}, (C) & $a=1, c_o=2, b=2, \alpha=1, \beta=4$ & 0.44 & 51.4 & 25.13 \\
Fig. \ref{newcase}, (D) & $a=1, c_o=2, b=-0.5, \alpha=1, \beta=4$ & 0.43 & 26.4 & 18.85 \\
Fig. \ref{b<0}, (A) & $a=1, c_o=0.5, b=-0.5, \alpha=1, \beta=1$ & 0.98 & 6.4 & 6.28 \\
Fig. \ref{b<0}, (B) & $a=1, c_o=1, b=-0.5, \alpha=1, \beta=1$ & 0.87 & 7 & 6.28 \\
Fig. \ref{b<0}, (C) & $a=1, c_o=2, b=-0.5, \alpha=1, \beta=1$ & 0.65 & 9.4 & 6.28 \\
Fig. \ref{b<0}, (D) & $a=1, c_o=5, b=-0.5, \alpha=1, \beta=1$ & 0.38 & 18.1 & 6.28 \\
Fig. \ref{b>0}, (A) & $a=1, c_o=0.5, b=1, \alpha=1, \beta=1$ & 0.94 & 12.9 & 12.57 \\
Fig. \ref{b>0}, (B) & $a=1, c_o=1, b=1, \alpha=1, \beta=1$ & 0.83 & 13.9 & 12.57 \\
Fig. \ref{b>0}, (C) & $a=1, c_o=2, b=1, \alpha=1, \beta=1$ & 0.63 & 16.7 & 12.57 \\
Fig. \ref{b>0}, (D) & $a=1, c_o=5, b=1, \alpha=1, \beta=1$ & 0.37 & 25.7 &12.57 \\
\hline
\end{tabular}}
\end{table}
\end{center}

In Table \ref{tab1} we show the numerical values of the energy $E$ and the radii of the boundary circles for all the critical domains illustrated throughout the paper. We also include the lower bound obtained in Lemma 4.2 of \cite{PP2} for topological discs,
\begin{equation}\label{ED}
E_D:=2\pi\left(2\sqrt{\alpha\,\beta}-\lvert b\rvert+b\right),
\end{equation}
in order to compare to the energy $E$ of the critical domains. This bound can only be attained if $H+c_o\equiv 0$ holds and the boundary is a circle of radius $\sqrt{\alpha/\beta}$. Moreover, either $b=0$ or $\kappa_n\equiv 0$ must hold along the boundary circle. The case $c_o=0$ was studied in Theorem 4.3 of \cite{PP2}. For positive spontaneous curvature, $c_o>0$, if a surface attains the lower bound \eqref{ED} then it must be a spherical cap of radius $c_o^{-1}$ bounded by a circle of radius $\sqrt{\alpha/\beta}$. If, in addition, $b=0$, these caps exist provided that $c_o^2\leq\beta/\alpha$ (see Proposition 4.1 of \cite{PP2}). On the other hand, if $b\neq 0$ then $\kappa_n\equiv 0$ must hold along the boundary circle. These circles do not exist in spherical caps and, hence, attaining the lower bound \eqref{ED} is impossible for $c_o>0$ and $b\neq 0$.

Observe that the computed energies of many examples appear remarkably close to the lower bound \eqref{ED}. A priori, this may suggest that these surfaces attain the minimum value. However, some of them are unstable (see Proposition \ref{instability}), as is the case of Figures \ref{geo} (A) and (D), and others are not minimizing, as for instance Figures \ref{bal} (C)-(F) (see Proposition \ref{nonminr}). In all these cases, for the same energy parameters, we are able to numerically find another critical surface with less energy and whose shape is similar to those shown in Figures \ref{b<0} and \ref{b>0}, for $b<0$ and $b>0$, respectively. For instance, for the fixed parameters of Figure \ref{bal} (C), we found a critical disc with $E=69.99$. Similarly, for the energy parameters of Figure \ref{bal} (D), we obtained a critical domain with energy $E=71.57$. In both cases, the energy of these domains is much closer to the lower bound $E_D$, \eqref{ED}. 

In light of this, we show in Table \ref{tab2} the conjectured minimizing critical discs for the energy $E$. For fixed energy parameters $a>0$, $\alpha>0$ and $\beta>0$, we show the deformation of these discs as the spontaneous curvature $c_o>0$ and the saddle-splay modulus $b$ vary. The first column of Table \ref{tab2}, i.e. $b=-0.5$, is composed precisely by the domains of Figure \ref{b<0}, while the last column ($b=1$) contains the domains of Figure \ref{b>0}. For the case $b=0$ (the middle column), the first domain ($c_o=0.5$) is a ``small" spherical cap of radius $c_o^{-1}$ bounded by a circle of radius $\sqrt{\alpha/\beta}$, which was proven in Proposition 4.1 of \cite{PP2} to be the absolute minimizer for these energy parameters (the rest of the sphere, i.e. the ``big" spherical cap, has the same energy and so is also an absolute minimizer). The second domain is a half sphere of radius $c_o^{-1}$ bounded by the equator, which is a geodesic circle of radius $\sqrt{\alpha/\beta}$. This domain is also an absolute minimizer (Proposition 4.1 of \cite{PP2}) and it represents the limiting domain for the condition $c_o^2\leq \beta/\alpha$. When $c_o$ is bigger no such spherical cap can exist and, hence, for the conjectured minimizing critical discs in the remaining cases we show the domains of Figure \ref{b=0}, (B) and (C), respectively.

\begin{center}
\begin{table}[H]
\caption{Conjectured minimizing discs for the energy $E_{a=1,c_o,b,\alpha=1,\beta=1}$ as the spontaneous curvature $c_o>0$ and the saddle-splay modulus $b$ increase.} \label{tab2}
\makebox[\textwidth][c]{
\def\arraystretch{2}
\begin{tabular}{| c | c | c | c |} 
\hline
& $b=-0.5$ & $b=0$ & $b=1$ \\
\hline
\rule{0pt}{4.3cm}\raisebox{7.3\height}{$c_o=0.5$} & \,\,\raisebox{-0.045\height}{\includegraphics[width=0.35\linewidth]{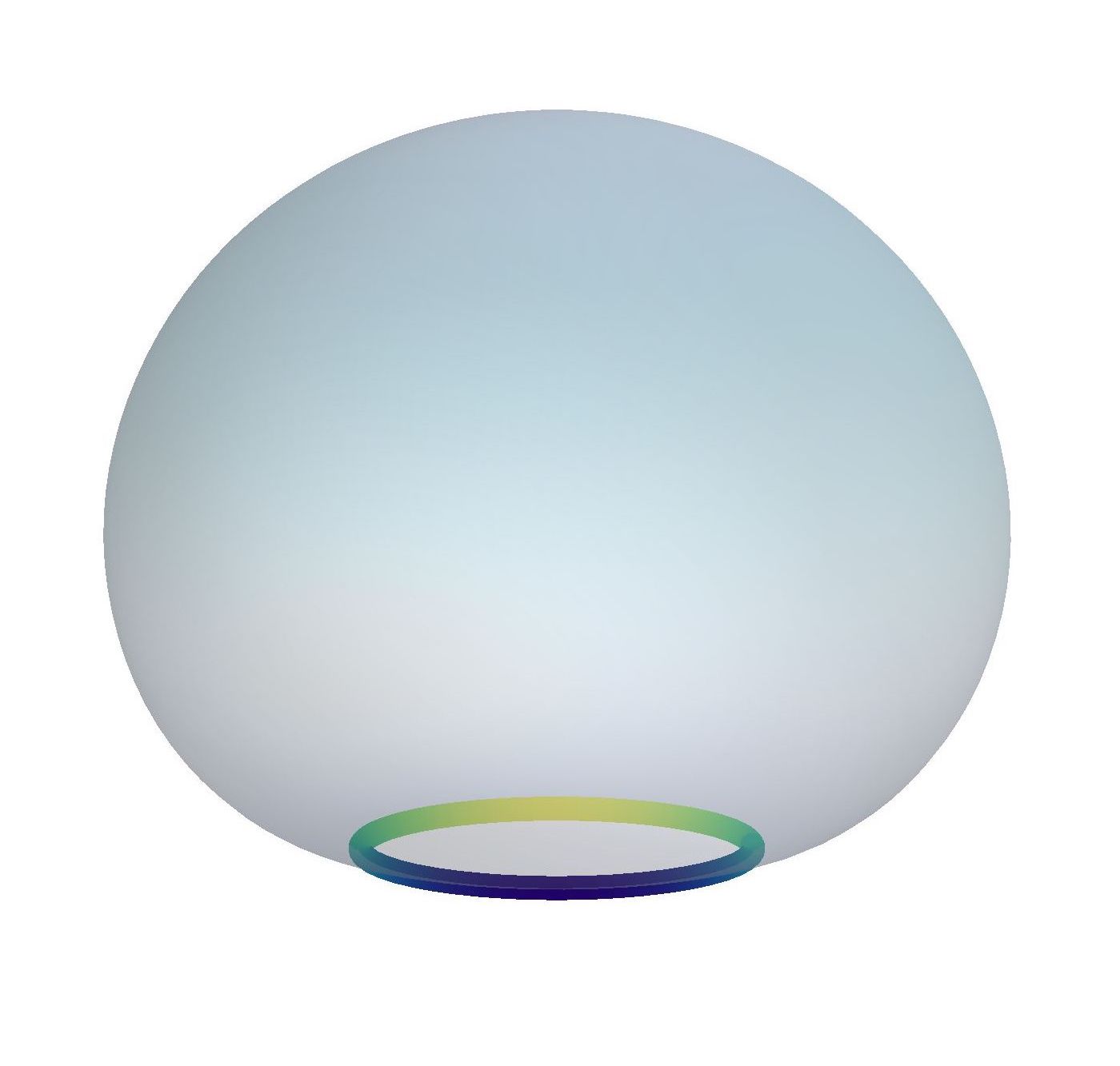}} & \raisebox{0.06\height}{\includegraphics[width=0.35\linewidth]{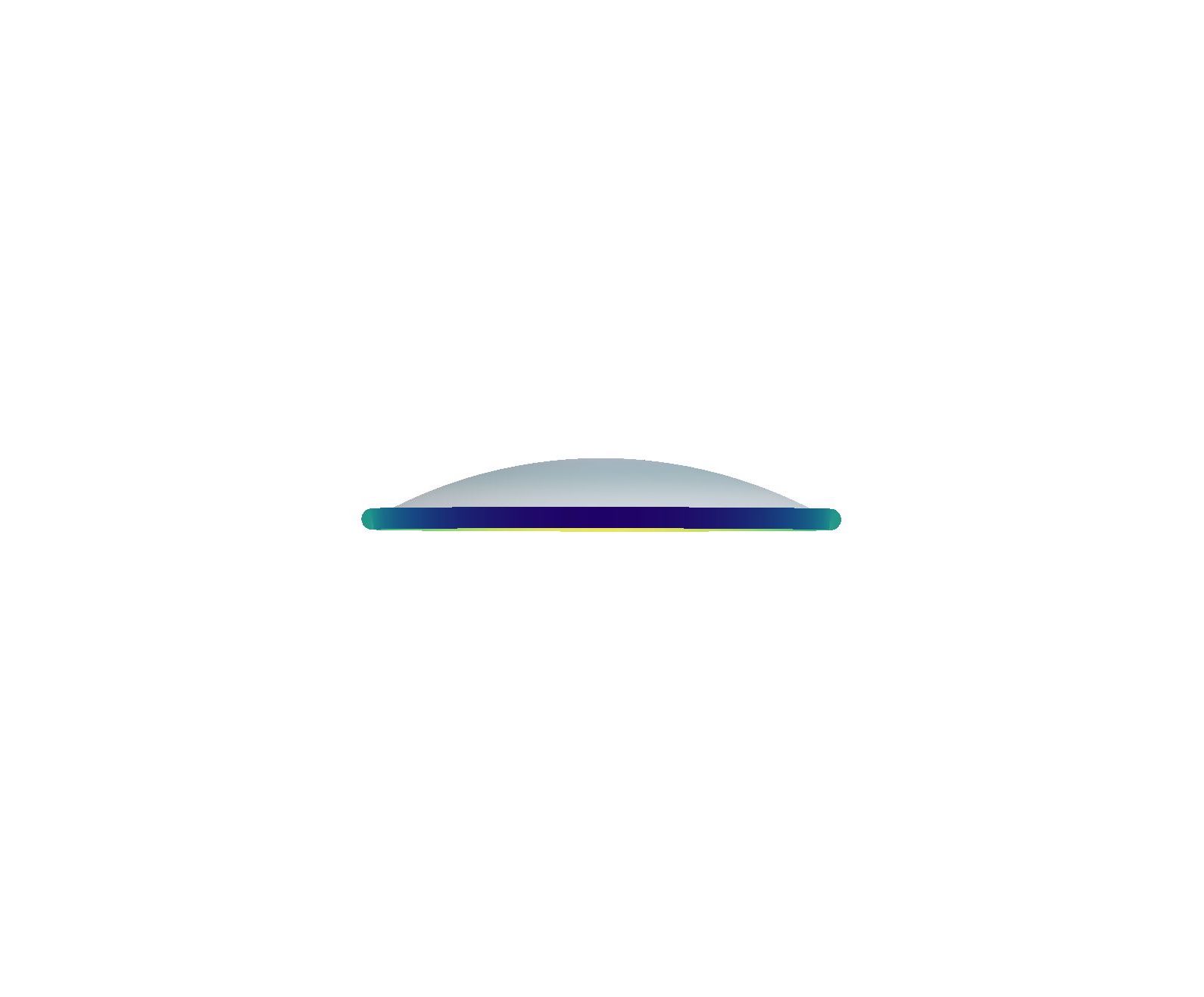}} & \quad\,\raisebox{0.19\height}{\includegraphics[width=0.28\linewidth]{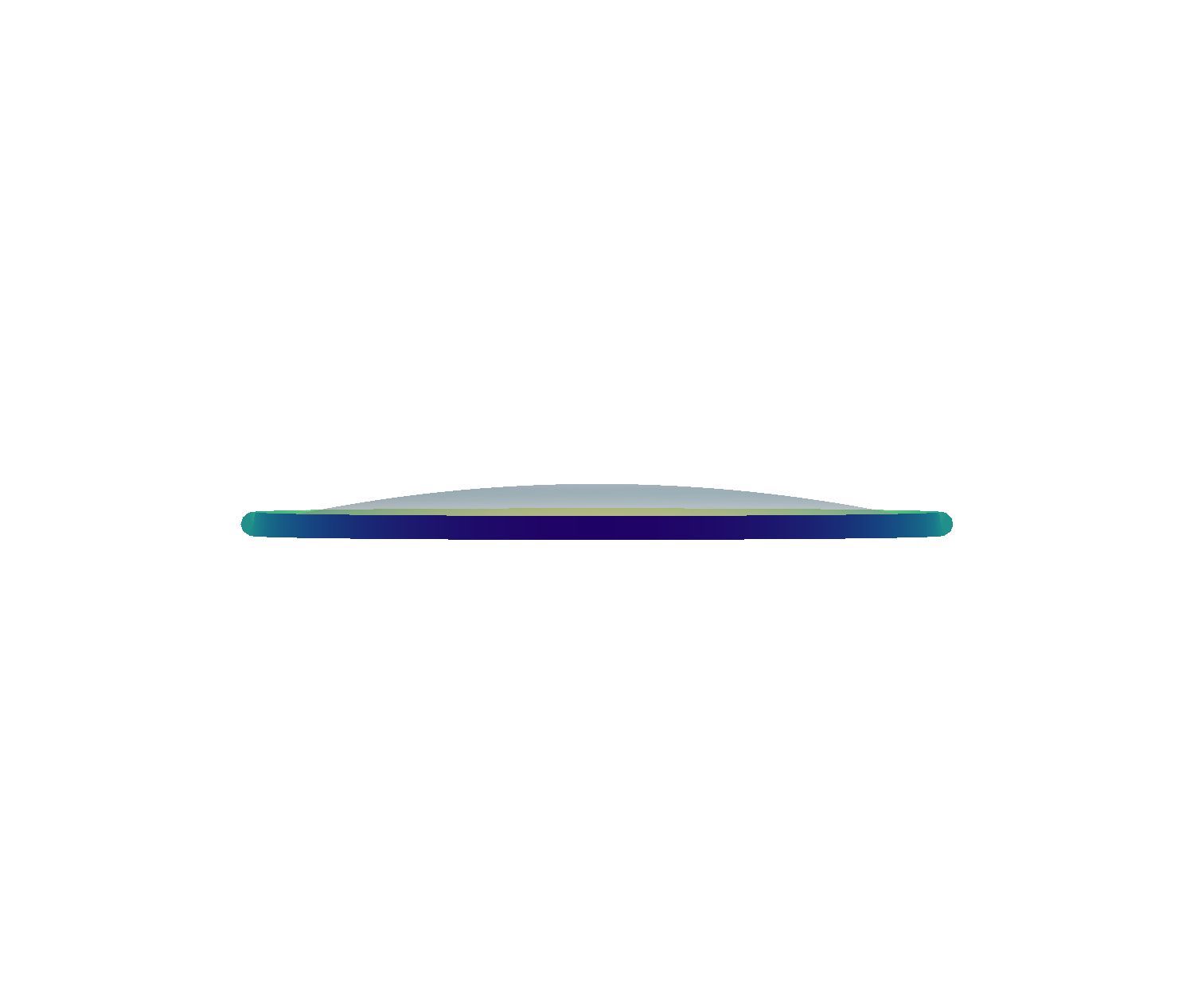}}\quad\, \\ 
\hline
\rule{0pt}{4.3cm}\raisebox{7.3\height}{$c_o=1$} & \,\,\raisebox{-0.045\height}{\includegraphics[width=0.35\linewidth]{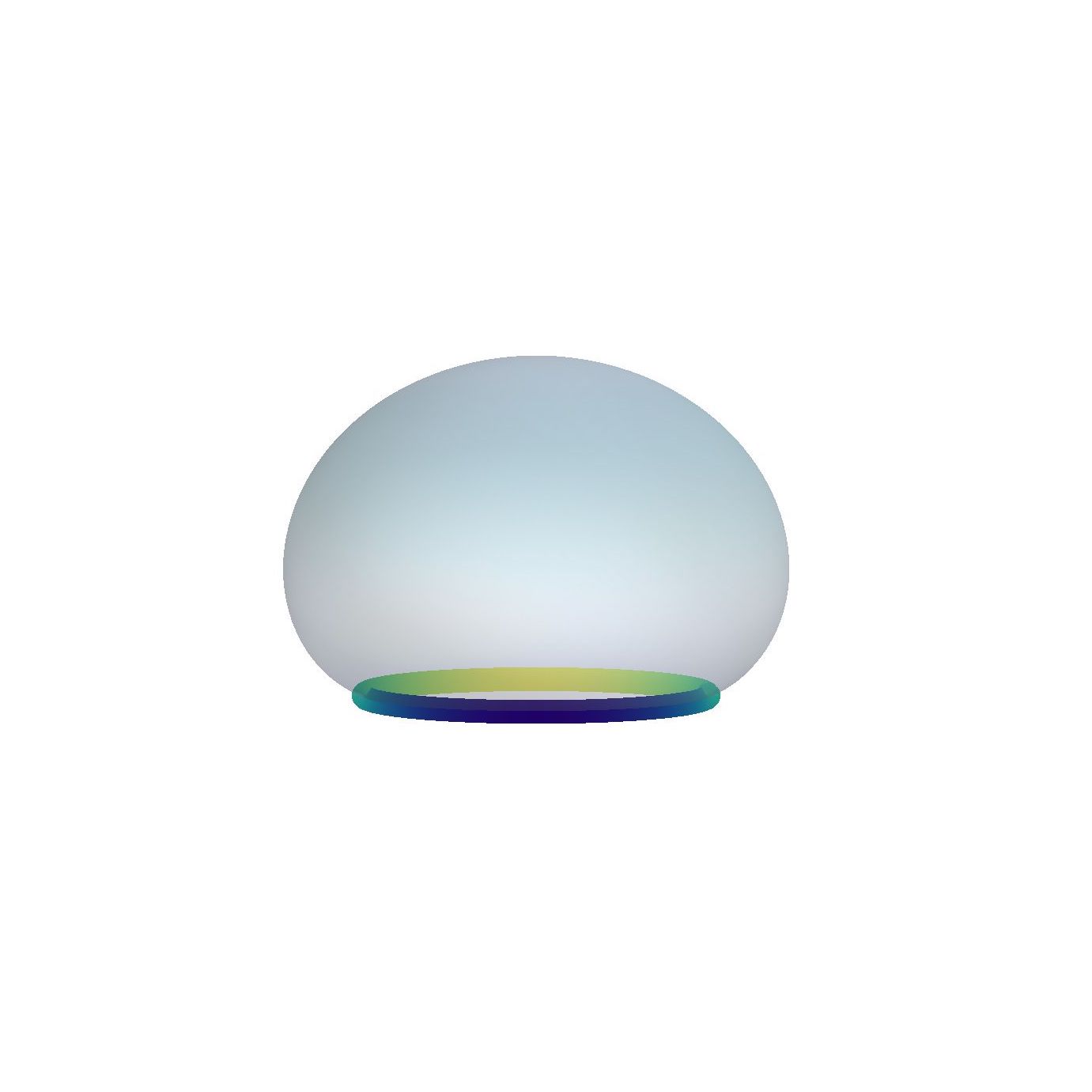}} & \raisebox{0.0\height}{\includegraphics[width=0.33\linewidth]{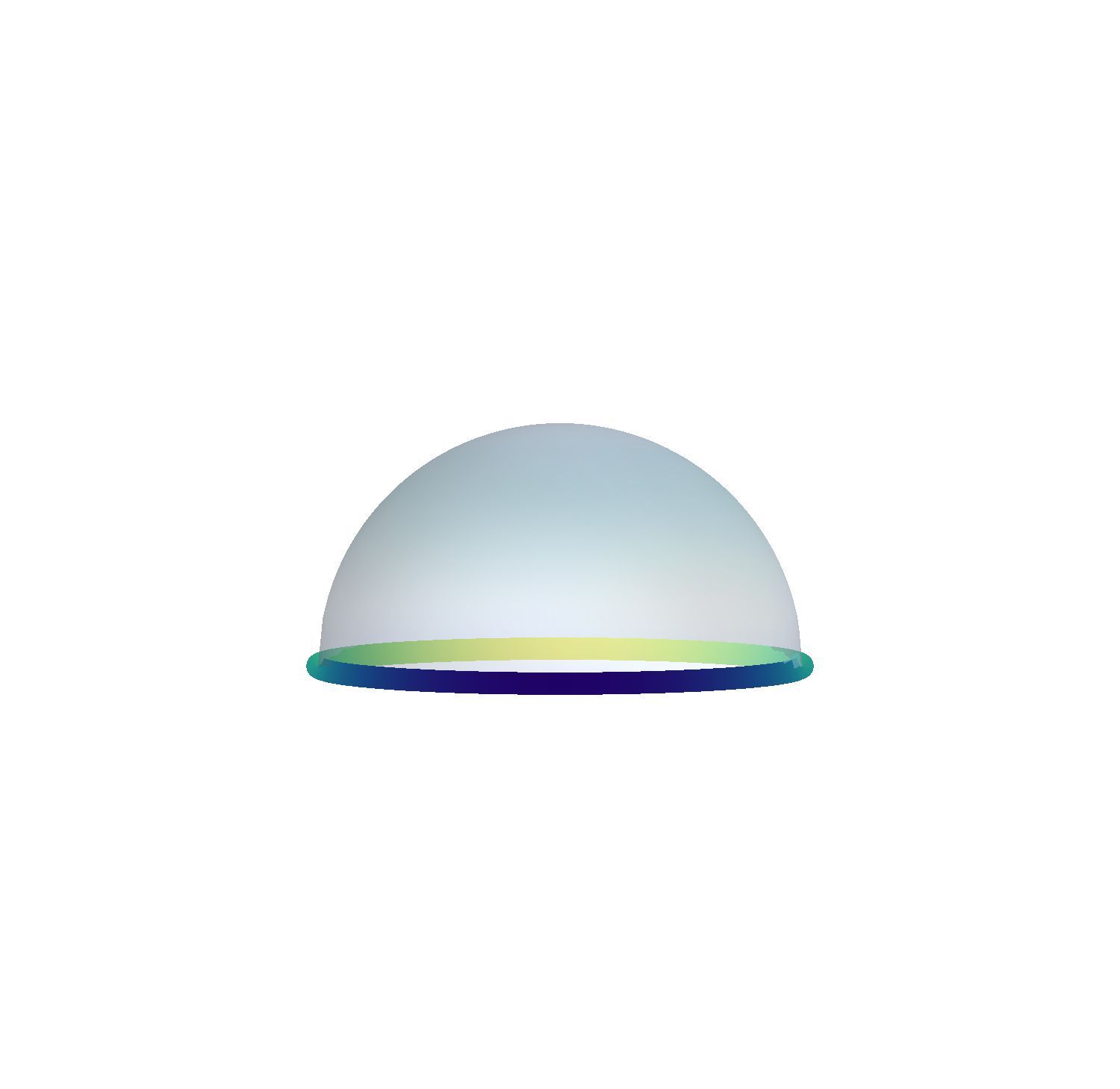}} & \quad\,\raisebox{0.19\height}{\includegraphics[width=0.28\linewidth]{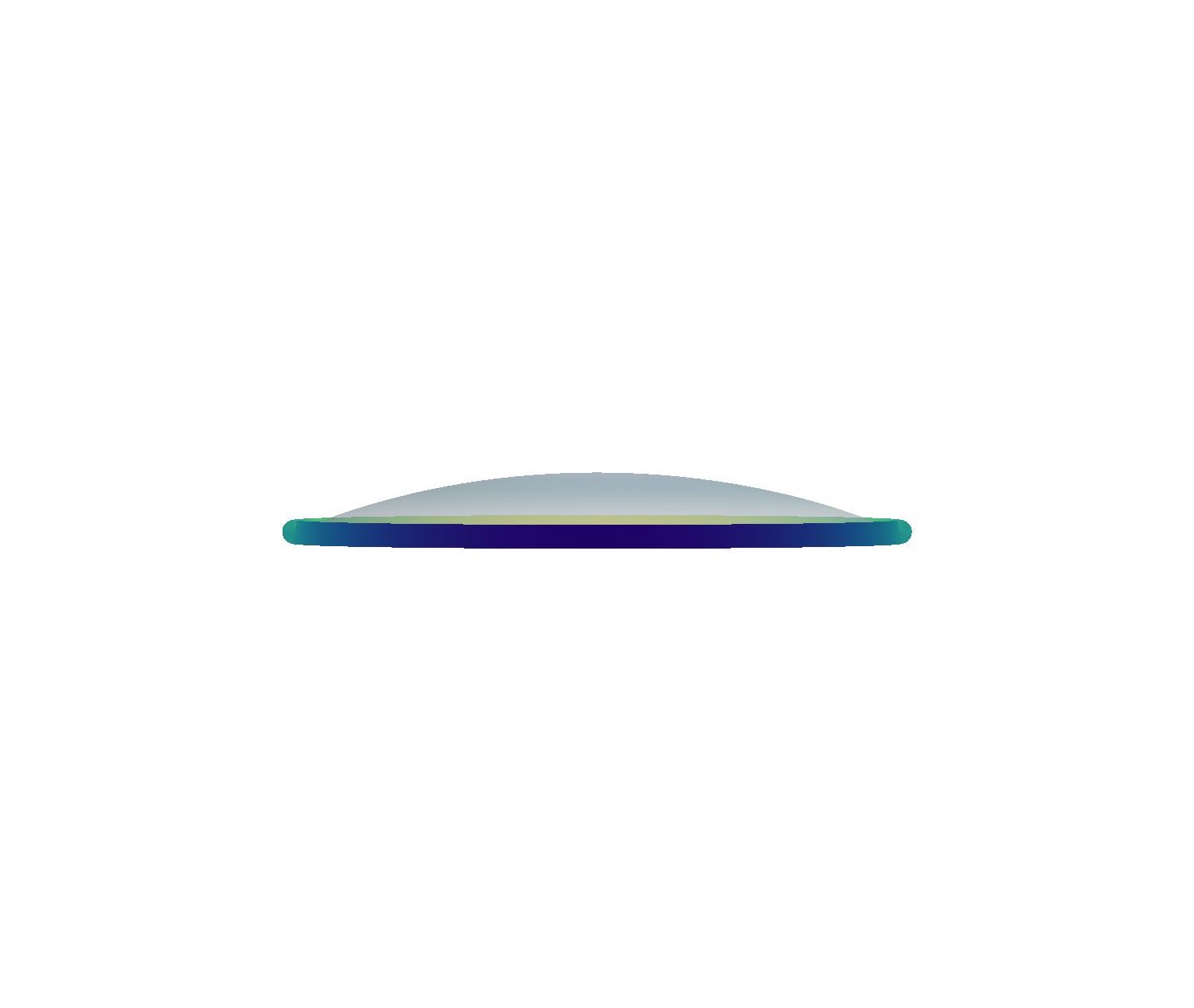}}\quad\, \\
\hline
\rule{0pt}{4.3cm}\raisebox{7.3\height}{$c_o=2$} & \raisebox{-0.015\height}{\includegraphics[width=0.35\linewidth]{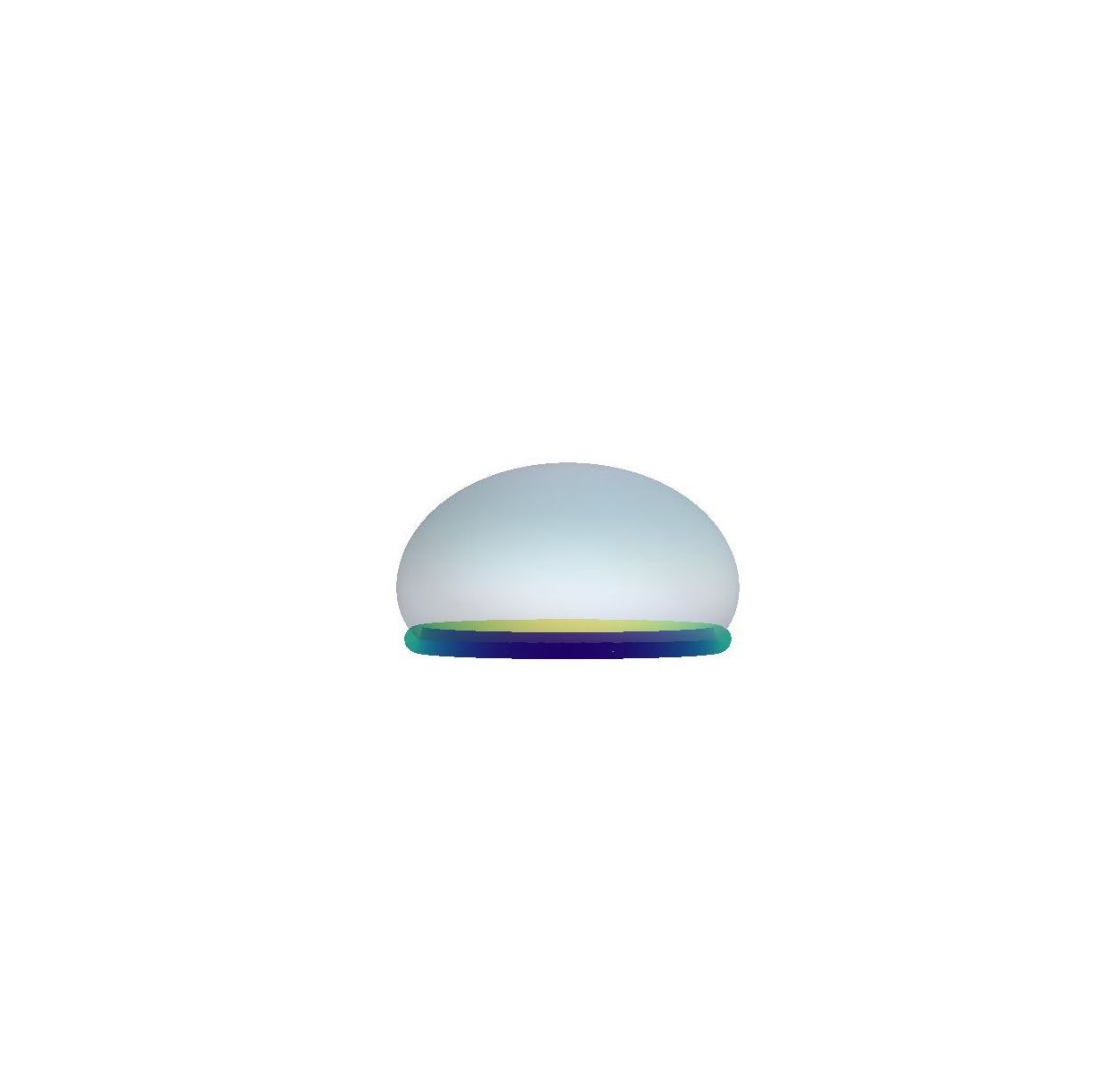}} & \raisebox{0.03\height}{\includegraphics[width=0.33\linewidth]{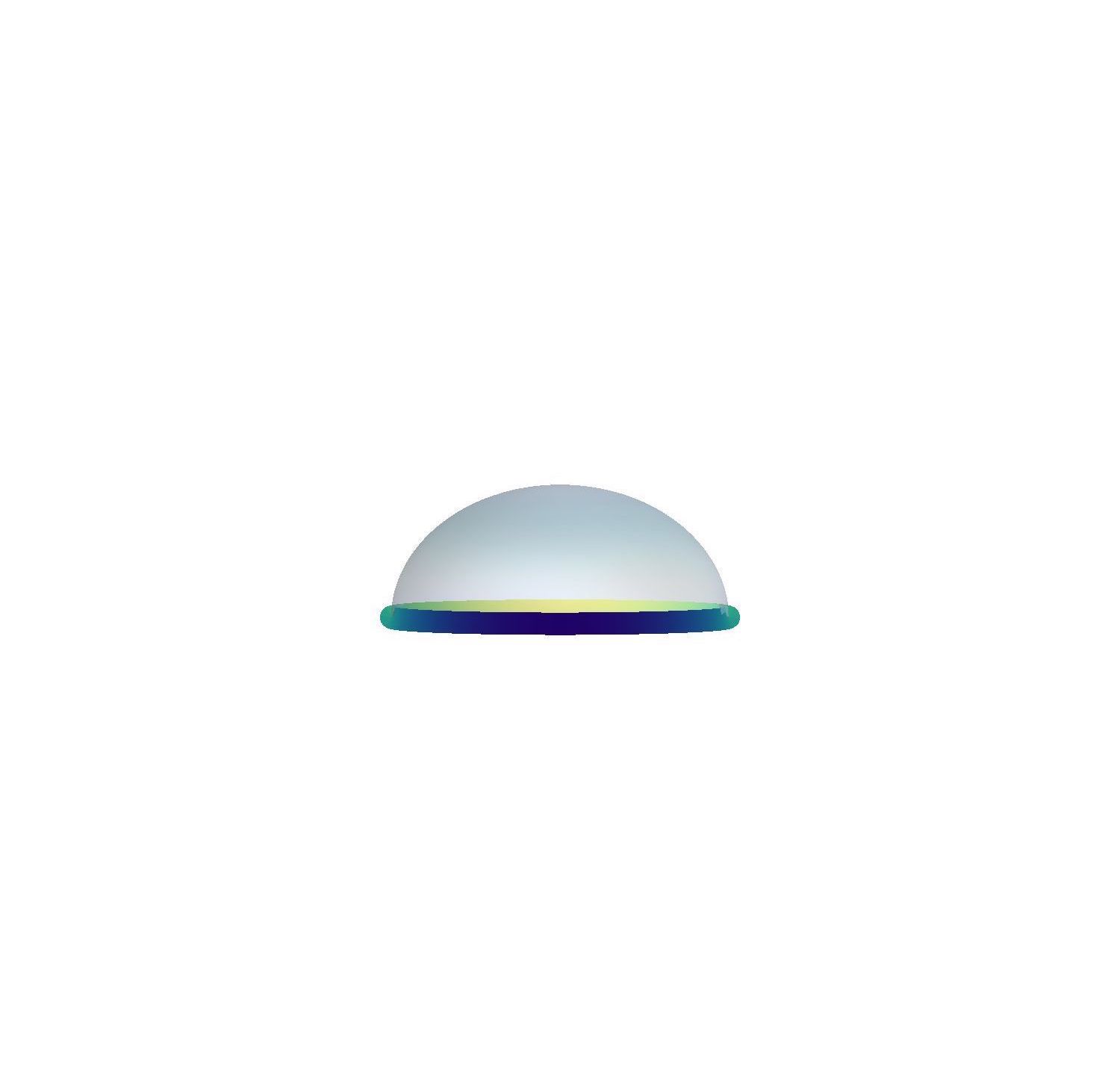}} & \quad\,\raisebox{0.19\height}{\includegraphics[width=0.28\linewidth]{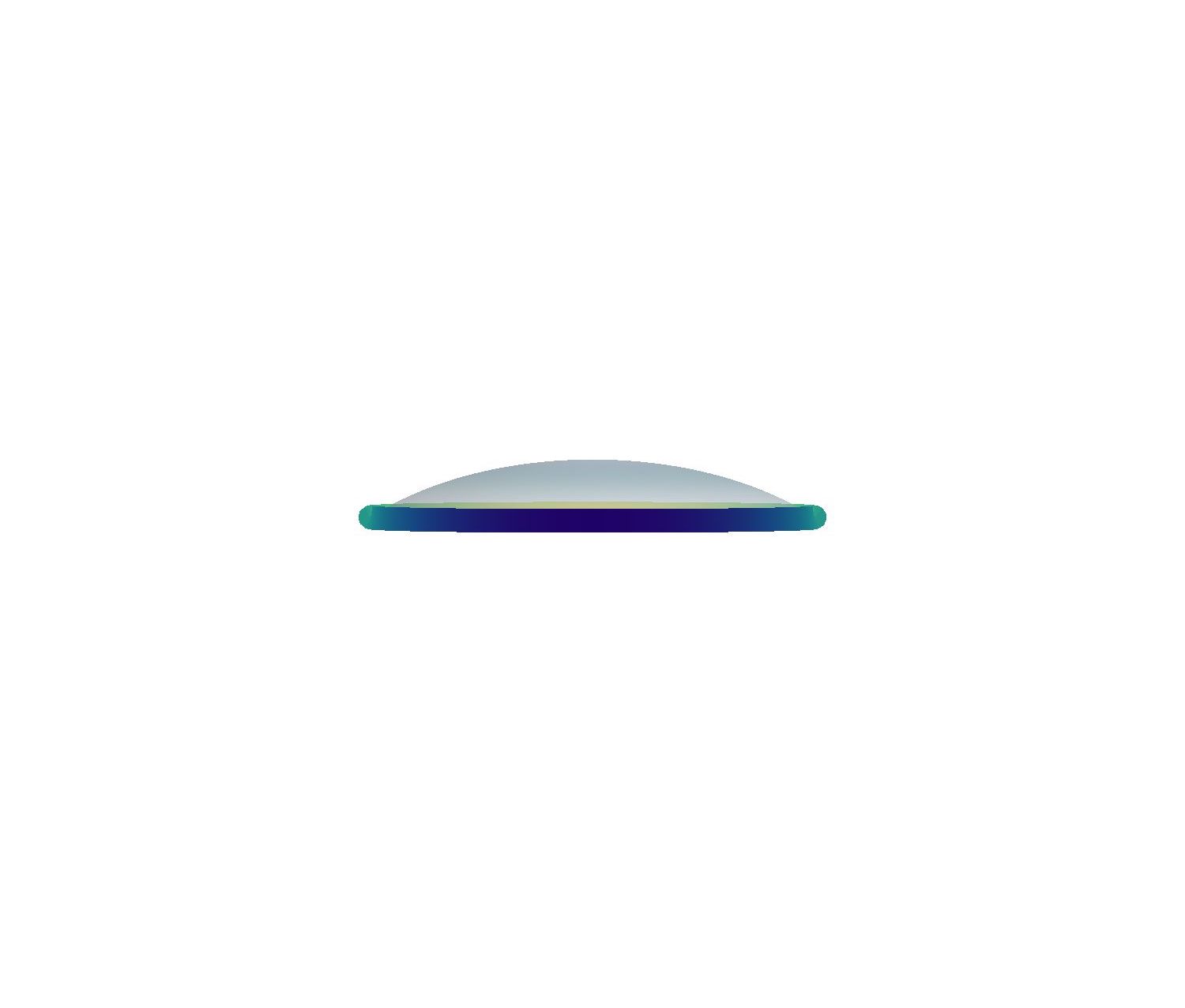}}\quad\, \\
\hline
\rule{0pt}{4.3cm}\raisebox{7.3\height}{$c_o=5$} & \raisebox{0.03\height}{\includegraphics[width=0.35\linewidth]{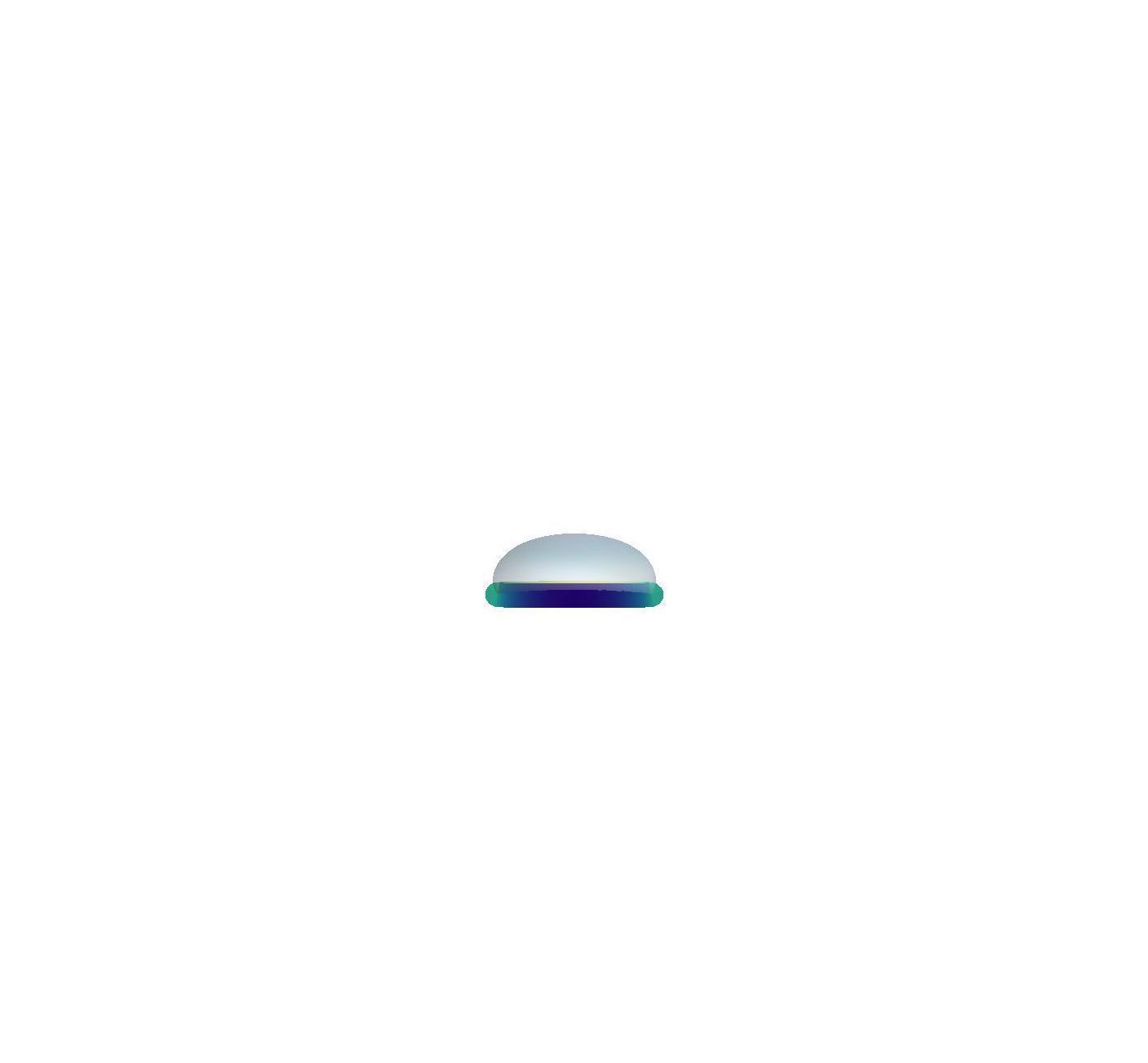}} & \raisebox{0.07\height}{\includegraphics[width=0.33\linewidth]{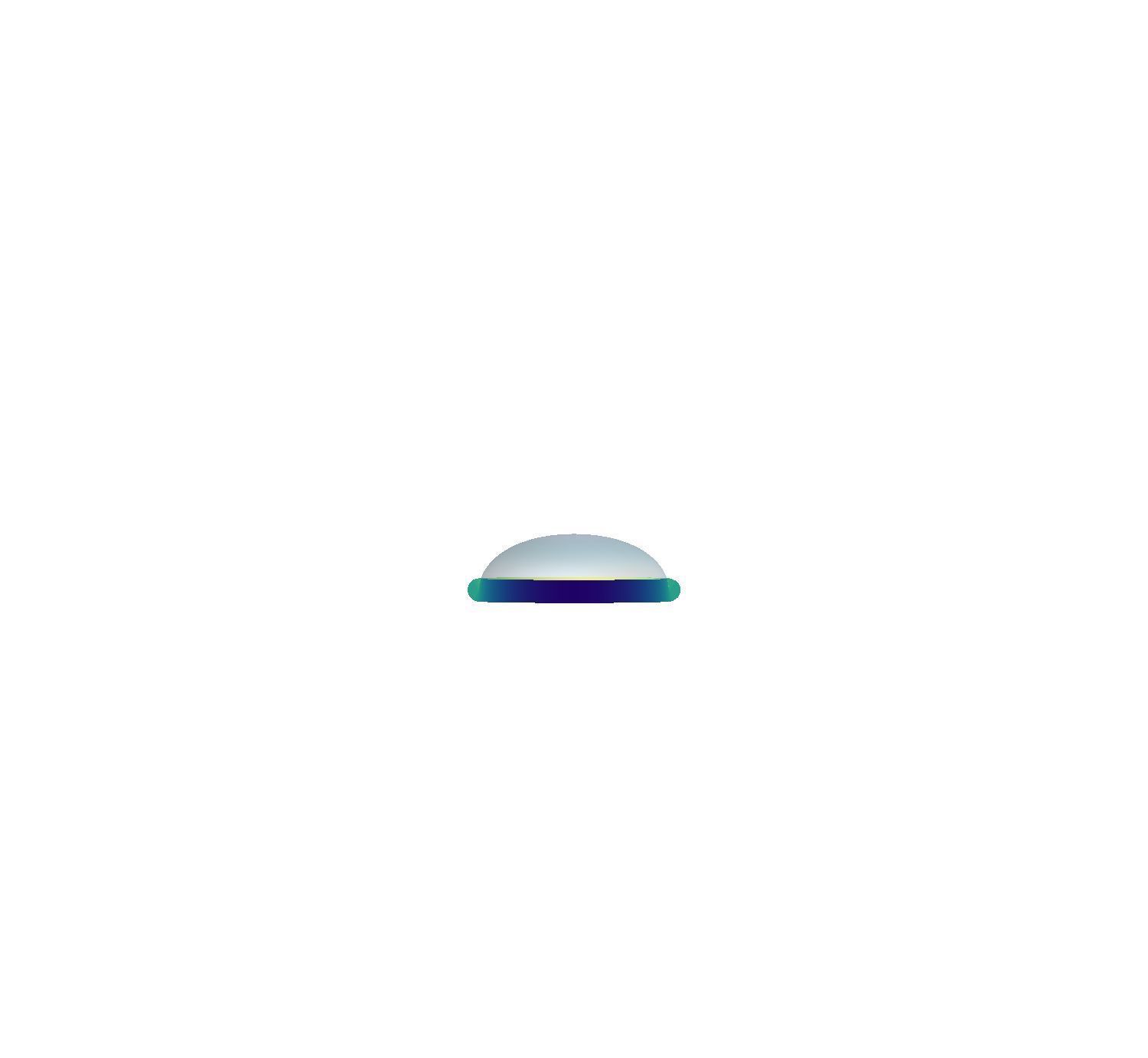}} & \quad\,\raisebox{0.19\height}{\includegraphics[width=0.28\linewidth]{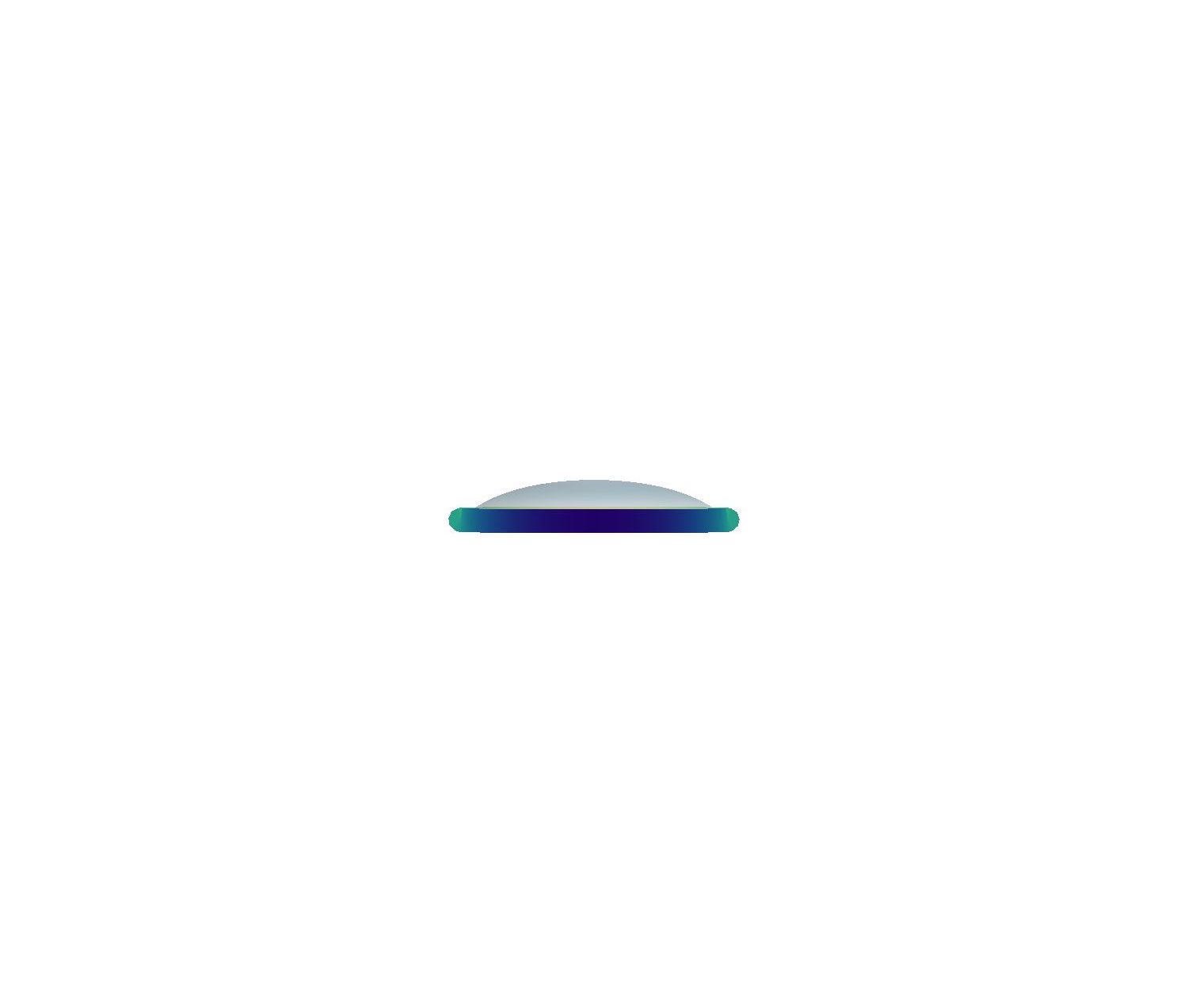}}\quad\, \\
\hline
\end{tabular}}
\end{table}
\end{center}

\bigskip

\begin{flushleft}
Bennett P{\footnotesize ALMER}\\
Department of Mathematics,
Idaho State University,
Pocatello, ID 83209,
U.S.A.\\
E-mail: palmbenn@isu.edu
\end{flushleft}

\bigskip

\begin{flushleft}
\'Alvaro P{\footnotesize \'AMPANO}\\
Department of Mathematics and Statistics, Texas Tech University, Lubbock, TX
79409, U.S.A.\\
E-mail: alvaro.pampano@ttu.edu
\end{flushleft} 


\begin{thebibliography}{777}

\bibitem{AB} M. Asgari and A. Biria, Free energy of the edge of an open lipid bilayer based on the interactions of its constituent molecules, \emph{Int. J. Nonlinear Mech.} \textbf{76} (2015), 135--143.

\bibitem{AGM} J. Arroyo, O. J. Garay and J. J. Menc\'ia, Elastic curves with constant curvature at rest in the hyperbolic plane, \emph{J. Geom. Phys.} \textbf{61} (2011), 1823--1844.

\bibitem{BMF} A. Biria, M. Maleki and E. Fried, Continuum theory for the edge of an open lipid bilayer, \emph{Adv. Appl. Mech.} \textbf{46} (2013), 1--68.

\bibitem{BR} D. H. Boal and M. Rao, Topology changes in fluid membranes, \emph{Phys. Rev. A} \textbf{46} (1992), 3037.

\bibitem{C} R. Capovilla, J. Guven and J. Santiago, Lipid membranes with an edge, \emph{Phys. Rev. E} \textbf{66} (2002), 021607.

\bibitem{DDG} K. Deckelnick, M. Doemeland, H. C. Grunau, Boundary value problems for a special Helfrich functional for surfaces of revolution: existence and asymptotic behaviour, \emph{Calc. Var. Partial Differ. Equ.} \textbf{60-1} (2021), 1--31.

\bibitem{Dierkes} U. Dierkes, Singular minimal surfaces. In: S. Hildebrandt and H. Karcher (eds), \emph{Geometric Analysis and Nonlinear Partial Differential Equations}, Springer, Berlin, 2003, pp. 177--193.

\bibitem{Euler} L. Euler, De Curvis Elasticis. In: \emph{Methodus Inveniendi Lineas Curvas Maximi Minimive Propietate Gaudentes,
Sive Solutio Problematis Isoperimetrici Lattissimo Sensu Accepti}, Additamentum 1 Ser. 1 24, Lausanne, 1744.

\bibitem{G17} T. Gibaud, C. N. Kaplan, P. Sharma, M. J. Zakhary, A. Ward, R. Oldenbourg, R. B. Meyer, R. D. Kamien, T. R. Powers, and Z. Dogic, Achiral symmetry breaking and positive Gaussian modulus lead to scalloped colloidal membranes, \emph{Proc. Natl. Acad. Sci. USA} \textbf{114-17} (2017), 3376--3384.

\bibitem{Helfrich} W. Helfrich, Elastic properties of lipid bilayers: theory and possible experiments. \emph{Zeit. Naturfor. C} \textbf{28} (1973), 693--703.

\bibitem{KP} M. Koiso and B. Palmer, Geometry and stability of bubbles with gravity, \emph{Indiana Univ. Math. J.} \textbf{54} (2005), 65--98.

\bibitem{KP2} M. Koiso and B. Palmer, On a variational problem for soap films with gravity and partially free boundary, \emph{J. Math. Soc. Japan} \textbf{57} (2005), 333--355.

\bibitem{Levien} R. Levien, The elastica: a mathematical history, Technical Report No. UCB/EECS-2008-103, University of Berkeley.

\bibitem{Rafa} R. L\'opez, Invariant singular minimal surfaces, \emph{Ann. Glob. Anal. Geom.} \textbf{53} (2018), 521--541.

\bibitem{Rafa2} R. L\'opez, Symmetry of stationary hypersurfaces in hyperbolic space, \emph{Geom. Dedicata} \textbf{119} (2006), 35--47.

\bibitem{L} A. E. H. Love, \emph{A Treatise on the Mathematical Theory of Elasticity}, Dover Publications 4th Ed., New York, 1944.

\bibitem{MF} M. Maleki and E. Fried, Stability of discoidal high-density lipoprotein particles, \emph{Soft Matter} \textbf{9-42} (2013), 9991--9998.

\bibitem{MS} A. Mondino and C. Scharrer, Existence and regularity of spheres minimising the Canham-Helfrich energy, \emph{Arch. Ration. Mech. Anal.} \textbf{236-3} (2020), 1455--1485.

\bibitem{Morrey} C. B. Morrey Jr., \emph{Multiple Integrals in the Calculus of Variations}, Springer, New York, 2009.

\bibitem{NOOY} H. Naito, M. Okuda and Z. C. Ou-Yang, Polygonal shape transformation of a circular biconcave vesicle induced by osmotic pressure, \emph{Phys. Rev. E} \textbf{54} (1996), 2816--2826.

\bibitem{Palmer} B. Palmer, The conformal Gauss map and the stability of Willmore surfaces, \emph{Ann. Glob. Anal. Geom.} \textbf{9} (1991), 305--317.

\bibitem{P} B. Palmer, Uniqueness theorems for Willmore surfaces with fixed and free boundaries, \emph{Indiana Univ. Math. J.} \textbf{49-4} (2000), 1581--1601.

\bibitem{PP2} B. Palmer and A. P\'ampano, Minimizing configurations for elastic surface energies with elastic boundaries, \emph{J. Nonlinear Sci.} \textbf{31} (2021).

\bibitem{RL} B. R\'ozycki and R. Lipowsky, Spontaneous curvature of bilayer membranes from molecular simulations: asymmetric lipid densities and asymmetric adsorption, \emph{J. Chem. Phys.} \textbf{142-5} (2015), 054101.

\bibitem{CZT}  Z. C. Tu, Compatibility between shape equation and boundary conditions of lipid membranes with free edges, \emph{J. Chem. Phys.} \textbf{132-8} (2010), 084111.

\bibitem{Tu2} Z. C. Tu, Geometry of membranes, \emph{J. Geom. Symmetry Phys.} \textbf{24} (2011), 45--75.

\bibitem{TOY} Z. C. Tu and Z. C. Ou-Yang, A geometric theory on the elasticity of bio-membranes, \emph{J. Phys. A: Math. Gen.} \textbf{37} (2004), 11407--11429.

\bibitem{TOY03} Z. C. Tu and Z. C. Ou-Yang, Lipid membranes with free edges, \emph{Phys. Rev. E} \textbf{68} (2003), 061915.

\bibitem{TOY2} Z. C. Tu and Z. C. Ou-Yang, Recent theoretical advances in elasticity of membranes following Helfrich's spontaneous curvature model, \emph{Adv. Colloid Interface Sci.} \textbf{208} (2014), 66--75.

\bibitem{W} N. Walani, J. Torres, and A. Agrawal, Anisotropic spontaneous curvatures in lipid membranes, \emph{Phys. Rev. E} \textbf{89-6} (2014), 062715.

\bibitem{Z} X. Zhou, An integral case of the axisymmetric shape equation of open vesicles with free edges, \emph{Int. J. Nonlinear Mech.} \textbf{106} (2018), 25--28.

\end{thebibliography}
\end{document}